\documentclass[12pt,twoside]{article}
\usepackage{a4wide}
\usepackage{amsmath,mathrsfs,amssymb,amsthm,amscd}
\usepackage[]{fontenc}
\usepackage[all]{xy} 
\numberwithin{equation}{section}

\theoremstyle{plain}
\newtheorem{theorem}{Theorem}[section]
\newtheorem{corollary}[theorem]{Corollary}
\newtheorem{proposition}[theorem]{Proposition}
\newtheorem{lemma}[theorem]{Lemma}

\theoremstyle{definition}
\newtheorem{definition}[theorem]{Definition}

\newtheorem{notation}[theorem]{Notation}

\theoremstyle{remark}
\newtheorem{remark}[theorem]{Remark}

\newcommand{\N}{\mathbb{N}}
\newcommand{\Z}{\mathbb{Z}}
\newcommand{\Q}{\mathbb{Q}}
\newcommand{\R}{\mathbb{R}}
\newcommand{\C}{\mathbb{C}}
\newcommand{\D}{\mathbf{D}}

\newcommand{\E}{\mathbb{E}}
\renewcommand{\H}{\mathbb{H}}
\renewcommand{\P}{\mathbb{P}}

\newcommand{\ZZ}{{\mathbb Z}}
\newcommand{\RR}{{\mathbb R}}

\newcommand{\CC}{{\mathbb C}}
\newcommand{\QQ}{{\mathbb Q}}

\newcommand{\e}{{\mathbf e}}

\newcommand{\calH}{\mathcal{H}}

\newcommand{\calM}{\mathcal{M}}

\newcommand{\calO}{\mathcal{O}}

\newcommand{\calT}{\mathcal{T}}

\newcommand{\calX}{\mathcal{X}}
\newcommand{\calY}{\mathcal{Y}}

\newcommand{\fraka}{\mathfrak a}
\newcommand{\frakb}{\mathfrak b}
\newcommand{\frakc}{\mathfrak c}
\newcommand{\frakd}{\mathfrak d}

\newcommand{\frakp}{\mathfrak p}

\newcommand{\frakt}{\mathfrak t}

\newcommand{\frakl}{\mathfrak l}
\newcommand{\frakg}{\mathfrak g}

\newcommand{\zxz}[4]{\begin{pmatrix} #1 & #2 \\ #3 & #4 \end{pmatrix}}
\newcommand{\abcd}{\zxz{a}{b}{c}{d}}
\newcommand{\kzxz}[4]{\left(\begin{smallmatrix} #1 & #2 \\ #3 & #4\end{smallmatrix}\right) }
\newcommand{\kabcd}{\kzxz{a}{b}{c}{d}}

\newcommand{\Gr}{\operatorname{Gr}}
\newcommand{\diva}{\operatorname{\widehat{div}}}
\newcommand{\cha}{\operatorname{\widehat{CH}}}
\newcommand{\CH}{\operatorname{CH}}
\newcommand{\za}{\operatorname{\widehat{Z}}}
\newcommand{\rata}{\operatorname{\widehat{Rat}}}

\newcommand{\dd}{\operatorname{d}}
\newcommand{\ddc}{\operatorname{dd}^c}
\newcommand{\Img}{\operatorname{Im}}

\newcommand{\cl}{\operatorname{cl}}

\newcommand{\Hom}{\operatorname{Hom}}
\newcommand{\End}{\operatorname{End}}

\newcommand{\Pic}{\operatorname{Pic}}

\newcommand{\Coker}{\operatorname{Coker}}
\newcommand{\Spec}{\operatorname{Spec}}
\newcommand{\Proj}{\operatorname{Proj}}
\newcommand{\Sch}{\operatorname{Sch}}
\newcommand{\supp}{\operatorname{supp}}

\newcommand{\amap}{\operatorname{a}}

\newcommand{\pr}{\operatorname{pr}}
\newcommand{\fin}{\operatorname{fin}}
\newcommand{\fh}{\operatorname{ht}}
\newcommand{\sym}{\text{\rm sym}}
\newcommand{\pre}{\text{\rm pre}}
\newcommand{\Pet}{\text{\rm Pet}}
\newcommand{\ohh}{\mathbb{H}}
\newcommand{\Ga}{{\Gamma}}

\newcommand{\DD}{{D}}

\newcommand{\eps}{\varepsilon}
\newcommand{\bs}{\backslash}
\newcommand{\norm}{\operatorname{N}}
\newcommand{\vol}{\operatorname{vol}}
\newcommand{\tr}{\operatorname{tr}}

\newcommand{\Sl}{\operatorname{SL}}
\newcommand{\Gl}{\operatorname{GL}}
\newcommand{\Orth}{\operatorname{O}}

\newcommand{\Gal}{\operatorname{Gal}}

\newcommand{\Div}{\operatorname{Div}}
\renewcommand{\div}{\operatorname{div}}

\newcommand{\dv}{\operatorname{div}}
\newcommand{\dega}{\operatorname{\widehat{deg}}}
\newcommand{\pica}{\operatorname{\widehat{Pic}}}
\newcommand{\cca}{\operatorname{\widehat{c}}_1}
\newcommand{\cc}{\operatorname{c}}
\newcommand{\XX}{{\mathscr{X}}}
\newcommand{\OK}{{\mathcal{O}_K}}
\newcommand{\GK}{{\Gamma_K}}

\long\def\symbolfootnote[#1]#2{\begingroup%
\def\thefootnote{\fnsymbol{footnote}}\footnote[#1]{#2}\endgroup} 

\def\?{\ ???\ \immediate\write16{}%
\immediate\write16{Warning: There was still a question mark . . . }%
\immediate\write16{}}

\title%[Borcherds products and arithmetic intersection theory]
{Borcherds products and arithmetic intersection theory on Hilbert
modular surfaces}
\author{Jan H.~Bruinier, Jose I.~Burgos Gil, and Ulf K\"uhn
}
\date{October 13, 2003}

\begin{document}

\maketitle

\symbolfootnote[0]{
Revision: September 9, 2004.
The first author was temporarily supported by a Heisenberg-Fellowship
of the DFG.  The second author was partially supported by DGI Grant
BFM2000-0799-c02-01.
}

\begin{abstract} 
We prove an arithmetic version of a theorem of Hirzebruch and Zagier
 saying that Hirzebruch-Zagier divisors on a Hilbert modular surface are the coefficients of an elliptic modular form of weight two.
Moreover, we determine the arithmetic self-intersection number of
the line bundle of modular forms equipped with its Petersson metric 
on a regular model of a Hilbert modular surface, and study Faltings
heights of arithmetic Hirzebruch-Zagier divisors. 
\end{abstract}

\tableofcontents
% h-intro3.tex

\section*{Introduction}
\addcontentsline{toc}{section}{Introduction}

%\item[intro] 
It is of special interest in Arakelov theory to determine intrinsic
arithmetic intersection numbers of varieties defined over number
fields and Faltings heights of its subvarieties.  Here we study
Hilbert modular surfaces associated to a real quadratic field $K$.
Since these Shimura varieties are non-compact, we consider 
suitable toroidal compactifications. As the natural metrics of
automorphic vector bundles on such varieties have singularities along
the boundary (see e.g.~\cite{Mu}, \cite{BKK2}), we work with the extended
arithmetic Chow rings constructed in \cite{BKK}.
  
  In their celebrated paper \cite{HZ} Hirzebruch and Zagier proved
  that the generating series for the cohomology classes of 
  Hirzebruch-Zagier divisors is a holomorphic modular form of
  weight $2$.  A different proof was recently obtained by Borcherds \cite{Bo3}.
  Moreover, it was shown by Franke and Hausmann that 
the product of this generating series with the first Chern class of the line
  bundle of modular forms, i.e., the generating series for the
  hyperbolic volumes, equals a particular holomorphic Eisenstein series
  $E(\tau)$ of weight $2$ \cite{Fra}, \cite{Hau}. 

In the present paper we define certain arithmetic Hirzebruch-Zagier divisors and
show that the generating series of their classes in the arithmetic
Chow ring of a suitable regular model of the Hilbert modular
surface is a holomorphic modular form (of the same level, weight, and character as in the case of Hirzebruch and Zagier). 
In addition we show that the
product of this generating series with the square of the first
arithmetic Chern class of the line bundle of modular forms equipped
with its Petersson metric is equal to a multiple of $E(\tau)$.
The factor of proportionality is the arithmetic self intersection
number of the line bundle of modular forms.

Since
Hilbert modular surfaces can also be viewed as  Shimura varieties associated to the orthogonal group of a rational quadratic space of signature $(2,2)$, these results are  related to the program described by Kudla 
in \cite{Ku}, \cite{Ku2}, \cite{Ku:ICM}, \cite{Ku:MSRI};
for a discussion of results in that direction see also to the references therein. 
However, notice that there are several slight technical differences. For instance in the present case the Shimura variety has to be compactified, the arithmetic Hirzebruch-Zagier divisors contain boundary components, and
we work with different Green functions.

We show that the arithmetic self intersection number of
the line bundle of modular forms is essentially given by the
logarithmic derivative at $s=-1$ of the Dedekind zeta function $\zeta_K(s)$ of $K$.  It equals the expression conjectured in \cite{Kue1}. We
refer to Theorem B below for the precise statement.

We also determine the Faltings heights of those Hirzebruch-Zagier divisors
which are disjoint to the boundary. It is well known that their
normalizations are isomorphic to compact Shimura curves associated with
quaternion algebras. This result will be used in a subsequent paper
\cite{KuKue} to determine the arithmetic self intersection number of
the Hodge bundle equipped with the hyperbolic metric on Shimura
curves.

Our formulas provide further evidence for the conjecture of Kramer, based on results obtained in \cite{kramer} and \cite{Kue1}, that
the arithmetic volume is essentially the derivative of the zeta value
for the volume of the fundamental domain, the
conjecture of Kudla on the constant term of
the derivative of certain Eisenstein series \cite{Ku2}, 
\cite{Ku:ICM}, \cite{Ku:MSRI},
\cite{KRY2}, and the conjecture of Maillot-Roessler
on special values of logarithmic derivatives of Artin  $L$-functions 
\cite{MaRo}, \cite{MaRo:annals}.

Our approach 
%to the results sketched above
requires various results on regular
models for Hilbert modular surfaces \cite{Ra}, \cite{DePa}, \cite{Pa}, 
as well as an extensive use of the theory of Borcherds products \cite{Bo2},
\cite{Bo3}, \cite{Br1}, \cite{Br2}.  Another central point is
the $q$-expansion principle \cite{Ra}, \cite{Ch}, which relates
analysis to geometry.

\medskip

We now describe the content of this paper in more detail.

\medskip
  
%\item[chap1] 
Since the natural metrics on automorphic line bundles have
singularities along the boundary, we cannot use arithmetic
intersection theory as presented in \cite{SABK} and have to work with
the extension developed in \cite{BKK}.  In particular, our arithmetic
cycles will be classes in the generalized arithmetic Chow ring
$\widehat{\textrm{CH}}^*(\XX,\mathcal{D}_{\pre})$, in which the Green
objects satisfy besides the usual logarithmic additional log-log
growth conditions. We recall some of its basic properties in Section
1.  Choosing representatives, the arithmetic intersection product of
two arithmetic cycles can be split into the sum of a geometric
contribution as in \cite{SABK} and an integral over a star product of
Green objects. The formulas for the latter quantity generalize those
of \cite{SABK} and may contain additional boundary terms.

We concentrate on the analytic
aspects in Sections $2$ and $3$. The geometric contribution is
considered in Sections $5$ and $6$.

\medskip

%\item[chap2] 
We begin Section 2 by recalling some facts on the analytic theory of
Hilbert modular surfaces.  Throughout the paper the discriminant $D$
of $K$ is assumed to be prime.  Let $\Ga$ be a subgroup of finite
index of the Hilbert modular group $\GK= \Sl_2(\OK)$, where $\OK$
denotes the ring of integers in $K$.  We consider an (at the outset
arbitrary) desingularization $\widetilde{X}(\Ga)$ of the Baily-Borel
compactification $X(\Ga)$ of $\Ga \bs \H^2$, where $\H$ is the upper
complex half plane.  Here the curve lemma due to Freitag is a main
technical tool to obtain adequate local descriptions.  We introduce
the line bundle of modular forms $\mathcal{M}_k(\Ga)$ of weight $k$
with its Petersson metric $\|\cdot\|_{\Pet}$, which is singular along
the normal crossing divisor formed by the exceptional curves of the
desingularization.  We show that the Petersson metric is a pre-log
singular hermitian metric in the sense of \cite{BKK}.

For any positive integer $m$, there is an algebraic divisor $T(m)$ on
$X(\Gamma)$, called the Hirzebruch-Zagier divisor of discriminant $m$.
These divisors play a central role in the study of the geometry of
Hilbert modular surfaces (see e.g.~\cite{Ge}).  We consider the Green function
$\Phi_m(z_1,z_2,s)$ associated to $T(m)$ introduced in \cite{Br1}.  A
new point here is the explicit description of the constant term, given
in Theorem \ref{satzformel2}, which involves $\zeta_K(s)$ and a
certain generalized divisor function $\sigma_m(s)$ (see
\eqref{def:sig}).  We then define the automorphic Green function $G_m$
for $T(m)$ to be the constant term in the Laurent expansion at $s=1$
of $\Phi_m(z_1,z_2,s)$ minus a natural normalizing constant, which is
needed to obtain compatibility with the theory of Borcherds products
as used in Section 4.  By means of the curve lemma we show that these
Green functions define pre-log-log Green objects
 in the sense of \cite{BKK}.

\medskip
  
%\item[chap3] 

In Section 3 we specialize the general formula for the star product to
the Green objects $\frakg(m_j)= (-2\partial \bar \partial
G_{m_j},G_{m_j})$ associated to suitable triples of Hirzebruch-Zagier
divisors $T(m_j)$ (see Theorem \ref{thm:123star}).  In this case,
using the curve lemma, it can be shown that certain boundary terms
vanish. As a consequence, the resulting formula only depends on the
Baily-Borel compactification $X(\Gamma)$ and is independent of the
choices made in the desingularization.  An analogous formula holds for
the star product associated to a Hirzebruch-Zagier divisor $T(m)$ and
the divisors of two Hilbert modular forms.  In that way we 
reduce the analytical contributions to the arithmetic intersection
numbers in question to integrals of Green functions $G_m$ over
$X(\Ga)$ and star products on Hirzebruch-Zagier divisors $T(m)$.
  
  The analysis of Section 2 allows us to determine the
  integral of $G_m$ in terms of the logarithmic derivatives
  of $\zeta_K(s)$ and $\sigma_m(s)$ at $s=-1$ (Corollary \ref{cor:greenint}).
  If $p$ is a prime that splits in $K$, then $T(p)$ is birational to the
  modular curve $X_0(p)$. We determine the star products on these
  $T(p)$ using the results of \cite{Kue2} in Theorem \ref{thm:modulstar}. 
 
\medskip

%\item[chap4] 
The proofs of our main results rely on the theory of Borcherds
products in a vital way. Therefore in Section 4 we recall some of
their basic properties with an emphasis on the construction given in
\cite{Br1}.  Borcherds products on Hilbert modular surfaces are
particular meromorphic Hilbert modular forms that arise as regularized
theta lifts of certain weakly holomorphic elliptic modular forms of
weight $0$. They enjoy striking arithmetic properties. For instance, a
sufficiently large power of any holomorphic Borcherds product has
coprime integral Fourier coefficients. Hence, by the $q$-expansion
principle, it defines a section of the line bundle of modular forms
over $\Z$.  Moreover, the divisors of Borcherds products are explicit
linear combinations of Hirzebruch Zagier divisors, dictated by the
poles of the input modular form.  By \cite{Br1}, the logarithm of
their Petersson norm is precisely given by a linear combination of the
Green functions $G_m$.
  
By a careful analysis of the obstruction space for the existence of
Borcherds products we prove that for any given linear combination of
Hirzebruch-Zagier divisors $C$ there are infinitely many Borcherds
products $F_1$, $F_2$, such that $C\cap \dv(F_1)\cap
\dv(F_2)=\emptyset$, and such that further technical conditions are
satisfied (Theorem \ref{densitytriple}).  This will give an ample
supply of sections of the line bundle of modular forms for which the
associated star products can be calculated.  In the rest of this
section we essentially show that the subspace of
$\Pic(X(\Gamma_K))\otimes_\Z\Q$ spanned by all Hirzebruch-Zagier
divisors, is already generated by the $T(p)$, where $p$ is a prime
which splits in $K$ (Theorem \ref{densityshim}). This fact can 
be viewed as an explicit moving lemma for Hirzebruch-Zagier divisors (cf.~Theorem \ref{thm:chowfinite}), and is crucial
in the proof of Theorem B.

\medskip
  
%\item[chap5] 
In Section 5 we recall the arithmetic theory of Hilbert modular
surfaces. Unfortunately, there are currently no references for
projective regular models defined over $\Spec \ZZ$. Therefore we will
work with regular models over the subring $\ZZ[\zeta_N,1/N]$ of the
$N$-th cyclotomic field $\QQ(\zeta_N)$.  More precisely, we consider
toroidal compactifications $ \widetilde{\calH}(N)$ of the moduli
scheme associated with the Hilbert modular variety for the principal
congruence subgroup $\GK(N)$ of arbitrary level $N\geq 3$.  We define
$\mathcal{T}_N(m) \subseteq \widetilde{\calH}(N)$ as the Zariski
closure of the Hirzebruch-Zagier divisor $T(m)$ on the generic fiber.
This definition is well-behaved with respect to pull-backs and is
compatible with the theory of Borcherds products (see Proposition \ref{prop:moduli-divisor}).  If $p$ is a split
prime, there exists a natural modular morphism from  the compactified moduli space of elliptic curves with a subgroup of order $p$
to the minimal compactification of $\calH(1)$ whose image is essentially $\mathcal{T}_1(p)$ (cf.~Proposition \ref{prop:tpmin}).  This allows us
to determine the geometric contribution of the arithmetic intersection
numbers in question using the projection formula (Proposition~\ref{finiteexplicit}).

\medskip

%\item[chap6] 
In Section 6 the arithmetic intersection theory of $
\widetilde{\calH}(N)$ is studied. We introduce the arithmetic
Hirzebruch-Zagier divisors
$\widehat{\mathcal{T}}_N(m)=(\mathcal{T}_N(m), \frakg_N(m))$, where
$\frakg_N(m)$ is the pull-back of the automorphic Green object
described in Section 2. The subgroup of the first arithmetic Chow group generated by
these divisors has finite rank 
%$\left[\frac{D+19}{24}\right]$
(see Theorem \ref{thm:chowfinite}).
We let
$\cca(\overline{\mathcal{M}}_k(\GK(N)))$ be the first arithmetic Chern
class of the pre-log singular hermitian line bundle
$\overline{\mathcal{M}}_k(\GK(N))$ of modular forms of weight $k$ with
the Petersson metric.

  Let $M^+_2(D,\chi_D)$ be the space of holomorphic modular forms of weight $2$ for
  the congruence subgroup $\Gamma_0(D)\subseteq \Sl_2(\Z)$ with
  character $\chi_D=\left(\frac{D}{\cdot}\right)$ satisfying the ``plus-condition'' as in \cite{HZ}.

\theoremstyle{plain}
\newtheorem*{theorema}{Theorem A}
\begin{theorema}%\label{thm:A}
The arithmetic generating series
\begin{align*}
  \widehat{A}(\tau)=\cca\big(\overline{\mathcal{M}}_{1/2}^\vee(\GK(N))\big) +
  \sum_{m>0} \widehat{\mathcal{T}}_N(m) q^m
\end{align*}
%where $q=\e(\tau)$ with $\tau\in \H$, 
is a holomorphic modular form in $M^+_2(D,\chi_D)$ with values in
$\cha^1(\widetilde{\calH}(N),\mathcal{D}_\pre)_\QQ$, i.e., an element
of
$M^+_2(D,\chi_D)\otimes_\Q\cha^1(\widetilde{\calH}(N),\mathcal{D}_\pre)_\QQ$.
\end{theorema}

In particular, via the natural morphism from the arithmetic Chow group
to cohomology, one recovers the classical result of Hirzebruch and
Zagier.

In the following all arithmetic intersection numbers and Faltings
heights are elements of $\RR_N= \RR/ \langle\sum_{p|N} \QQ \log
p\rangle$, where the sum extends over all prime divisors of $N$.  This
is due to the fact that, as an arithmetic variety,
$\widetilde{\calH}(N)$ is defined over $\ZZ[\zeta_N,1/N]$.
%rather than over $\Z[\zeta_N]$.

\newtheorem*{theoremb}{Theorem B}
\begin{theoremb}%\label{thm:B}
     We have the following identities of arithmetic intersection
  numbers in $\RR_N$:
\begin{align*}
\widehat{A}(\tau) \cdot 
\cca(\overline{\mathcal{M}}_k(\GK(N)))^2
= \frac{k^2}{2} d_N  \zeta_K(-1)  \left(
    \frac{\zeta_K'(-1)}{ \zeta_K(-1)} + \frac{\zeta'(-1)}{\zeta(-1)} +
    \frac{3}{2} + \frac{1}{2} \log(D) \right)
\cdot E(\tau).
\end{align*}
Here $d_N= [\GK:\GK(N)] [\QQ(\zeta_N):\QQ]$, and $E(\tau)$ denotes the Eisenstein series defined in \eqref{eis}. 
In particular, the arithmetic
  self intersection number of
  $\overline{\mathcal{M}}_k(\GK(N))$ is  given by:
\begin{align*}
  \overline{\mathcal{M}}_k(\GK(N))^3
    & = -k^3 d_N
\zeta_K(-1)  \left(
    \frac{\zeta_K'(-1)}{ \zeta_K(-1)} + \frac{\zeta'(-1)}{\zeta(-1)} +
    \frac{3}{2} + \frac{1}{2} \log(D) \right).
\end{align*}
\end{theoremb}

\newtheorem*{theoremc}{Theorem C}
\begin{theoremc}%\label{thm:C}
If $T(m)$ is a
  Hirzebruch-Zagier divisor which is
  disjoint to the boundary, then in $\RR_N$  the Faltings height of its model
  $\mathcal{T}_N(m)$ with respect to
  $\overline{\mathcal{M}}_k(\GK(N))$ equals 
\begin{align*}
%\label{eq:fh_tm}
  \fh_{\overline{\mathcal{M}}_k(\Gamma_K(N))} (\mathcal{T}_N(m))&= -
  (2k)^2 d_N \vol(T(m)) \left( \frac{\zeta'(-1)}{\zeta(-1)} +
  \frac{1}{2} + \frac{1}{2} \frac{\sigma_m'(-1)}{\sigma_m(-1)} \right).
\end{align*}
\end{theoremc}

Observe that, after normalizing by the degree $d_{N}$, the
previous formulas do not depend on the level $N$. Indeed our
calculations essentially only depend on the coarse moduli space with
no level-structure. 

From the Theorems B and C we can also obtain equalities in $\mathbb{R}$ and
not just in $\R_N$.  Following a suggestion of S.~Kudla, we consider
the whole tower of Hilbert modular surfaces $\{\widetilde
{\calH}(N)\}_{N\ge 3}$ as a substitute of the stack $\widetilde
{\calH}(1)$. In fact it is enough to choose two integers $N,M\ge 3$
such that $(M,N)=1$ and view the pair $(\widetilde {\calH}(N),
\widetilde {\calH}(M))$ as a Galois \v Cech covering of $\widetilde
{\calH}(1)$. Then, defining an arithmetic cycle on this covering as a
pair of arithmetic cycles, the first on $\widetilde {\calH}(N)$ and
the second on $\widetilde {\calH}(M)$, that agree after pull-back to
$\widetilde {\calH}(NM)$, one obtains an arithmetic Chow ring that has
an arithmetic degree with values in $\R$. Using this ad hoc arithmetic
intersection theory we can give a meaning to the formulas in Theorems
B and C on $\R$ (see Theorem \ref{thm:BR} and Theorem \ref{thm:CR}).

This theory can be viewed as a substitute of an
arithmetic intersection theory on stacks for our 
particular situation. 
Notice that any reasonable arithmetic
intersection theory for the coarse moduli space without level
structure would have to inject into the ad hoc theory so that all
arithmetic degrees and heights would agree.

One could also follow different strategies to obtain 
equalities in $\mathbb{R}$. 
%There are at least two other strategies to obtain equalities in
%$\mathbb{R}$. 
First, we
remark that the the proofs Theorems B and C would directly carry over 
to other level structures. For instance, Pappas constructed regular
models
%$\calH_{00}^{Fil}(\mathcal{A})$ 
over $\Z$ of Hilbert modular surfaces associated to congruence
subgroups $\Gamma_{00}(\mathcal{A})\subset\GK$ \cite{Pa}.  If there
exists a toroidal compactification over $\Z$ of
$\calH_{00}^{Fil}(\mathcal{A})$ (compare loc.~cit.~p.~51), one can
work on that arithmetic variety and derive formulas as in Theorems B and C,
%\ref{thm:B} and \ref{thm:C},
but now as equalities in $\R$.
%% Alternatively we could assume the existence of arithmetic Chow groups
%% for singular arithmetic varieties of some kind, including arithmetic
%% Baily-Borel compactifications for Hilbert modular surfaces, equipped
%% with operational first arithmetic Chern classes for hermitian line
%% bundles.  As our computations involve only such Chern classes,
%% we would also obtain equalities in $\R$.
%% Ultimately, an arithmetic intersection theory for divisors on algebraic 
%% stacks would solve this problem.
%
Alternatively one could try to work directly on the 
coarse moduli space. This  requires  arithmetic Chow groups
for singular arithmetic varieties of some kind (e.g.~for arithmetic
Baily-Borel compactifications of Hilbert modular surfaces), equipped
with  arithmetic Chern class operations for hermitian line
bundles. The existence of such a theory would also lead to equalities in $\R$. 
Ultimately, an arithmetic intersection theory for divisors on algebraic 
stacks would solve the problem.

Many arguments of the present paper can  (probably) be
generalized to Shimura varieties of type $\Orth(2,n)$, where one would
like to study the arithmetic intersection theory of Heegner divisors (also referred to as special cycles, see e.g.~\cite{Ku:MSRI}, \cite{Ku6}).
In fact, the work on this paper already motivated generalizations of
partial aspects.  For instance, automorphic Green functions are
investigated in \cite{BrKue} and the curve lemma in \cite{Br3}.
However, progress on the whole picture will (in the near future) be
limited to the exceptional cases, where regular models of such Shimura
varieties are available.

\medskip
 
%Acknowledgments:  
While working on this manuscript, we benefited from encouraging and
stimulating disscussions with many colleagues.  We would like to thank
them all. In particular, we thank E.~Freitag, E.~Goren, W.~Gubler,
J.~Kramer, S.~Kudla, K.~Ono, G.~Pappas, and T.~Wedhorn for their help.

%%% Local Variables: 
%%% mode: latex
%%% TeX-region: "h-intro2-region"
%%% TeX-master: "h-master2"
%%% End: 

% h-arakelov.tex\\
\section{Arithmetic Chow rings with pre-log-log forms} 

In this section we recall the basic definitions and properties of the arithmetic Chow ring $\widehat{\textrm{CH}}^*(\mathscr{X},\mathcal{D}_{\pre})$ constructed in \cite{BKK}.
This arithmetic Chow ring is an extension of the classical construction due to Gillet and Soul\'e (see e.g.~\cite{SABK}), in which the Green objects satisfy besides the usual logarithmic additional log-log growth conditions.
Such singularities naturally occur if one works with automorphic vector bundles \cite{BKK2}.
%
%In \cite{BKK} several extensions of the classical
%arithmetic Chow ring are defined.  
%In this article we will work with the
%arithmetic Chow ring $\widehat{\textrm{CH}}^*(\mathcal{X},\mathcal{D}_{\pre})$,
%in which the Green objects satisfy certain logarithmic, respectively log-log, growth conditions.  Since there are various notions of
%logarithmic, respectively log-log, singularities in the literature, we
%recall the precise definitions in this section. For a comparison to  
%other definitions we refer the interested reader to loc.~cit.

\subsection{Differential forms with growth conditions}

\begin{notation}
Let $X$ be a complex algebraic manifold of dimension $d$ and 
$D$ a normal crossing divisor of $X$. 
We denote by $\mathscr{E}^*_X$ the sheaf of smooth complex differential 
forms on $X$.
Moreover, we write $U=X\setminus D$, 
and let $j:U\to X$ be the inclusion. 

Let $V$ be an open coordinate subset of $X$ with coordinates 
$z_{1},\dots,z_{d}$; we put $r_{i}=|z_{i}|$. We say that $V$ 
\emph{is adapted to $D$}, if the divisor $D$ is locally given  
by the equation $z_{1}\cdots z_{k}=0$. We assume that the coordinate
neighborhood $V$ is small enough; more precisely, we will assume
that all the coordinates satisfy $r_{i}<1/e^{e}$, which implies 
that $\log 1/r_{i}>e$ and $\log(\log 1/r_{i})>1$.

If $f$ and $g$ are two complex functions, 
we will write $f\prec g$, if there exists a constant $C>0$ 
such that $|f(x)|\le C |g(x)|$ for all $x$ in the domain of 
definition under consideration.
\end{notation} 
\begin{definition}  
\label{def:loglog} 
We say that a smooth complex function $f$ on $X\setminus D$ has 
\emph{log-log growth along $D$}, if we have
\begin{equation}
\label{eq:loglog1}
f(z_{1},\dots,z_{d}) \prec \prod_{i=1}^{k}
\log(\log(1/r_{i}))^{M}
\end{equation}
for any coordinate subset $V$ adapted to $D$ and some positive 
integer $M$. The \emph{sheaf of differential forms on $X$ with 
log-log growth along $D$} is the sub algebra of $j_{\ast}\mathscr
{E}^{\ast}_{U}$ generated, in each coordinate neighborhood $V$ 
adapted to $D$, by the functions with log-log growth along $D$ 
and the differentials
\begin{alignat*}{2} 
&\frac{\dd z_{i}}{z_{i}\log(1/r_{i})},\,\frac{\dd\bar{z}_{i}}
{\bar{z}_{i}\log(1/r_{i})},&\qquad\text{for }i&=1,\dots,k, \\
&\dd z_{i},\,\dd\bar{z}_{i},&\qquad\text{for }i&=k+1,\dots,d.
\end{alignat*}
A differential form with log-log growth along $D$ will be called 
a \emph{log-log growth form}.
\end{definition}

\begin{definition} \label{def:preloglog}
A log-log growth form $\omega$ such that $\partial\omega$,
$\bar{\partial}\omega$ and $\partial\bar{\partial}\omega$ 
are also log-log growth forms is called a \emph{pre-log-log 
form}. The sheaf of pre-log-log forms is the sub algebra of 
$j_{\ast}\mathscr{E}^{\ast}_{U}$ generated by the pre-log-log 
forms. We will denote this complex by $\mathscr{E}^{\ast}_{X}
\langle\langle D\rangle\rangle_{\pre}$.
\end{definition}

In \cite{BKK}, Proposition 7.6, it is shown that pre-log-log forms
are integrable and the currents associated to them do not have
residues.

The sheaf $\mathscr{E}^{\ast}_{X}\langle\langle D\rangle\rangle_
{\pre}$ together with its real structure, its bigrading, and the usual
differential operators $\partial$, $\bar{\partial}$ is easily checked
to be a sheaf of Dolbeault algebras.  We call it the \emph{Dolbeault
  algebra of pre-log-log forms}.  Observe that it is the maximal subsheaf
of Dolbeault algebras of the sheaf of differential forms with log-log
growth.

\begin{definition}  
\label{def:log}
We say that a smooth complex function $f$ on $U$ has \emph{log 
growth along $D$}, if we have
\begin{equation}
\label{eq:log1}
f(z_{1},\dots,z_{d})\prec \prod_{i=1}^{k}\log(1/r_{i})^{M}
\end{equation}
for any coordinate subset $V$ adapted to $D$ and some positive 
integer $M$. The \emph{sheaf of differential forms on $X$ with 
log growth along $D$} is the sub algebra of $j_{\ast}\mathscr
{E}^{\ast}_{U}$ generated, in each coordinate neighborhood $V$ 
adapted to $D$, by the functions with log growth along $D$ and
the differentials
\begin{alignat*}{2}
&\frac{\dd z_{i}}{z_{i}},\,\frac{\dd\bar{z}_{i}}{\bar{z}_{i}},
&\qquad\text{for }i&=1,\dots,k, \\
&\dd z_{i},\,\dd\bar{z}_{i},&\qquad\text{for }i&=k+1,\dots,d.
\end{alignat*}
A differential form with log growth along $D$ will be called a 
\emph{log growth form}.
\end{definition}

\begin{definition}
A log growth form $\omega$ such that $\partial\omega$, $\bar
{\partial}\omega$ and $\partial\bar{\partial}\omega$ are also 
log growth forms is called a \emph{pre-log form}. The sheaf 
of pre-log forms is the sub algebra of $j_{\ast}\mathscr{E}^
{\ast}_{U}$ generated by the pre-log forms. We will denote 
this complex by $\mathscr{E}^{\ast}_{X}\langle D\rangle_{\pre}$.
\end{definition}

The sheaf $\mathscr{E}^{\ast}_{X}\langle D\rangle_{\pre} $ together
with its real structure, its bigrading and the usual differential
operators $\partial$, $\bar{\partial}$ is easily checked to be a sheaf
of Dolbeault algebras. We call it the \emph{Dolbeault algebra of
  pre-log forms}. It is the maximal subsheaf of Dolbeault algebras of
the sheaf of differential forms with log growth.

For the general situation which we are interested in, we need a
combination of the concepts of pre-log-log and pre-log forms.

\begin{notation}
Let $X$, $D$, $U$ and $j$ be as above. Let $D_{1}$ and $D_{2}$ be
normal crossing divisors, which may have common components, and such
that $D=D_{1}\cup D_{2}$. We denote by $D_{2}'$ the union of the
components of $D_{2}$ which are not contained in $D_{1}$. We will say
that the open coordinate subset $V$ is adapted to $D_{1}$ and $D_{2}$,
if $D_{1}$ has the equation $z_{1}\cdots z_{k}=0$, $D_{2}'$ has the equation
$z_{k+1} \cdots z_{l}=0$, and $r_{i}=|z_{i}|<1/e^{e}$ for
$i=1,\dots,d$.
\end{notation}

\begin{definition} 
We define the \emph{sheaf of differential forms with log growth along $D_{1}$ 
and log-log growth along $D_{2}$} to be the sub algebra of $j_{\ast}
\mathscr{E}^{\ast}_{U}$ generated by differential forms with log 
growth along $D_{1}$ and log-log growth along $D_{2}$. 

A differential form with log growth along $D_{1}$ and log-log 
growth along $D_{2}$ will be called a \emph{mixed growth form}, 
if the divisors $D_{1}$ and $D_{2}$ are clear from the context.
\end{definition}

\begin{definition} 
Let $X$, $D=D_{1}\cup D_{2}$, $U$ and $j$ be as before. A 
mixed growth form $\omega$ such that $\partial\omega$, $\bar
{\partial}\omega$ and $\partial\bar{\partial}\omega$ are also 
mixed growth forms is called a \emph{mixed form}. The sheaf of
mixed forms is the sub algebra of $j_{\ast}\mathscr{E}^{\ast}_
{U}$ generated by the mixed forms. We will denote this complex 
by $\mathscr{E}^{\ast}_{X}\langle D_{1}\langle D_{2}\rangle
\rangle_{\pre}$.
\end{definition}

The sheaf $\mathscr{E}^{\ast}_{X}\langle D_{1}\langle D_{2}\rangle 
\rangle_{\pre}$ together with its real structure, its bigrading and 
the usual differential operators $\partial$, $\bar{\partial}$ is 
easily checked to be a sheaf of Dolbeault algebras. We call it the 
\emph{Dolbeault algebra of mixed forms}. Observe that 
we have by definition 
\begin{displaymath}
\mathscr{E}^{\ast}_{X}\langle D_{1}\langle D_{2}\rangle\rangle_{\pre}=
\mathscr{E}^{\ast}_{X}\langle D_{1}\langle D_{2}'\rangle\rangle_{\pre}.
%=\mathscr{E}^{\ast}_{X}\langle D_{1}\rangle_{\pre}\land \mathscr{E}^
%{\ast}_{X}\langle\langle D_{2}\rangle\rangle_{\pre}. 
\end{displaymath}

\subsection{Pre-log-log Green objects}

\begin{notation}
Let $X$ be a complex algebraic manifold of dimension $d$ and $D$  
a normal crossing divisor. We denote by $\underline{X}$ the 
pair $(X,D)$. If $W\subseteq X$ is an open subset, we write 
$\underline{W}=(W,D\cap W)$.

In the sequel we will consider all operations adapted to the pair 
$\underline{X}$. For instance, if $Y\subsetneq X$ is a closed 
algebraic subset and $W=X\setminus Y$, then an embedded resolution 
of singularities of $Y$ in $\underline{X}$ is a proper modification 
$\pi:\widetilde{X}\to X$ such that $\pi\big|_{\pi^{-1}
(W)}:\pi^{-1}(W)\to W$ is an isomorphism, and
\begin{displaymath}
\pi^{-1}(Y),\,\pi^{-1}(D),\,\pi^{-1}(Y\cup D)
%,\,\pi^{-1}(Y\cap D)  
\end{displaymath}
are normal crossing divisors on $\widetilde{X}$. Using Hironaka's 
theorem on the resolution of singularities \cite{Hi}, one 
can see that such an embedded resolution of singularities exists.

Analogously, a normal crossing compactification of $\underline{X}$
will be a smooth compactification $\overline{X}$ such that the 
adherence $\overline{D}$ of $D$, the subset $B_{\overline{X}}=
\overline{X}\setminus X$, and the subset $B_{\overline{X}}\cup
\overline{D}$
%, $B_{\overline{X}}\cap\overline{D}$ 
are normal 
crossing divisors. 
\end{notation}

\begin{definition}
Given a diagram of normal crossing compactifications of $\underline{X}$ 
\begin{displaymath}
\xymatrix{
\overline{X}'\ar[r]^{\varphi}&\overline{X} \\
&X\ar[ul]\ar[u],}
\end{displaymath}
with divisors $B_{\overline{X}'}$ and $B_{\overline{X}}$ at 
infinity, respectively, then by functoriality of mixed forms
there is an induced morphism
\begin{displaymath}
\varphi^{\ast }:\mathscr{E}^{\ast}_{\overline{X}}\langle 
B_{\overline{X}}\langle\overline{D}\rangle\rangle_{\pre}
\longrightarrow\mathscr{E}^{\ast}_{\overline{X}'}\langle 
B_{\overline{X}'}\langle\overline{D}'\rangle\rangle_{\pre}.
\end{displaymath}
In order to have a complex which is independent of the choice of a
particular compactification, we take the limit over all possible
compactifications.  Namely, we define the \emph{complex
  $E^{\ast}_{\pre}(\underline{X})$ of differential forms on $X$,
  pre-log along infinity and pre-log-log along $D$} as
\begin{displaymath}
E^{\ast}_{\pre}(\underline{X})=\lim_{\longrightarrow}\Gamma
(\overline{X},\mathscr{E}^{\ast}_{\overline{X}}\langle 
B_{\overline{X}}\langle\overline{D}\rangle\rangle_{\pre}),
\end{displaymath}
where the limit is taken over all normal crossing compactifications 
$\overline{X}$ of $\underline{X}$.
\end{definition}

Let $X$ be a smooth real variety and $D$ a normal crossing divisor 
defined over $\mathbb{R}$; as before, we write $\underline{X}=(X,D)$. 
For any $U\subseteq X$, the complex $E^{\ast}_{\pre}(\underline{U})$ 
is a Dolbeault algebra with respect to the wedge product.

\begin{notation}
For any Zariski open subset $U\subseteq X$, we put
\begin{displaymath}
\mathcal{D}^{\ast}_{\pre,\underline{X}}(U,p)=(\mathcal{D}^{\ast}_
{\pre,\underline{X}}(U,p),\dd_{\mathcal{D}_\pre})=(\mathcal{D}^{\ast}
(E_{\pre}(\underline{U}_{\mathbb{C}}),p)^{\sigma},\dd_{\mathcal{D}_\pre}),
\end{displaymath}
where $\mathcal{D}^*(E_{\pre}(\underline{U}_{\mathbb{C}}),p)$ is the
Deligne algebra associated to the Dolbeault algebra
$E^{\ast}_{\pre}(\underline{U})$, and $\sigma$ is the antilinear
involution $\omega \mapsto \overline{F_\infty (\omega)}$ (see
\cite{BKK} Definition 7.17). When $\underline{X}=(X,D)$ is clear from
the context we write $\mathcal{D}^{\ast}_{\pre}(U,p)$ instead of 
$\mathcal{D}^{\ast}_{\pre,\underline{X}}(U,p)$. 
\end{notation}

Let
$U\to X$ be an open immersion and $Y=X\setminus U$. We
write
\begin{displaymath}
    \widehat{H}^{n}_{\mathcal{D}_{\pre},Y}(X,p)=
  \widehat{H}^{n}(\mathcal{D}^{\ast}_{\pre,\underline{X}}(X,p),
\mathcal{D}^{\ast}_{\pre,\underline{X}}(U,p)), 
\end{displaymath}
where the latter groups are truncated relative cohomology groups  
(see \cite{BKK} Definition 2.55).  
%as defined in \cite{Bu}.
Recall that a class $\mathfrak{g} \in
\widehat{H}^{n}_{\mathcal{D}_{\pre},Y}(X,p)$ is represented by a pair
$\mathfrak{g}=(\omega ,\widetilde g)$, with $\omega \in 
{\rm Z}(\mathcal{D}_{\pre}^{n}(X,p))$ a cocycle and $\widetilde g \in \widetilde
{\mathcal{D}}_{\pre}^{n-1}(U,p)=\mathcal{D}_{\pre}^{n-1}(U,p)\big / \Img
\dd_{\mathcal{D}_{\pre}} $ , such that $\dd_{\mathcal{D}_{\pre}} \widetilde g
=\omega$.  There are morphisms
\begin{align*}
 \omega :\widehat{H}^{n}_{\mathcal{D}_{\pre},Y}(X,p)&\longrightarrow
  {\rm Z}(\mathcal{D}_{\pre}^{n}(X,p)),
\end{align*}
given by $\omega(\mathfrak{g}) = \omega( \omega, g) = \omega$ 
and
\begin{align*}
  \cl:\widehat{H}^{n}_{\mathcal{D}_{\pre},Y}(X,p)
  &\longrightarrow 
  H^{n}_{\mathcal{D}_\pre,Y}(X,p),
\end{align*}
given by sending the class of the pair $(\omega, g)$ in
$\widehat{H}^{n}_{\mathcal{D}_{\pre},Y}(X,p)$ to its class $[\omega,g]$ in the
cohomology group $H^n(\mathcal{D}_{\pre}^{\ast}(X,p),\mathcal{D}_{\pre}^{\ast}(U,p))$.

\begin{definition} Let $y$ be a $p$-codimensional algebraic cycle on $X$, with
  $\supp y=Y$. A \emph{pre-log-log Green object for $y$} is an element
    $\mathfrak{g}_{y}\in 
    \widehat{H}^{2p}_{\mathcal{D}_{\pre},Y}(X,p)$ such that 
    $$\cl(\mathfrak{g}_{y})=\cl(y)\in
    H^{2p}_{\mathcal{D}_{\pre},Y}(X,p),$$
    here the class $\cl(y)$ is
    given by the image of the class of the cycle $y$ in real
    Deligne-Beilinson cohomology via the natural morphism
    $H^{2p}_{\mathcal{D},Y}(X,\RR(p))\to
    H^{2p}_{\mathcal{D}_{\pre},Y}(X,p)$.
  \end{definition}
  The surjectivity of the morphism $\cl$ implies that any algebraic
  cycle has a pre-log-log Green object. For the convenience of the
  reader we now recall that
$$
{\rm Z}(\mathcal{D}_{\pre}^{2p}(X,p))
= \left \{ \omega \in E_{\pre}^{p,p}(\underline{X}) \cap
  E_{\pre,\mathbb{R}}^{2p}(\underline{X},p)  \mid \dd \omega = 0 \right \}
$$ 
and 
$$
\widetilde{\mathcal{D}}_{\pre}^{2p-1}(U,p)
= \left \{  g \in E_{\pre}^{p-1,p-1}(\underline{U}) \cap
  E^{2p-2}_{\pre,\mathbb{R}}(\underline{U},p-1) \right \} \Big / ( \Img
\partial + \Img \bar \partial),
$$
where $\dd = \partial + \bar \partial$ and $\partial$, respectively $\bar
\partial$, are the usual holomorphic, respectively anti holomorphic,
derivatives.  Then a pre-log-log Green object for $y$ as above
is represented by a pair
\begin{align*}
  (-2 \partial \bar \partial g_y, \widetilde{g}_y) \in
  {\rm Z}(\mathcal{D}_{\pre}^{2p}(X,p)) \oplus
  \widetilde{\mathcal{D}}_{\pre}^{2p-1}(U,p).
\end{align*}

\begin{proposition} 
\label{prop:302}
%(See \cite{BKK} Proposition 7.20.)
Let $\underline{X}=(X,D)$, where $X$ is a proper smooth real variety
and $D$ is a fixed normal crossing divisor. Let  $y$ be a $p$-codimensional
cycle on $X$ with support $Y$.
%Then we have the following statements: 

i)
If the class of a cycle $(\omega,g)$ in $H^{2p}_{\mathcal{D}_
{\pre},Y}(X,p)$ is equal to the class of $y$, then  
\begin{align}
\label{eq:eqgreenwlog}
-2\partial\bar{\partial}[g]=[\omega]-\delta_{y}. 
\end{align}
%Here $[\cdot]$ denotes the current associated to a differential form as defined in \cite{BKK} (5.32).
Here the currents $[\cdot]$ and $\delta_{y}$ are normalized as in
 \cite{BKK} (5.32) and Definition 5.35.

ii)
Assume that $y=\sum_{j}n_{j}Y_{j}$ with irreducible subvarieties 
$Y_{j}$ and certain multiplicities $n_{j}$. If the cycle $(\omega,g)$
represents the class of $y$, then the equality
\begin{equation}
\label{eq:dcwlg} 
-\lim_{\varepsilon\rightarrow 0}\int_{\partial B_{\varepsilon}(Y)}
\alpha\dd^{c}g=\frac{(2\pi i)^{p-1}}{2}\sum_{j}n_{j}\int_{Y_{j}}
\alpha 
\end{equation}
holds for any differential form $\alpha$; here $\dd^{c}=\frac{1}{4\pi
  i}(\partial-\bar \partial) $ and $B_{\varepsilon}(Y)$ 
is an $\varepsilon$-neighborhood of $Y$ such that the orientation 
of $\partial B_{\varepsilon}(Y)$ is induced from the orientation of 
$B_{\varepsilon}(Y)$.
\hfill
$\square$
\end{proposition} 
%\begin{proof} See 
%For the claims (i) and (ii) we refer to
% \cite{BKK} Proposition 7.20.
%\end{proof}

%A pre-log-log form  
%satisfying \eqref{eq:basicgreen} is called a 
%\emph{basic pre-log-log Green form}.

\subsection{Star products of pre-log-log Green objects}
Let $\underline{X}=(X,D)$ be a proper smooth real variety of
dimension $d$ with fixed normal crossing divisor $D$. Moreover, 
let $Y,Z$ be closed subsets of $X$. Then it is shown in \cite{BKK}
that the product $\bullet$ of Deligne-Beilinson cohomology induces a
star product
\begin{align*}   
* : \widehat{H}^{n}_{\mathcal{D}_{\pre},Y}(X,p) \times
    \widehat{H}^{m}_{\mathcal{D}_{\pre},Z}(X,q) \longrightarrow
    \widehat{H}^{n+m}_{\mathcal{D}_{\pre},Y\cap Z}(X,p+q),
\end{align*}
which is graded commutative, associative and compatible with the
morphisms $\omega$ and $\cl$.
We now recall how to calculate the star product of Green objects in
the cases being of importance to us.  For this we fix a cycle $y\in
{\rm Z}^p(X_\RR)$ with support $Y$ and a cycle $z\in {\rm Z}^q(X_\RR)$
with support $Z$ such that $y$
and $z$ intersect properly.

Let $\mathfrak{g}_y =(\omega_{y},\widetilde{g}_{y}) \in
\widehat{H}^{2p}_{\mathcal{D}_{\pre},Y}(X,p)$ and $\mathfrak{g}_z
=(\omega_{z},\widetilde{g}_{z}) \in
\widehat{H}^{2q}_{\mathcal{D}_{\pre},Z}(X,q)$ be 
Green objects for $y$ and $z$, respectively.
%We now recall how to calculate the star product of Green objects in
%the cases being of importance to us.  For this we fix $\mathfrak{g}_y$
%a Green object for a cycle $y\in {\rm Z}^p(X_\RR)$ with support $Y$ and $\mathfrak{g}_Z$ a
%Green object for a cycle $Z\in {\rm Z}^p(X_\RR)$.  We assume $Y$ and $Z$ are
%both irreducible and intersect properly.  
Adapting the argument of \cite{Bu}, we can find an embedded resolution
of singularities of $Y\cup Z$, $\pi :\widetilde {X}_\RR\to X_\RR$,
which factors through embedded resolutions of $Y$, $Z$ and $Y\cap Z$.
In particular, we can assume that
\begin{displaymath}
      \pi ^{-1}(Y),\ 
      \pi ^{-1}(Z),\ 
      \pi ^{-1}(Y\cap Z),\ 
      \pi ^{-1}(Y\cap \DD) \text{, and }  
      \pi ^{-1}(Z\cap \DD) 
\end{displaymath}
are also normal crossings divisors. Let us denote by $\widehat {Y}$
the normal crossings divisor formed by the components of $\pi
^{-1}(Y)$ that are not contained in $\pi ^{-1}(Y\cap Z)$.  Analogously,
we denote by $\widehat {Z}$ the normal crossing divisor formed by the
components of $\pi ^{-1}(Z)$ that are not contained in $\pi
^{-1}(Y\cap Z)$. Then $\widehat {Y}$ and $\widehat {Z}$ are closed
subsets of $\widetilde {X}$ that do not meet.
Therefore, there exist two 
smooth, $F_{\infty}$-invariant functions $\sigma_{_{YZ}}$ and 
$\sigma_{_{ZY}}$ satisfying $0\leq\sigma_{_{YZ}},\sigma_{_{ZY}}
\leq 1$, $\sigma_{_{YZ}}+\sigma_{_{ZY}}=1$, $\sigma_{_{YZ}}=1$ 
in a neighborhood of $\widehat{Y}$, and $\sigma_{_{ZY}}=1$ in 
a neighborhood of $\widehat{Z}$.
%Let $\mathfrak{g}_{y}=(\omega_{y},\widetilde{g}_{y})$, $\mathfrak
%{g}_{z}=(\omega_{z},\widetilde{g}_{z})$, and $\sigma_{_{YZ}}$, 
%$\sigma_{_{ZY}}$ be as above. 
Finally, in the group $\widehat{H}^{2p+2q}
_{\mathcal{D}_{\pre},Y\cap Z}(X,p+q)$, we then have the identity 
\begin{align}\label{eq:starprodformel1}     
\mathfrak{g}_{y}*\mathfrak{g}_{z}=\left(\omega_{y}\land\omega_{z},
(-2\sigma_{_{ZY}}g_{y}\land\partial\bar{\partial}g_{z}-2\partial\bar
{\partial}(\sigma_{_{YZ}}g_{y})\land g_{z})^\sim
\right)\,.  
\end{align}
Moreover, $\mathfrak{g}_{y}*\mathfrak{g}_{z}$ is a pre-log-log Green
object for the algebraic cycle $y \cdot z$, given by the usual geometric
intersection of $y$ and $z$. 

In order to compute the arithmetic degree of an arithmetic 
intersection, we need formulas for the push-forward of certain 
$*$-products in the top degree of truncated cohomology groups.  
We make the convention that for 
$\mathfrak{g}=(\omega,\widetilde g)\in \widehat{H}^{2d+2}
_{\mathcal{D}_{\pre},\emptyset}(X,d)$  we write $\int_{X}
\mathfrak{g}$ instead of $ \int_{X} g$.
 
\begin{theorem}
\label{STAR-PRODUCT} Let $\underline{X}$ be as before and
assume that $D=D_{1} \cup D_{2}$, where $D_{1}$ and $D_{2}$ are normal
crossing divisors of $X$ satisfying $D_{1}\cap D_{2}=\emptyset$.  Let
$y$ and $z$ be cycles of $X$ with corresponding Green objects
$\frakg_y$ and $\frakg_z$ as above.
Assume that $p+q=d+1$, 
i.e., $Y\cap Z=\emptyset$, and that $Y\cap D_{2}=\emptyset$ 
and $Z\cap D_{1}=\emptyset$. Then
\begin{align*}
\frac{1}{(2\pi i)^{d}}\int\limits_{X}\mathfrak{g}_{y}\ast 
\mathfrak{g}_{z} &=\lim_{\varepsilon\to 0}\Bigg(\frac{1}{(2\pi i)^{d}}\!\!\int
\limits_{X\setminus B_{\varepsilon}(D)}\!\!g_{y}\wedge\omega_{z}- 
\frac{2}{(2\pi i)^{d-1}}\!\!\int\limits_{\partial\left(B_{\varepsilon}
(D_{1})\right)}\!\!\left(g_{z}\wedge\dd^{c}g_{y}-g_{y}\wedge\dd^{c}
g_{z}\right)\Bigg) \\
&\phantom{=}{}+\frac{1}{(2\pi i)^{q-1}}\int\limits_{\overline{y\setminus 
(y\cap D_{1})}} g_{z},
\end{align*}
where $B_{\varepsilon}(D_j)$ denotes an $\varepsilon$-neighborhood 
of $D_{j}$ ($j=1,2$), and $B_{\varepsilon}(D)=B_{\varepsilon}(D_{1})
\cup B_{\varepsilon}(D_{2})$. 
\end{theorem}
\begin{proof}   
Let $\widetilde{X}$ be an embedded resolution of $Y\cup Z$ as
described above. We write $Y'$, respectively $Z'$, 
for the strict transforms of $Y$, respectively $Z$; we note that $Y'
\cap Z'=\emptyset$. Furthermore, we write $Y''$ for the strict 
transform of the closure of $Y\setminus(Y\cap D)$. Choosing $\sigma_
{_{YZ}}$ and $\sigma_{_{ZY}}$ as above, we
may assume that $\sigma_{_{YZ}}$ has value $1$ in a neighborhood 
of $D_{1}$, since $D_{1}\cap Z=\emptyset$, and vanishes in a
neighborhood of $D_2$. 

Using $-2\partial\bar 
{\partial}= (4\pi i)\ddc$, we get by means of
\eqref{eq:starprodformel1}:
\begin{equation} \label{eq:starprodformel}
\mathfrak{g}_{y}\ast\mathfrak{g}_{z}=\left(\omega_{y}\wedge 
\omega_{z},4\pi i\left(\ddc (\sigma_{_{YZ}}g_{y})\wedge 
g_{z}+\sigma_{_{ZY}}g_{y}\wedge \ddc g_{z}\right)^{\sim}
\right). 
\end{equation}
In order to perform the following calculations, we put
\begin{displaymath}
X_{\varepsilon}=\widetilde{X}\setminus\left(B_{\varepsilon}(D)
\cup B_{\varepsilon}(Y'')\cup B_{\varepsilon}(Z')\right), 
\end{displaymath}
where $B_{\varepsilon}(\cdot)$ denotes an $\varepsilon$-neighborhood 
of the quantities in question. On $X_{\varepsilon}$ 
one can split up the integral in question
by means of formula (7.32) of \cite{BKK}:  
\begin{align}
\label{eq:keyformulall}
&\int_{X_{\varepsilon}}\left(\ddc (\sigma_{_{YZ}}g_{y})
\wedge g_{z}+\sigma_{_{ZY}}g_{y}\wedge \ddc g_{z}\right)
\notag \\ 
&\qquad =\int_{X_{\varepsilon}}g_{y}\wedge \ddc g_{z}+\int_
{X_{\varepsilon}}\dd \left(g_{z}\wedge\dd^{c}(\sigma_{_{YZ}}
g_{y})-\sigma_{_{YZ}}g_{y}\wedge\dd^{c}g_{z}
\right).
\end{align}
Applying Stokes' theorem to the latter integral and using the
properties of the function $\sigma_{_{YZ}}$ we obtain for sufficiently small
$\varepsilon>0$: 
\begin{align*}
\int_
{X_{\varepsilon}}\dd \left(g_{z}\wedge\dd^{c}(\sigma_{_{YZ}}
g_{y})-\sigma_{_{YZ}}g_{y}\wedge\dd^{c}g_{z}
\right)
%\notag\\
& =-\int\limits_{\partial\left(B_{\varepsilon}(D_{1})\cup 
B_{\varepsilon}(Y'')\right)}\left(g_{z}\wedge\dd^{c}g_{y}-g_{y} 
\wedge\dd^{c}g_{z}\right) \\
& =-\int\limits_{\partial\left(B_{\varepsilon}(D_{1})\right)}
\left(g_{z}\wedge\dd^{c}g_{y}-g_{y}\wedge\dd^{c}g_{z}\right)\notag \\
&\phantom{=}{}-\!\!\!\! \int\limits_{\partial\left(B_{\varepsilon}(Y'')\setminus
(B_{\varepsilon}(Y'')\cap B_{\varepsilon}(D_{1}))\right)}\!\!\!\!\left(g_{z} 
\wedge\dd^{c}g_{y}-g_{y}\wedge\dd^{c}g_{z}\right).
\end{align*}
Taking into account that $(\omega_{y},\widetilde{g}_{y})$ is a Green 
object for $y$ and that $g_{z}$ is a smooth $(n-p,n-p)$-form on $B_
{\varepsilon}(Y'')\setminus(B_{\varepsilon}(Y'')\cap B_{\varepsilon}
(D_{1}))$, we derive from \eqref{eq:dcwlg}:
\begin{align}\label{keyformulaxxl}
-\int\limits_{\partial\left(B_{\varepsilon}(Y'')\setminus(B_
{\varepsilon}(Y'')\cap B_{\varepsilon}(D_{1}))\right)}\hspace{-1cm}
\left(g_{z}\wedge\dd^{c}g_{y}-g_{y}\wedge\dd^{c}g_{z}\right)=\frac
{(2\pi i)^{p-1}}{2}\int\limits_{\overline{y\setminus(y\cap 
D_{1})} } g_{z}+f(\eps).
\end{align}    
Here $f(\eps)$ is a continuous function with $\lim_{\eps\to 0} f(\eps)=0$.
Combining \eqref{eq:starprodformel}, \eqref{eq:keyformulall}, and \eqref{keyformulaxxl}, we finally find  
\begin{align*}
\int\limits_{X_{\varepsilon}}\mathfrak{g}_{y}\ast\mathfrak{g}_{z}
%&= 
%\lim_{\varepsilon\to 0}\int\limits_{X_{\varepsilon}}4\pi i\left
%(\ddc (\sigma_{_{YZ}}g_{y})\wedge g_{z}+\sigma_{_{ZY}}g_{y} 
%\wedge \ddc g_{z}\right) \\
&=\lim_{\varepsilon\to 0}\Bigg(\int\limits_{X\setminus B_{\varepsilon}
(D)}g_{y}\wedge(-2\partial\bar{\partial}g_{z})+(2\pi i)^{p}\int\limits_
{\overline{y\setminus(y\cap D_{1})}}g_{z} \\
&\phantom{=}{}-4\pi i\int\limits_{\partial\left(B_{\varepsilon}(D_{1})  
\right)}\left(g_{z}\wedge\dd^{c}g_{y}-g_{y}\wedge\dd^{c}g_{z}\right)+
f(\eps)\Bigg)\,.  
\end{align*}
Hence the claim follows.
\end{proof}

\subsection{Arithmetic Chow rings with pre-log-log forms}
Let $K$ be a number field, $A$
  a subring of $K$ with field of fractions $K$, and  $\Sigma$ a complete set of complex
  embeddings of $K$ into $\CC$.
 Let $\XX$ be an arithmetic variety over $A$ of (relative)
dimension $d$ over $S=\Spec A$. We let $\XX_\infty = \coprod_{\sigma \in
  \Sigma} \XX_\sigma(\CC)$ and set $\XX_\RR = (\XX_\infty, F_\infty)$.
Let $D_{K}$ be a fixed normal crossing divisor of $\XX_K $. We denote
  by $D$ the induced normal crossing divisor on $\XX_{\mathbb{R}}$. In
  this paragraph we recall basic 
properties of the arithmetic Chow groups
$\cha^{\ast}(\XX,\mathcal{D}_{\pre})$ defined in
\cite{BKK}.
%Here, we made use of the
%convention that whenever $\underline{\XX}$ is clear from the context,
%we write $\mathcal{D}_{\pre}$ instead of $\mathcal{D}_
%{\pre,\underline{\XX}}$.

Let $ {\rm Z}^p(\XX)$ be the group formed by cycles on $\XX$
 of codimension $p$. Given $y \in {\rm Z}^p(\XX)$, we write 
$y_\infty = \coprod_{\sigma \in \Sigma}y_\sigma(\CC)$ and  let 
$Y =\supp y_\infty$.  We define 
$$
 \widehat{H}^{2p}_{\mathcal{D}_{\pre},Y}(\XX,p)
=  \widehat{H}^{2p}_{\mathcal{D}_{\pre},Y}(\XX_\infty,p)^{F_\infty}
$$
and put
\begin{align*}
 \widehat{H}^{2p}_{\mathcal{D}_{\pre},\mathscr{Z}^p}(\XX,p)
= \lim_{\overset{\longrightarrow}{ \mathscr{Z}^p}}
\widehat{H}^{2p}_{\mathcal{D}_{\pre},Y}(\XX,p),
\end{align*}
where in the limit $\mathscr{Z}^p$ is the set of cycles on $\XX_\infty$ of
codimension $\ge p$ ordered by inclusion.  

\begin{definition} The \emph{group of $p$-codimensional arithmetic
    cycles on $\XX$ } is the group
\begin{align*}
  \za^p(\XX, \mathcal{D}_{\pre})
 &= \left\{ (y, \mathfrak{g}_y)  \in {\rm Z}^p(\XX) \oplus
  \widehat{H}^{2p}_{\mathcal{D}_{\pre},\mathscr{Z}^p}(\XX,p) \, \mid
  \cl( \mathfrak{g}_y) = \cl(y_\infty) \right \}.
\end{align*}
Let $w$ be a codimension $p-1$ irreducible subvariety of $\XX$ and let
$h\in k(w)^*$. Write $h_\infty$ for the induced function on $w_\infty$
and set $Y = \supp (\dv (h_\infty))$. Then there is a distinguished
pre-log-log Green object $\mathfrak{g}(h) \in
\widehat{H}^{2p}_{\mathcal{D}_{\pre},Y}(\XX,p)$ for $\dv (h_\infty)$. We
point out that $\mathfrak{g}(h)$ does only depend on the class of
$h_\infty$ in $H^{2p-1}_{\mathcal{D}_\pre}(X\setminus Y, p)$.  We write
$\diva(h) = (\dv (h), \mathfrak{g}(h))$ for this arithmetic cycle and
denote by $\rata^p(\XX,\mathcal{D}_{\pre})$ the subgroup of $
\za^p(\XX,\mathcal{D}_{\pre})$ generated by arithmetic cycles of the form
$\diva(h)$.  Then, the \emph{$p$-th arithmetic Chow group of $\XX$ with
  log-log growth along $D$} is defined by
    \begin{displaymath}
      \cha^p(\XX,\mathcal{D}_{\pre})=\za^p(\XX,\mathcal{D}_{\pre}) 
 \big/ \rata^p(\XX,\mathcal{D}_{\pre}). 
    \end{displaymath}
 \end{definition}

\begin{theorem}
There exists an arithmetic intersection  product 
$$
\cha^p(\XX, \mathcal{D}_{\pre}) \otimes \cha^q(\XX,
\mathcal{D}_{\pre} ) \longrightarrow \cha^{p+q}(\XX,
\mathcal{D}_{\pre} )\otimes_\Z\QQ,
$$
and 
$$\cha^*(\XX, \mathcal{D}_{\pre})_{\mathbb{Q}}=\bigoplus_{p\ge 0}
\cha^p(\XX, \mathcal{D}_{\pre})\otimes_\Z \mathbb{Q}
$$ 
equipped with this product
has the structure of a commutative associative ring. 
% Moreover the the inclusion of \ref{prop:cha-iso} of
%  arithmetic Chow ring  $\cha^*( \XX)_{\mathbb{Q}}$ defined by
%  Gillet-Soul\'e in $\cha^*(\underline
%   {\XX})_{\mathbb{Q}}$ is an monomorphism of rings.
\hfill$\square$
\end{theorem}

We call $\cha^*(\XX, \mathcal{D}_{\pre})_{\mathbb{Q}}$ the
\emph{arithmetic Chow ring of $\XX$ with log-log growth along $D$}. For a detailed description of the arithmetic intersection product we refer to \cite{BKK} Theorems 4.18 and 4.19.
We now briefly discuss the special case $p+q=d+1$.

Let $(y,\mathfrak{g}_y) \in \za^p(\XX,\mathcal{D}_{\pre})$ and
$(z,\mathfrak{g}_z) \in \za^q(\XX,\mathcal{D}_{\pre})$ be such that
$y_\infty$ and $z_\infty$ have proper intersection on $\XX$.  Since 
$p+q=d+1$, this means $ y_\infty \cap z_\infty = \emptyset$, and the
intersection of $y$ and $z$ defines a class $[y\cdot z]_{\fin}$ in the Chow group with finite support $\operatorname{CH}^{d+1}_{\fin} (\XX)$. 
One obtains
\begin{align}\label{eq:chaprod}
\left[ (y,\mathfrak{g}_y)\cdot(z,\mathfrak{g}_z)\right]=\left[ [y\cdot z]_{\fin}, \mathfrak{g}_y *
\mathfrak{g}_z\right]
\in \cha^{d+1}(\XX,\mathcal{D}_{\pre})_\Q.
\end{align}

\begin{definition} Let $K$ be a number field, $\OK$ 
  its ring of integers, and $\Sigma$ a complete set of complex
  embeddings of $K$ into $\CC$. Then $\Spec \OK$ is an arithmetic
  variety, and due to the product formula for $K$, we have as in
  \cite{SABK} a well defined \emph{arithmetic degree map}
\begin{align}
\label{eq:ar-deg}
\dega:\cha^{1}(\Spec \OK,\mathcal{D}_{\pre})\longrightarrow\mathbb{R},
\end{align}
induced by the assignment
\begin{displaymath}
\left(\sum_{\mathfrak{p_{j}}\in S}n_{j}\mathfrak{p}_{j}, 
\sum_{\sigma\in\Sigma}(0,\widetilde{g}_{\sigma})\right)
\mapsto\sum_{\mathfrak{p_{j}}\in S}n_{j}\log\big|\mathcal
{O}_{K}\big/\mathfrak{p}_{j}\big|+\sum_{\sigma\in\Sigma} 
g_{\sigma}.
\end{displaymath}
\end{definition}
In particular, this map is a group homomorphism, which is an
isomorphism in the case $K=\mathbb{Q}$; it is common to identify
$\cha^{1}(\Spec(\Z),\mathcal{D}_{\pre})$ with $\mathbb{R}$.

\begin{remark}\label{rem:localized}
If $N$ is an integer and $A=\OK[1/N]$, then we obtain a
homomorphism
\begin{align*}
\dega:\cha^{1}(S,\mathcal{D}_{\pre})\longrightarrow\RR_N,
\end{align*}
where $\R_N = \R \Big/ \Big< \sum_{p|N} \Q \cdot\log(p) \Big>$, and the sum extends over all prime divisors $p$ of $N$.
\end{remark}

\begin{theorem}
\label{thm:pushforward-point}
 If $\XX$ is a $d$-dimensional projective arithmetic variety
over $A$ and
$\pi:\XX\to S$ the structure morphism, then
there is a push forward morphism
\begin{displaymath}
\pi_{\ast}:\cha^{d+1}(\XX,\mathcal{D}_{\pre})\longrightarrow
\cha^{1}(S,\mathcal{D}_{\pre}).
\end{displaymath}
%% ii)  Let $f : \underline{\XX} \longrightarrow
%% \underline{\mathscr{Y}}$ be a morphism of regular, projective schemes
%% flat over $\Spec \ZZ$ such that the normal crossing divisor
%% $D_X\subset \mathscr{\XX}(\CC)$ maps to the divisor $D_Y \subset
%% \mathscr{Y}(\CC)$. Then there is a Pull-back morphism
%% \begin{eqnarray*}
%% f^*:  \widehat{ CH }^{p}(  \mathscr{Y},\mathcal{D}_\pre ) 
%% \longrightarrow 
%% \widehat{ CH }^{p}(\XX,\mathcal{D}_\pre )_\QQ.
%% \end{eqnarray*} 
\hfill$\square$
\end{theorem}

\begin{remark} 
Let $A$ be as in Remark \ref{rem:localized}, then the morphism 
$$
\dega \pi_*  : \cha^{d+1}(\XX,\mathcal{D}_{\pre}) 
\longrightarrow \RR_N,
$$
is induced by the assignment
%on the level of representatives given by
$$
\sum_{P \in {\rm Z}^{d+1}(\XX)} n_P (P, \mathfrak{g}_P) 
\mapsto
\sum_{P \in {\rm Z}^{d+1}(\XX)} n_P  \bigg(\log \# |k(P)| +\frac{1}{(2 \pi i)^d}
 \int_{\XX_\infty} \mathfrak{g}_P\bigg); 
$$
here we used the convention that for
$\mathfrak{g}_P=(\omega_P,\widetilde{g}_P)$ we write $ \int_{\XX_\infty}
\mathfrak{g}_P$ instead of $ \int_{\XX_\infty} g_P$.
\end{remark}

In order to ease notation, we sometimes write for $\alpha \in \cha^{d+1}(\XX,\mathcal{D}_{\pre})$ simply $\alpha$ instead of $\dega \pi_* (\alpha)$.

%% \begin{lemma} \label{lem:project} Let 
%%   $\iota: \mathcal{T} \longrightarrow \mathcal{X}$ be a morphism of
%%   arithmetic varieties $f: \mathcal{T}\longrightarrow \Spec \ZZ $ and
%%   $\mathcal{X}\longrightarrow \Spec \ZZ $. If $\alpha \in
%%   \operatorname{Z}^{d+1-p}(\mathcal{X})$ is such that $\alpha_\QQ.
%%   \iota_*\mathcal{Z}_\QQ = \emptyset$, then we have the equality
%%   $$
%%   f_*\left( \alpha.  \iota_*\mathcal{Z} \right) = g_* \iota^*
%%   \alpha \in \CH^1(\Spec \ZZ).
%% $$
%% \end{lemma}
%% \begin{proof}
%%   The assertion follows essentially by the projection formula (see
%%   \cite{SABK}\?).  The morphism $\iota: \mathcal{Z} \longrightarrow
%%   \mathcal{X}$ is proper.  Using that $Z_\QQ$ embeds into
%%   $\mathcal{X}$ we find that $\iota_* \mathcal{Z}$ equals
%%   $\iota(\mathcal{Z})$ (see \cite{SABK}\?)).  Applying the projection
%%   formula (see \cite{SABK}\?) yields
%%   $$
%%   \iota_* \iota^*\alpha = \iota_*( \iota^*\alpha.\mathcal{Z}) =
%%   \alpha. \iota(\mathcal{Z}) \in \CH^{d+1}(\mathcal{X}).
%%   $$
%%   Finally, taking into account $g_* = f_* i_*$ (see \cite{SABK}\?),
%%   we derive the claim.
%% \end{proof}

\subsection{Pre-log singular hermitian line bundles and Faltings
 heights} 

Let $\XX$ and $D$ be as in the previous section.

\begin{definition} \label{def:pre-log-log-hlb}
  Let $\mathscr{L}$ be a line bundle on $\XX$ equipped with a
  $F_{\infty}$-invariant singular hermitian metric $\|\cdot\|$ on the
  induced line bundle 
  $\mathscr{L}_\infty$ over $\XX_\infty$. If there is an
  analytic trivializing cover $\{U_{\alpha },s_{\alpha }\}_{\alpha}$
  such that, for all $\alpha $,
\begin{align} \label{eq:pre-log-log-metric}
  - \log\| s_{\alpha }\| \in \Gamma (U_{\alpha
  },\mathscr{E}^0_{\XX_{\infty}}
\left<\left< D \right>\right>_\pre),
\end{align}
then the metric is called a \emph{pre-log singular hermitian metric}.
The pair $(\mathscr{L},\|\cdot\|)$ is called a
\emph{pre-log singular hermitian line bundle} and denoted by
$\overline{\mathscr{L}}$.
\end{definition}

\begin{lemma} \label{lemm:1}
If $\overline{\mathscr{L}}$ is a pre-log singular hermitian
line bundle on $\mathscr{X}$, then, for any rational section $s$ of $\mathscr{L}$, 
\begin{align} \label{eq:ca1}
\big(\dv(s), ( 2\partial \bar \partial \log\|s\|,
 -\log\|s\|)\big) \in \za^1(\XX,\mathcal{D}_\pre).
\end{align}  
\end{lemma}
\begin{proof}
  Let $\|\cdot\|_{0}$ be a $F_{\infty}$-invariant smooth
  hermitian metric on the 
  line bundle   $\mathscr{L}_\infty$. Since the quotient
  $\|s\|/\|s\|_{0}$  does not depend on the section $s$, Definition
  \ref{def:pre-log-log-hlb} implies that
  \begin{displaymath}
    f=-\log(\|s\|/\|s\|_{0})
  \end{displaymath}
  is a pre-log-log function. 
Consequently, $\amap(f)=( 2\partial \bar \partial f, -f)$ is a  pre-log-log Green object for the empty divisor (see \cite{BKK} Section~7.7).
Since $( 2\partial \bar \partial
 \log\|s\|_{0}, -\log\|s\|_{0})$ is a Green object for $\dv(s)$ and 
 \begin{displaymath}
   ( 2\partial \bar \partial \log\|s\|,
 -\log\|s\|)=
( 2\partial \bar \partial \log\|s\|_{0},
 -\log\|s\|_{0})-\amap(f),
 \end{displaymath}
 we obtain that  $( 2\partial \bar \partial
 \log\|s\|, -\log\|s\|)$  is a pre-log-log Green object for
 $\dv(s)$. 
\end{proof}

% In the particular case, where $L = \mathcal{O}_X$ equipped with the
% absolute value $|\cdot|$ , we get for any rational section $h \in
% k(\XX)^*$
% \begin{align*}
% \mathfrak{g}(h) = ( 0, \log|h|)
% \end{align*}
% and $\diva (g) =\left(\dv g, \mathfrak{g}(g)\right)$. 
% Since
% $
% \rata^1(\XX)= \left\{ \diva(g) \,\big | \,g \in  k(\XX)^* \right\}
% $
% we  
It is easy to see that the class of \eqref{eq:ca1} only depends on
the pair $ (\mathscr{L},\|\cdot\|)$.  We denote it  by
$\cca(\overline{\mathscr{L}})$ and call it the \emph{ first arithmetic
  Chern class of $\overline{\mathscr{L}}$}.

\begin{definition} 
  The \emph{arithmetic Picard group} $\pica(\XX,\mathcal{D}_\pre)$ is the
  group of isomorphy classes of pre-log singular hermitian
  line bundles, where the group structure is given by the tensor product. 
\end{definition}

We have an inclusion $\pica(\XX) \subseteq
\pica(\XX,\mathcal{D}_\pre)$, where $\pica(\XX)$ is the arithmetic
Picard group defined by Gillet and Soul\'e. Moreover, the morphism
\begin{align}\label{eq:pic-cha}
  \cca : \pica(\XX,\mathcal{D}_\pre) \longrightarrow
  \cha^1(\XX,\mathcal{D}_\pre)
\end{align}
given by equation \eqref{eq:ca1} is an isomorphism. Finally, given a
pre-log singular hermitian line bundle $\overline{\mathscr{L}}$ on 
an arithmetic variety of (relative) dimension $d$ over $A$, we write
\begin{align*}
\overline{\mathscr{L}}^{d+1}=
\dega \pi_* \left( \cca(\overline{\mathscr{L}})^{d+1}\right),
\end{align*}
and call it the \emph{arithmetic self-intersection number of 
$\overline{\mathscr{L}}$}.

Putting $U=\XX_{\mathbb{R}}\setminus D$, we write ${\rm Z}^{q}_{U} (\XX)$
for the group of the $q$-codimensional cycles $z$ of $\XX$ such that
$z_{\mathbb{R}}$ intersects $D$ properly. 
We introduce $U_{K}=\XX_{K}\setminus D_{K}$ and observe that
there is a natural injective map $ {\rm
  Z}^{q}(U_{K})\to {\rm Z}^{q}_{U}(\XX)$.  In \cite{BKK} 
a height pairing
\begin{equation}
\label{eq:sing-height}
(\cdot\mid\cdot):\,\cha^{p}(\XX,\mathcal{D}_{\pre})\otimes{\rm Z}^{q}_
{U}(\XX)\longrightarrow\cha^{p+q-d}(S,\mathcal{D}_{\pre})_{\mathbb{Q}}
\end{equation}
is defined. 
%%This generalization includes in particular the
%logarithmic heights for points considered by Faltings in
%\cite{Faltings:EaVZ}, \cite{Faltings86:ftavnf}.  
Observe that, since the height of a cycle whose generic part is supported
in $D$, might be infinite, one cannot expect that the height pairing
due to Bost, Gillet, and Soul\'e unconditionally generalizes to a height
pairing between the arithmetic Chow groups $\cha^{p}(\XX,\mathcal
{D}_{\pre})$ and the \emph{whole} group of cycles ${\rm Z}^{q}(\XX)$.

We now let $\XX$, $D$, $U$ be as before and $p$, $q$ integers satisfying
$p+q=d+1$.
Let $z\in{\rm Z}^{q}_{U}(\XX)$ be an irreducible, reduced cycle and 
$\alpha\in\cha^{p}(\XX,\mathcal{D}_{\pre})$. We represent $\alpha$ 
by the class of an arithmetic cycle $(y,\mathfrak{g}_y)$, where 
$y$ is a $p$-codimensional cycle such that $y_{K}$ intersects $z_
{K}$ properly, and where $\mathfrak{g}_{y}=(\omega_{y},\widetilde
{g}_{y})$ is a pre-log-log Green object for $y$. We have
\begin{align}\label{eq:heightdirac}
(\alpha\mid z)=\left[\pi_{\ast}([y\cdot z]_{\fin}),\left(0,
\widetilde{\pi_{\#}(g_{y}\wedge\delta_{z})}\right)\right]\in
\cha^{1}(S,\mathcal{D}_{\pre})_{\mathbb{Q}}.
\end{align}
Here, the quantity $\widetilde{\pi_{\#}(g_{y}\wedge\delta_{z})}$ 
has to be understood as follows: Let $Z=\supp z_{\mathbb{R}}$ 
and $\imath:\widetilde{Z}\to Z$ be a resolution of 
singularities of $Z$ adapted to $D$. Since $y_{K}\cap z_{K}=
\emptyset$, the functoriality of pre-log-log forms shows that
$\imath^{\ast}(g_{y})$ is a pre-log-log form on $\widetilde{Z}$, 
hence it is locally integrable on $\widetilde{Z}$, and we have
\begin{displaymath}
\widetilde{\pi_{\#}(g_{y}\wedge\delta_{z})}=\frac{1}{(2\pi i)^
{p-1}}\int_{\widetilde{Z}}\imath^{\ast}(g_{y}).
\end{displaymath}
The pairing \eqref{eq:sing-height} is now obtained by linearly 
extending the above definitions. 

If we choose a basic pre-log-log
Green form $g_{z}$ for $z$ and put $\mathfrak{g}_{z}=(-2\partial
\bar{\partial}g_{z},\widetilde{g}_z)$, then the height pairing 
\eqref{eq:sing-height} satisfies 
\begin{equation}
\label{eq:log-sing-heightpairing}
(\alpha\mid z)=\pi_{\ast}\left(\alpha\cdot[z,\mathfrak{g}_{z}]
\right)+\amap\left(\pi_{\#}\left(\widetilde{[\omega(\alpha)
\wedge g_{z}]_{\XX_\R}}\right)\right) \in \cha^{1}(S,\mathcal{D}_{\pre})_{\mathbb{Q}}. 
\end{equation}
The height pairing \eqref{eq:sing-height} is of particular interest,
when $\alpha=\cca(\overline{\mathscr{L}})^{p}$ for some pre-log singular hermitian
line bundle $\overline{\mathscr{L}}$ on $\XX$. We call the real number
\begin{align}
\operatorname{ht}_{\overline{\mathscr{L}}}(z)=\dega\left(\cca(\overline{\mathscr{L}})^
{p}\mid z\right)
\end{align}
the \emph{Faltings height of $z$ (with respect to $\overline{\mathscr{L}}$)}.

%%% Local Variables: 
%%% mode: latex 
%%% TeX-region: "arakelov-region" 
%%% TeX-master: "h-master2"
%%% End:

% h-green4.tex:\\
\section{Complex theory of Hilbert modular surfaces}
\label{sec:h-green}

We begin by recalling some basic facts on Hilbert modular surfaces.
This mainly serves to fix notation. For a detailed account we refer to
\cite{Fr} and \cite{Ge}.

Let $K$ be a real quadratic field with discriminant $D$. Let $\OK$ be
its ring of integers and $\frakd=(\sqrt{D})$ the different. We write
$x\mapsto x'$ for the conjugation in $K$, $\tr(x)=x+x'$ for the trace,
and $\norm(x)=xx'$ for the norm of an element. Given an $a \in K$ we will write $a \gg 0$, if $a$ is totally positive. Furthermore, we denote by $\eps_0>1$ the fundamental unit of $K$.

Let $\H=\{z\in \C;\; \Im(z)>0\}$ be the  upper complex half
plane. The group $\Sl_2(\R)\times \Sl_2(\R)$ acts on the product
$\H^2$ of two copies of $\H$ via Moebius
transformations on both factors. As usual we identify
$\Sl_2(K)$ with a subgroup of $\Sl_2(\R)\times \Sl_2(\R)$ by the
embedding $M\mapsto (M,M')$, where $M'=\kzxz{a'}{b'}{c'}{d'}$ denotes
the conjugate of the matrix $M=\kabcd$.  
If $\fraka$ is a fractional ideal of $K$, we write 
\[
\Gamma(\fraka)=\left\{ \kabcd\in \Sl_2(K);\quad a,d\in \OK,\;b\in \fraka^{-1},\; c\in \fraka\right\} 
\]
for the {\em Hilbert modular group} corresponding to $\fraka$. Moreover, we briefly write $\GK=\Gamma(\OK)$, and $\GK(N)$ for the principal congruence subgroup of level $N$, i.e., the kernel of the natural homomorphism $\GK\to \Sl_2(\OK/N\OK)$.
Throughout we use
$z=(z_1,z_2)$ as a standard variable on $\H^2$. We denote its real
part by $(x_1,x_2)$ and its imaginary part by $(y_1,y_2)$.

Let $\Gamma\leq \Sl_2(K)$ be a subgroup which is commensurable with $\GK$. 
The quotient $\Gamma\bs\H^2$ is called the {\em Hilbert modular surface}
associated with $\Gamma$. It is a non-compact normal complex space which can be
compactified by adding the cusps of $\Gamma$, i.e., the $\Gamma$-classes of 
$\P^1(K)$.
By the theory of Baily-Borel, the quotient
\[
X(\Gamma)=\Gamma\bs(\H^2\cup \P^1(K))
\]
together with the Baily-Borel topology can be given the structure of a
normal projective algebraic variety over $\C$. It is called the
Baily-Borel compactification of 
$\Gamma\bs\H^2$.  
Recall that the cusps of $\Gamma(\fraka)$ are in bijection with the ideal classes of $K$ by mapping $(\alpha:\beta)\in\P^1(K)$ to the ideal 
$\alpha\OK+\beta \fraka^{-1}$. So in particular the cusp $\infty$ corresponds 
to the principal class and $0$ to the class of $\fraka^{-1}$. 
For any point $\xi \in \H^2\cup \P^1(K)$, we denote by $\Gamma_\xi$  the stabilizer of $\xi$ in $\Gamma$. If $\xi\in \H^2$, then the quotient $G=\Gamma_\xi/\{\pm 1\}$ 
is a finite cyclic group.  If $|G|>1$, then $\xi$
is called an elliptic fixed point. 
Notice that $\GK$ always has elliptic fixed points of
order $2$ and $3$. On the other hand $\GK(N)$ acts fixed point freely if $N\geq 3$.

If $\fraka$ is a fractional ideal from the principal genus of $K$, there is a fractional ideal $\frakc$ and a totally positive $\lambda\in K$ such that $\fraka=\lambda\frakc^2$. If $M$ denotes a matrix in
$\kzxz{ \frakc^{-1}}{  \frakc^{-1}}{\frakc}{ \frakc}\cap \Sl_2(K)$,
then 
\begin{align}\label{eq:conj}
\zxz{\lambda^{-1}}{0}{0}{1} M \Gamma_K M^{-1} \zxz{\lambda}{0}{0}{1} 
= \Gamma(\fraka).
\end{align}
This induces an isomorphism of algebraic varieties over $\C$,
\begin{align}\label{eq:aiso}
\Gamma_K\bs \H^2\longrightarrow \Gamma(\fraka)\bs\H^2,\quad z\mapsto \kzxz{\lambda^{-1}}{0}{0}{1} M z.
\end{align}
It extends to an isomorphism $X(\GK)\to X(\Gamma(\fraka))$ mapping the cusp $\frakb$ of $X(\GK)$ to the cusp $\frakb/\frakc$ of $X(\Gamma(\fraka))$. In particular the cusp $\frakc$ is mapped to the cusp $\infty$ of $X(\Gamma(\fraka))$, and $\frakc^{-1}$ mapped to the cusp $0$.

\subsection{Desingularization and the curve lemma}
\label{sect:curvelemma}

Throughout we write $\e(z)= e^{2 \pi i z}$. We denote by
$\E(\delta)=\{ q\in \C; \; |q|<\delta\}$ the $\delta$-disc around the
origin and put $\E=\E(1)$. Moreover, we write
$\dot{\E}=\{ q\in \E; \; q\neq 0\}$.

The singular locus $X(\Gamma)^{sing}$ 
of $X(\Gamma)$ consists of the
cusps and the elliptic fixed points.
Throughout we will work with desingularizations of $X(\Gamma)$ such that
the pull-back of the singular locus is a divisor with normal
crossings.  Given such a desingularization 
\[
\pi: \widetilde{X}(\Gamma)
\longrightarrow X(\Gamma),
\]
we denote this divisor by
$$
\D_\Gamma= \pi^*(X(\Gamma)^{sing}).
$$
We now present a local
description of $\widetilde{X}(\Gamma)$
using the ``curve lemma'' due to Freitag.

%{\bf Embedded desingularization of the cusps.} 

Let $\kappa \in \P^1(K)$  be a cusp of $\Gamma$ 
and $g \in \Sl_2(K)$ with
$\kappa = g\infty$. By replacing  $\Gamma$ by the commensurable group
$g^{-1} \Gamma g$, we may assume that $\kappa=\infty$.
%Consider the cusp $\infty $ of $g^{-1} \Gamma
%g$ instead of the cusp $\kappa$ of $\Gamma$.  Therefore it suffices to consider the cusp $\infty$. 
There is a complete $\Z$-module $\frakt$ of $K$ and a finite index subgroup 
$\Lambda$ of the totally positive units of $\OK$ acting on $\frakt$ such that 
$\Gamma_\infty$ has finite index in the semidirect product 
$\frakt\rtimes\Lambda$ (see \cite{Ge} Chapter II.1). 
In particular, if $\gamma\in \Gamma_\infty$, then
$\gamma (z_1,z_2)=(\eps z_1+\mu,\eps'z_2+\mu')$ for some $\mu\in
\frakt$ and some totally positive unit $\eps\in \Lambda$.  A
fundamental system of open neighborhoods of $\infty\in X(\Gamma)$ is
given by
\begin{align}\label{eq:W}
  V_C =\Gamma_\infty\bs \{ (z_1,z_2) \in \H^2 ;\quad y_1 y_2  > C
  \}\cup \infty, \qquad C>0.
\end{align}

% Let $M=
% \Z + w_0 \Z$ with $w_0 \in K$ a quadratic irrationality and let
% $V=\{\eps_V^n : \, n\in\Z\}$ be a subgroup of the group of the totally
% positive units of $K$ such that any $\eps \in V$ stabilizes $M$, i.e.
% $\eps M \subset M$. Then we put
% $$
% G(M,V)= \left \{ \begin{pmatrix} \eps & \mu \\ 0 &1
% \end{pmatrix} 
% \subset \Gl_2(K) ;\quad \eps \in V,\, \mu \in M \right \}.
% $$
% Define an action of $\kzxz{\eps}{\mu}{0}{1} \in G(M,V)$ on $\H^2$
% by $(z_1,z_2) \mapsto (\eps z_1+ \mu, \eps' z_2 + \mu')$. It is known
% that  there exist a quadratic module $G(M,V)$ whose action coincides with the action of $\Gamma_\infty \subset \Gamma$. 

%If $T$ is a divisor on $X(\Gamma)$ such that $\infty \in T$, then,
%according to Hironaka's theory \cite{Hi}, there exists a desingularization $\pi:
%\widetilde{X}(\Gamma)\to X(\Gamma)$ of $X(\Gamma)$ with
%respect to the divisor $T$ such that $\pi^* T$ is a divisor with
%normal crossings.

%In order to get a local description we denote by $\E(\delta)=\{ q\in
%\C; \; |q|<\delta\}$ the $\delta$-disc around the origin and put
%$\E=\E(1)$.  Recall, if $f$ and $g$ are complex numbers or complex
%valued functions, then we write $f\prec g$, if there exists a constant
%$C>0$ such that $|f| \le C |g|$ on the domain under consideration.

Let $C>0$. 
Let $a\in \widetilde{X}(\Gamma)$ be a point with $\pi(a)=\infty$, and 
$U\subset \widetilde{X}(\Gamma)$ be a small open
neighborhood of $a$ such that $\pi(U)\subset V_C$.
% such that the restriction of $\pi$
%induces an isomorphism $U \setminus \D_\Gamma \to W\setminus\{\infty\}$.  
Possibly replacing $U$ by a smaller neighborhood, after a biholomorphic change of coordinates we may 
assume that $U=\E^2$ is the product of two unit discs, $a=(0,0)$, and $\pi^* \infty=\dv(q_1^\alpha q_2^\beta)$ on $U$.
%$\pi^* \dv(F)=\dv(q_1^\alpha q_2^\beta)$ on $U$, 

The desingularization map induces a holomorphic map $\E^2\to V_C$, which
we also denote by $\pi$. If we restrict it to $\dot{\E}^2$, we get a holomorphic 
map $\dot{\E}^2 \to \Gamma_\infty\bs \H^2$.  
Lifting it to the universal covers, 
%$\H \to \dot{\E}$ and $\H^2\to \Gamma_\infty\bs\H^2$% 
we obtain a commutative diagram
 \begin{align}\label{curvecusp}
 \xymatrix{  \H^2\ar[r]^{\Pi}\ar[d]   &    \H^2\ar[d]\\
 \dot{\E}^2 \ar[r]  &\Gamma_\infty\bs \H^2.}
 \end{align}
 Here $\Pi$ is a holomorphic function satisfying
 \begin{align}  \label{per1}
 \Pi(\tau_1 + 1,\tau_2)=  
\left(\begin{matrix} \eps_1 & 0 \\ 0 & \eps_1'\end{matrix}\right)
 \Pi(\tau_1,\tau_2)+ \left(\begin{matrix} \mu_1 \\ \mu_1'
\end{matrix}\right)\\
\label{per2}
 \Pi(\tau_1,\tau_2+1)=  
\left(\begin{matrix} \eps_2 & 0 \\ 0 & \eps_2'\end{matrix}\right)
 \Pi(\tau_1,\tau_2)+ \left(\begin{matrix} \mu_2 \\ \mu_2'
\end{matrix}\right)
\end{align}
with $\mu_1,\mu_2\in \frakt$ and totally positive units $\eps_1,\eps_2 \in \Lambda$.

 \begin{lemma}\label{l2}
   Let $\Pi:\H^2\to\H^2$ be a holomorphic function satisfying
   \eqref{per1} and \eqref{per2}.  Then $\eps_1=\eps_2=1$ and there is
   a holomorphic function $H: {\E}^2\to \H^2$, such that
 \begin{equation}\label{sh1}
 \Pi(\tau_1,\tau_2)=\begin{pmatrix}
 \Pi_1(\tau_1,\tau_2)\\\Pi_2(\tau_1,\tau_2)
 \end{pmatrix}
 = \begin{pmatrix}\mu_1\tau_1+\mu_2\tau_2\\
 \mu_1'\tau_1+\mu_2'\tau_2
 \end{pmatrix}+H(q_1,q_2),
 \end{equation}
 where $q_1=\e( \tau_1)$ and $q_2=\e( \tau_2)$.  Moreover, $\mu_1$,
 $\mu_2$, and their conjugates are non-negative.
 \end{lemma}

 \begin{proof}
   Applying the ``curve lemma'' (\cite{Fr2} Satz 1, Hilfssatz 2) to the
   functions $\H\to\H$ given by $\tau_i\mapsto \Pi_j(\tau_1,\tau_2)$
   with $i,j \in \{1,2\}$, one finds that
   $\Pi(\tau_1,\tau_2)$ has the form \eqref{sh1} with $\mu_1,
   \mu_2\in \frakt$.  The fact that $\Im (\Pi_j(\tau_1,\tau_2))$ is
   positive implies that $\mu_1$, $\mu_2$, and their conjugates have to
   be non-negative.
 \end{proof}

\begin{remark}
 The properties of the Baily-Borel topology on $X(\Gamma)$ imply that
 the exceptional divisor $\pi^*(\infty)$ 
%over the cusp $\infty$ 
contains the component $\{q_j=0\}$, if and only if $\mu_j$ is totally positive.  
%Hence $\mu_1,\mu_2$ are both totally positive,
% if $a$ is an intersection point of $\pi^*(\infty)$. If $a$ is a
% regular point of $\pi^*(\infty)$, then without loss of
% generality $\mu_1$ is totally positive and $\mu_2=0$.
\end{remark}

% {\bf Desingularization of the elliptic fixed points.}  

Let now $\xi\in \H^2$ be an elliptic fixed point of $\Gamma$, $G$ the
 cyclic group $\Gamma_\xi/\{\pm 1\}$, and $n=|G|$. 
%Let $U_\xi
%\subseteq \H^2$ be a sufficiently small neighborhood of $\xi$ such
% that $\Gamma_\xi$ acts on $U_\xi$ and such that
Let  $V\subset\H^2$ be a small open neighborhood of $\xi$ on which
$\Gamma_\xi$ acts. Then $\Gamma_\xi\bs V$ is an open neighborhood of $\xi\in \Gamma\bs\H^2$.
Let $a\in \widetilde{X}(\Gamma)$ be a point with $\pi(a)=\xi$, and 
$U\subset \widetilde{X}(\Gamma)$ be an open
neighborhood of $a$ such that $\pi(U)\subset \Gamma_\xi\bs V$. Without loss of generality we may 
%After a biholomorphic change of coordinates we may
assume that $U=\E^2$ is the product of two unit discs, $a=(0,0)$, and that $\pi^* \xi=\dv(q_1^\alpha q_2^\beta)$ on $U$.
%$\pi^* \dv(F)=\dv(q_1^\alpha q_2^\beta)$ on $U$, 
%and 
%$\pi(U)\subset \Gamma_\xi\bs W_\xi$,
%where $W_\xi\subseteq \H^2$ 
%is a small open neighborhood of $\xi$ on which
%$\Gamma_\xi$ acts.

The desingularization map induces a holomorphic map $\E^2\to \Gamma_\xi\bs V$. 
Arguing as in \cite{Fr3} (Hilfssatz 5.19, pp.~200) and using the fact that $|G|=n$, we get a  
commutative diagram of holomorphic maps
%
%We may restrict it to a holomorphic map 
%$\dot{\E}^2 \to \Gamma_\xi\bs (W_\xi\setminus \{\xi\})$,  
%lift it to a holomorphic map of the covers $\H^2 \to W_\xi\setminus \{\xi\}$, 
%% (2-dimensional analog of Forster Riemannsche Flaechen, Satz 4.17+ Satz 4.9)
%and use the fact that $|G|=n$, to get a  
%commutative diagram of holomorphic maps 
%
\begin{align}\label{curveelliptic}
\xymatrix{\E^2\ar[r]\ar[d]_{q_j\mapsto q_j^n} &V \ar[r]\ar[d]& \H^2 \ar[d] \\
\E^2 \ar[r] & \Gamma_\xi\bs V \ar[r]&\Gamma_\xi\bs \H^2.}
\end{align}
From this one derives an analogue of the curve lemma for the elliptic
fixed points. 

%
% Using the explicit description of the action of $G$ on
%$V$, the power series expansion of $\Pi$ could be described more
%explicitly, but we do not need that.

\bigskip
 
There is a (unique up to a positive multiple) symmetric $\Sl_2(\R)\times \Sl_2(\R)$-invariant K\"ahler metric on $\H^2$.
Its corresponding $(1,1)$-form is given by
\begin{equation}\label{defomega}
\omega =\frac{1}{4 \pi} \left( \frac{\dd x_1   \dd y_1}{y_1^2}
+ \frac{\dd x_2   \dd y_2}{y_2^2} \right).
\end{equation}

\begin{proposition}\label{omegasing}
The form $\omega$ induces a pre-log-log form on $\widetilde{X}(\Gamma)$ with respect to $\D_\Gamma$. 
\end{proposition}

\begin{proof}
We show that if $\kappa$ is a cusp of $\Gamma$, and $a\in \widetilde{X}(\Gamma)$ with $\pi(a)=\kappa$, then $\pi^*\omega$ satisfies the growth conditions of Definition \ref{def:loglog} in a small neighborhood of $a$.
The corresponding assertion for the elliptic fixed points is easy and will be left to the reader.

Without loss of generality we may assume that $\kappa=\infty$ and that $\pi$ looks locally near $a$ as in \eqref{curvecusp}.

By means of the $\Gamma_\infty$-invariant function $\log(y_1 y_2)$, we may write $\omega=\frac{1}{2\pi i}\partial\bar\partial \log(y_1 y_2)$.
Using the notation of Lemma \ref{l2} we see that 
\begin{align}\label{eq:omegasing0}
\Pi^*\begin{pmatrix} y_1\\ y_2
\end{pmatrix}
=-\frac{1}{2\pi }\begin{pmatrix}\mu_1\log |q_1|+\mu_2\log |q_2|\\
\mu_1'\log |q_1|+\mu_2'\log |q_2|
\end{pmatrix}+
\Im \begin{pmatrix}H_1(q_1,q_2)\\ H_2(q_1,q_2)
\end{pmatrix}.
\end{align}
Consequently,
\begin{align}\label{eq:omegasing1}
\pi^*(\log(y_1 y_2)) &= \log\left( -\tfrac{1}{4\pi}g_1(q_1,q_2)\right)
+\log\left(-\tfrac{1}{4\pi} g_2(q_1,q_2)\right),
\end{align}
where 
\begin{align*}
g_1(q_1,q_2)= \mu_1\log |q_1|^2+\mu_2\log |q_2|^2 
-4\pi \Im  H_1(q_1,q_2),\\ 
g_2(q_1,q_2)= \mu_1'\log |q_1|^2+\mu_2'\log |q_2|^2 
-4\pi\Im  H_2(q_1,q_2).
\end{align*}
Hence we find 
\begin{align}
\nonumber
2\pi i \pi^*(\omega)&=\partial\bar\partial \pi^*(\log(y_1 y_2))\\
\nonumber
&=-\frac{1}{ g_1(q_1,q_2)^2}\left( \mu_1 \frac{\dd q_1}{q_1} +  \mu_2 \frac{\dd q_2}{q_2} -4\pi\partial \Im H_1\right)
\left( \mu_1 \frac{\dd \bar q_1}{\bar q_1} +  \mu_2 \frac{\dd \bar q_2}{\bar q_2} -4\pi\bar \partial \Im H_1\right) \\
\nonumber
&\phantom{=}{}-\frac{1}{ g_2(q_1,q_2)^2}\left( \mu_1' \frac{\dd q_1}{q_1} +  \mu_2' \frac{\dd q_2}{q_2} -4\pi \partial \Im H_2\right)
\left( \mu_1' \frac{\dd \bar q_1}{\bar q_1} +  \mu_2' \frac{\dd \bar q_2}{\bar q_2} -4\pi\bar \partial \Im H_2\right) \\
\label{eq:omegasing2}
&\phantom{=}{}-4\pi\frac{\partial\bar\partial \Im H_1}{g_1(q_1,q_2)} -4\pi \frac{\partial\bar\partial \Im H_2}{g_2(q_1,q_2)}. 
\end{align}
Since $H$ is holomorphic, this differential form has 
log-log growth  along $\pi^*(\infty)$.
\end{proof}

It follows that the volume form
\begin{equation}\label{dmu}
\omega^2 = \frac{ 1}{8 \pi^2}
\frac{\dd x_1  \dd y_1}{y_1^2} \,\frac{\dd x_2  \dd y_2}{y_2^2}
\end{equation}
is also a pre-log-log form.
It is well known that
\begin{align}
  \int_{\widetilde{X}(\GK)}\omega^2 = \frac{ 1}{8 \pi^2} \int_{
      \GK\bs \H^2} \frac{\dd x_1 \dd y_1}{y_1^2}\, \frac{\dd x_2 \dd y_2}{y_2^2} =
  \zeta_K(-1)
\label{eq:ZETAVOLUMEN}
\end{align}
(see e.g.~\cite{Ge}, p.~59). Here $\zeta_K(s)$ is the Dedekind zeta
function of $K$.

\subsection{Hilbert  modular forms and the Petersson metric}

Let $k$ be an integer and $\chi$ a character of $\Gamma$. A meromorphic
function $F$ on $\H^2$ is called a
{\em Hilbert modular form} of weight $k$
(with respect to $\Gamma$ and $\chi$), if it satisfies
\begin{align}\label{eq:trafolaw}
F (\gamma z_1, \gamma'z_2) 
= \chi(\gamma)(c z_1 + d)^k (c' z_2 + d' )^k F(z_1,z_2) 
\end{align}
for all $\gamma=\kabcd\in \Gamma$. If $F$ is even holomorphic on $\H^2$, it is called a {\em holomorphic Hilbert modular form}.
Then, by the K\"ocher principle, $F$ is
automatically holomorphic at the cusps.  We denote the vector space
of holomorphic Hilbert modular forms of weight $k$ (with respect to $\Gamma$ and
trivial character) by $M_k(\Gamma)$.  A  holomorphic Hilbert modular form $F$ has a Fourier expansion at the cusp
$\infty$ of the form
\begin{equation}\label{FFourier}
F(z_1,z_2) = a_0 +\sum_{\substack{\nu \in \frakt\\ \nu \gg 0}} 
a_\nu \, \e(\nu z_1 + \nu' z_2), 
\end{equation}
and analogous expansions at the other cusps.    
The sum runs through
all totally positive $\nu$ in a suitable complete $\Z$-module $\frakt$ of $K$. For instance, if $\Gamma=\Gamma(\fraka)$, then $\frakt$ is equal to $\fraka\frakd^{-1}$. Any Hilbert modular form is the quotient of two holomorphic forms. We will say that a Hilbert modular form has rational Fourier coefficients, if it is the quotient of two holomorphic Hilbert modular forms with rational Fourier coefficients.

Meromorphic (holomorphic) modular forms of weight $k$ can be interpreted as rational (global) sections of the sheaf $\calM_k(\Gamma)$ of modular forms. If we write $p:\H^2\to \Gamma\bs \H^2$ for the canonical projection, then the sections over an open subset $U\subset\Gamma\bs \H^2$ are holomorphic functions on $p^{-1}(U)$, which satisfy the transformation law \eqref{eq:trafolaw}. This defines a coherent analytic sheaf on $\Gamma\bs \H^2$, which is actually algebraic.
By the Koecher principle, it extends to an algebraic sheaf on $X(\Gamma)$. By the theory of Baily-Borel, there is a positive integer $n_\Gamma$ such that $\calM_k(\Gamma)$ is a line bundle if $n_\Gamma|k$, and 
\begin{align*}
  X(\Gamma) \cong {\rm Proj} \Bigg( \bigoplus_{\substack{k\geq 0,\,
        n_\Gamma \mid k}} M_k(\Gamma)\Bigg)
\end{align*}
(see \cite{Ge} p.~44, and \cite{Ch} p.~549).  
The {\em line  bundle of modular forms} of weight $k$ (divisible by $n_\Gamma$) on $\widetilde{X}(\Gamma)$ is defined as the pull-back
$\pi^*\calM_{k}(\Gamma)$.
%\begin{align*}
%\pi^*\calM_{k}(\Gamma).
%\end{align*}
By abuse of notation we will also denote it by $\calM_{k}(\Gamma)$. In the same way, if $F$ is a Hilbert modular form of weight $k$, we simply write $F$ for the section $\pi^*(F)$ on $\widetilde{X}$.
The divisor $\div(F)$ of $F$ decomposes into
\begin{align*}
\dv(F)= \div(F)' +  \sum_{j}
n_j E_j.
\end{align*}
Here $\dv(F)'$ denotes the strict transform of the divisor of the
modular form $F$ on $X(\Gamma)$, and $E_j$
are the irreducible components of the exceptional divisor $\D_\Gamma$.
The multiplicities $n_j$ are determined by the orthogonality relations
\begin{align*}
  \dv(F).E_j =0.
\end{align*}

\begin{definition} \label{def:petersson}
If $F\in \calM_k(\Gamma)(U)$ is a rational section over an open subset $U\subset\Gamma\bs\H^2$, we define its {\em Petersson metric} by 
\begin{align*}
\|F(z_1,z_2)\|^2_{\Pet} 
%|F(z_1,z_2)|^2 
%\left( \log|q_1|^2 \log|q_2|^2\right)^k
=|F(z_1,z_2)|^2 (16 \pi^2 y_1 y_2)^k.
\end{align*}
\end{definition}
This defines a hermitian metric on the line bundle of modular forms of weight $k$ on $\Gamma\bs\H^2$.
We now study how it extends to $\widetilde{X}(\Gamma)$.
%, which is not smooth along $\D_\Gamma$. 

\begin{proposition} \label{prop:petersson}
The Petersson metric on the line bundle $\calM_k(\Gamma)$ of modular 
forms on $\widetilde{X}(\Gamma)$ is a pre-log singular hermitian metric (with  respect to  $\D_\Gamma$).
\end{proposition}

\begin{proof} 
We have to verify the conditions of  Definition \ref{def:pre-log-log-hlb} locally for the points of $\D_\Gamma$. Here we only consider the the points above the cusps of $X(\Gamma)$. The corresponding assertion for the elliptic fixed  points (if there are any) is left to the reader.

Let $\kappa$ be a cusp, and $a\in \widetilde{X}(\Gamma)$ with $\pi(a)=\kappa$.
Moreover, let $F$ be a trivializing section of $\calM_k(\Gamma)$ over a small neighborhood of $a$. We have to show that $\log \|F\|_{\Pet}$  satisfies the growth conditions of 
Definition \ref{def:preloglog}.
Without loss of generality we may assume that $\kappa=\infty$ and that $\pi$ looks locally near $a$ as in \eqref{curvecusp}.
It suffices to show that $\pi^*\log(y_1 y_2)$ is a pre-log-log form near $a$.

That $\pi^*\log(y_1 y_2)$ and $\partial\bar\partial\pi^*\log(y_1 y_2)$ have log-log growth along $\D_\Gamma$ follows from \eqref{eq:omegasing1} and \eqref{eq:omegasing2} in the proof of Lemma \ref{omegasing}. Using the notation of that Lemma, we see that
\begin{align*}
\bar\partial\pi^*\log(y_1 y_2)&=
\frac{1}{g_1(q_1,q_2)}\left( \mu_1 \frac{\dd \bar q_1}{\bar q_1} +  \mu_2 \frac{\dd \bar q_2}{\bar q_2} -4\pi\bar \partial \Im H_1\right)\\
&\phantom{=}{}+\frac{1}{g_2(q_1,q_2)}\left( \mu_1' \frac{\dd \bar q_1}{\bar q_1} +  \mu_2' \frac{\dd \bar q_2}{\bar q_2} -4\pi\bar \partial \Im H_2\right).
\end{align*}
Since $H$ is holomorphic, we may infer that $\bar\partial\pi^*\log(y_1 y_2)$ has log-log growth along $\D_\Gamma$. Analogously, we see that $\partial\pi^*\log(y_1 y_2)$ has log-log growth.
\end{proof}

%\begin{rem}\label{CHERNFORM} 
The first Chern form 
$\cc_1(\calM_k(\Gamma)),\|\cdot\|_{\Pet})$ 
of the line bundle $\calM_k(\Gamma)$ equipped with
the Petersson metric is given by
\[
\cc_1( {\calM}_k(\Gamma) ),\| \cdot \|_{\Pet}) =  2\pi i k\cdot \omega.
\] 
  
\begin{definition}\label{def:greenhilbert}
If $F$ is a Hilbert modular form for $\Gamma$, then we denote the Green object for $\div(F)$ by 
\[
\frakg (F) = (2\partial \bar\partial \log\|F\|_{\Pet}, -\log\|F\|_{\Pet}).
\]
\end{definition}

\begin{remark} \label{GRADGEOM} 
  In view of \eqref{eq:ZETAVOLUMEN} the geometric self intersection number $\calM_k(\Gamma)^2$ of
  the line bundle of modular forms of weight $k$ is equal to 
$k^2  \, [\GK: \Gamma]\,\zeta_K(-1)$.
\end{remark}

% \begin{remark} Recall the formula
% \begin{align*}
% \zeta_K(-1)= \frac{1}{60} \sum_{x \in \Z} \sigma_1\left( \frac{D
%     -x^2}{4}\right),
% \end{align*}   
% where $\sigma_1(x)= \sum_{d|x} d$, if $x \in \N$, and $\sigma_1(x)=0$,
% if $x \notin \N$.  The geometric self intersection number of a line
% bundle on an algebraic variety is always an integer.  Thus ${\cal
%   M}_k(\GK)_\infty$ can only be a line bundle, if $k^2$ is divisible
% by the denominator of $\zeta_K(-1)$. This gives a congruence condition
% on the weight, in accordance with the condition $n_\GK |k$ in Theorem
% \ref{BAILYBOREL}.
% \end{remark}

\subsection{Green functions for Hirzebruch-Zagier divisors}
\label{sect:2.3}

From now on we assume that the discriminant $D$ of the real quadratic
field $K$ is a prime. This implies that the fundamental unit $\eps_0$
has norm $-1$, and that narrow equivalence of ideal classes coincides with wide equivalence.
We write $\chi_D$ for the quadratic character
associated with $K$ given by the Legendre symbol $\chi_D(x)=(\frac{D}{x})$. 
The Dirichlet
$L$-function corresponding to $\chi_D$ is denoted by $L(s,\chi_D)$. 
Moreover, we write $\zeta(s)$ for the Riemann zeta function.

We consider the rational quadratic space $V$ of signature $(2,2)$ of matrices $A=\kzxz{a}{\nu}{\nu'}{b}$ with $a,b\in \Q$ and $\nu\in K$, with the quadratic form $q(A)=\det(A)$.
For a fractional ideal $\fraka$ of $K$ we consider the lattices
\begin{align*}
L(\fraka)&=\left\{\kzxz{a}{\nu}{\nu'}{b};\quad \text{$a\in \norm(\fraka)\Z$, $b\in \Z$ and $\nu\in \fraka$}\right\},\\
L'(\fraka)&=\left\{\kzxz{a}{\nu}{\nu'}{b};\quad \text{$a\in \norm(\fraka)\Z$, $b\in \Z$ and $\nu\in \fraka\frakd^{-1}$}\right\}.
\end{align*}
Notice that the dual of $L(\fraka)$ is $\frac{1}{\norm(\fraka)}L'(\fraka)$.
The group $\Sl_2(K)$ acts on $V$ by $\gamma.A= \gamma A {\gamma'}^t$ for $\gamma\in\Sl_2(K)$. Under this action  $\Gamma(\fraka^{-1})$ preserves the lattices $L(\fraka)$ and $L'(\fraka)$.  In particular, one obtains an injective homomorphism $\Gamma(\fraka^{-1})/\{\pm 1\}\to \Orth(L(\fraka))$ into the orthogonal group of $L(\fraka)$. 

Let $m$ be a positive integer.  Recall that the subset
\begin{equation}\label{deftm}
\bigcup_{\substack{A=\kzxz{a}{\nu}{\nu'}{b}\in L'(\fraka) \\\det(A)=m\norm(\fraka)/D}}
 \{(z_1,z_2)\in \H^2;\quad az_1 z_2 +\nu z_1+\nu'z_2+b = 0\}
\end{equation}
defines a $\Gamma(\fraka)$-invariant divisor $T_\fraka(m)$ on $\H^2$, the {\em
  Hirzebruch-Zagier divisor} of discriminant $m$. It is the inverse
image of an algebraic divisor on the quotient $\Gamma(\fraka)\bs\H^2$, which will
also be denoted by $T_\fraka(m)$. Here we understand that all irreducible
components of $T_\fraka(m)$ are assigned the multiplicity $1$. (There is no ramification in codimension $1$.) The divisor $T_\fraka(m)$
is non-zero, if and only if $\chi_D(m)\neq -1$. Since $D$ is prime, $T_\fraka(m)$ is
irreducible, if $m$ is not divisible by $D^2$ (see \cite{HZ}, \cite{Ge} Chapter
V). Moreover, $T_\fraka(m)$ and $T_\fraka(n)$ intersect properly, if and only if $mn$ is not a square.

Since there is only one genus, there is a fractional ideal $\frakc$ and a totally positive $\lambda\in K$ such that $\fraka=\lambda\frakc^2$. If $M$ denotes a matrix in
$\kzxz{ \frakc^{-1}}{  \frakc^{-1}}{\frakc}{ \frakc}\cap \Sl_2(K)$, then 
\[
\left(\kzxz{\lambda^{-1}}{0}{0}{1} M\right)^{-t}.L(\OK)=\frac{1}{\norm(\frakc)}L(\fraka),\qquad \left(\kzxz{\lambda^{-1}}{0}{0}{1} M\right)^{-t}.L'(\OK)=\frac{1}{\norm(\frakc)}L'(\fraka).
\]
This implies that the isomorphism \eqref{eq:aiso} takes $T_\OK(m)$ to $T_\fraka(m)$.

We will be mainly interested in the case that $\fraka=\OK$. To lighten the notation we will briefly write $L=L(\OK)$, $L'=L'(\OK)$, and $T(m)=T_\OK(m)$.

\begin{definition}
Let $m$ be a positive integer with $\chi_D(m)\neq -1$.
If $m$ is the norm of an ideal in $\OK$, then $T(m)$ is a non-compact divisor on $\GK\bs\H^2$, birational to a linear combination of modular curves. In this case we say that $T(m)$ is {\em isotropic}.
If $m$ is not the norm of an ideal in $\OK$, then $T(m)$ is a compact divisor on $\GK\bs\H^2$, birational to a linear combination of Shimura curves. In that case we say that $T(m)$ is {\em anisotropic}.
\end{definition}

These notions are compatible with the description of $\GK\bs\H^2$ and the divisors $T(m)$ as arithmetic quotients corresponding to orthogonal groups of type $\Orth(2,2)$ and  $\Orth(2,1)$, respectively (see Section \ref{sect:density}).  Here isotropic (anisotropic) Hirzebruch-Zagier divisors are given by isotropic (anisotropic) rational quadratic spaces of signature (2,1).

%\subsection{Green functions for Hirzebruch-Zagier divisors}

In \cite{Br1} a certain Green function $\Phi_m(z_1,z_2,s)$ was
constructed, which is associated to the divisor $T(m)$.  We briefly
recall some of its properties.  For $s\in\C$ with $\Re(s)>1$ the
function $\Phi_m(z_1,z_2,s)$ is defined by
\begin{equation}\label{Phim}
\Phi_m(z_1,z_2,s)  = \sum_{\substack{ a,b\in\Z \\ \lambda\in\frakd^{-1}\\ab- \norm(\lambda)=m/D }} Q_{s-1}\left(1+\frac{ |az_1 z_2 +\lambda z_1+\lambda' z_2+b|^2}{ 2y_1 y_2 m /D }\right).
\end{equation}
Here $Q_{s-1}(t)$ is the Legendre function of the second kind
(cf.~\cite{AS} \S8), defined by
\begin{equation}\label{phis}
Q_{s-1}(t)=\int\limits_0^\infty (t+\sqrt{t^2-1} \cosh u)^{-s}\dd u\qquad (t>1,\, 
\Re(s)>0).
\end{equation}
Notice that for $D=m=1$ the function $\Phi_m(z_1,z_2,s)$ equals the
resolvent kernel function (or hyperbolic Green function) for $\Sl_2(\Z)$ (cf.~\cite{He}). Therefore
$\Phi_m(z_1,z_2,s)$ can be viewed as a generalization to Hilbert
modular surfaces.

The sum in (\ref{Phim}) converges normally for $\Re(s)>1$ and
$(z_1,z_2)\in \H^2-T(m)$. This implies that $\Phi_m(z_1,z_2,s)$ is
invariant under $\GK$. It has a Fourier expansion
\begin{equation}\label{fourier1}
\Phi_m(z_1,z_2,s)=u_0(y_1,y_2,s)+\sum_{\substack{\nu\in \frakd^{-1} \\ \nu\neq 0}}u_\nu(y_1,y_2,s)\e(\nu x_1+\nu'x_2),
\end{equation} 
which converges for $y_1 y_2>m/D$ and $(z_1,z_2)\notin T(m)$.  As a
function in $s$ the latter sum over $\nu\neq 0$ converges normally for
$\Re(s)>3/4$.  The constant term is a meromorphic function in $s$ with
a simple pole at $s=1$.  A refinement of these facts can be used to
show that $\Phi_m(z_1,z_2,s)$ has a meromorphic continuation in $s$ to
$\{s\in\C;\;\Re(s)>3/4\}$. Up to a simple pole in $s=1$ it is
holomorphic in this domain (cf.~\cite{Br1} Theorem 1).

We denote by
\[
\Delta_1=y_1^2 \left( \frac{\partial^2}{\partial x_1^2} +
  \frac{\partial^2}{\partial y_1^2}\right),\qquad
\Delta_2=y_2^2 \left( \frac{\partial^2}{\partial x_2^2} +
  \frac{\partial^2}{\partial y_2^2}\right)
\]
the $\Sl_2(\R)\times\Sl_2(\R)$ invariant Laplace operators on $\H^2$.
% and put $\Delta=\Delta_1+\Delta_2$.
The differential equation for $Q_{s-1}(t)$ implies that
\[
\Delta_j \Phi_m(z_1,z_2,s) = s(s-1) \Phi_m(z_1,z_2,s).
\]

Because $Q_{s-1}(t)=-\tfrac{1}{2} \log (t-1)+O(1)$ for $t\to 1$, the
function $\Phi_m(z_1,z_2,s)$ has a logarithmic singularity along
$T(m)$ of type $-\log|f|^2$, where $f$ is a local equation.

\bigskip

The Fourier expansion (\ref{fourier1}) of $\Phi_m(z_1,z_2,s)$ was
determined in \cite{Br1} using a similar argument as Zagier in
\cite{Za1}. It follows from identity (19), Lemma 1, and Lemma 2 of \cite{Br1} that
the constant term  is given by
%
%It can be obtained from the identity
%\begin{equation}\label{Phim2}
%\Phi_m(z_1,z_2,s) = \Phi^0_m(z_1,z_2,s) + 2 \sum_{\nu\in\frakd^{-1}}\left[\sum_{a=1}^{\infty} G_a(m,\nu)b_s^{m/Da^2}(\nu,y_1,y_2) \right]\e(\nu x_1+\nu'x_2).
%\end{equation}
%Here
%\[
%\Phi_m^0(z_1,z_2,s)= \sum_{\substack{ \lambda\in\frakd^{-1}\\ \norm(\lambda)=-m/D }} \sum_{b\in\Z} Q_{s-1}\left(1+\frac{ |\lambda z_1+\lambda'z_2+b|^2}{2y_1 y_2 m/D}\right),
%\]
%and the $b^A_s(\nu,y_1,y_2)$ denote the Fourier coefficients of the
%auxiliary function
%\begin{equation}\label{HAs}
%H^A_s(z_1,z_2) = \sum_{\theta\in\OK}  Q_{s-1}\left(1+\frac{ |(z_1+\theta)(z_2+\theta')+A|^2}{  2y_1 y_2 A}\right).
%\end{equation}
%Moreover, $G_a(m,\nu)$ is the finite exponential sum
%\begin{equation}\label{Ga}
%G_a(m,\nu)=\sum_{\substack{\lambda\in \frakd^{-1}/a\OK \\ \norm(\lambda)\equiv -m/D \pod{a\Z}}} e\left(\frac{\tr(\nu\lambda)}{a}\right).
%\end{equation}
%Notice that our $G_a(m,\nu)$ equals $G_a(-m,\nu)$ in the notation of
%\cite{Za1}.  The Fourier expansion of $\Phi_m^0(z_1,z_2,s)$ is
%determined in \cite{Br1} Lemma 1, and the coefficients
%$b^A_s(\nu,y_1,y_2)$ are computed \cite{Br1} Lemma 2.  We are mainly
%interested in the constant term $u_0(y_1,y_2,s)$ of the Fourier
%expansion (\ref{fourier1}).  According to (\ref{Phim2}) and the
%results cited above it equals
%
\begin{align}
  \nonumber
  u_0(y_1,y_2,s)
%&=\text{constant term of $\Phi_m^0(z_1,z_2,s)$} +2 \sum_{a=1}^\infty G_a(m,0)b_s^{m/Da^2}(0,y_1,y_2)\\
%  \nonumber
  &= \frac{2\pi}{(2s-1)} \sum_{\substack{ \lambda\in\frakd^{-1}\\ \norm(\lambda)=-m/D}} \max(|\lambda y_1|,|\lambda' y_2|)^{1-s}\min(|\lambda y_1|,|\lambda' y_2|)^s\\
\label{formel1}
&\phantom{=}{}+ \frac{\pi \Gamma(s-1/2)^2}{\sqrt{D}\Gamma(2s)}(4m/D)^s (y_1 y_2)^{1-s}\sum_{a=1}^\infty G_a(m,0)a^{-2s}.
\end{align}
Here $G_a(m,\nu)$ is the finite exponential sum
\begin{equation}\label{Ga}
G_a(m,\nu)=\sum_{\substack{\lambda\in \frakd^{-1}/a\OK \\ \norm(\lambda)\equiv -m/D \pod{a\Z}}} e\left(\frac{\tr(\nu\lambda)}{a}\right).
\end{equation}
Notice that our $G_a(m,\nu)$ equals $G_a(-m,\nu)$ in the notation of \cite{Za1}.
For the purposes of the present paper we need to compute $u_0(y_1,y_2,s)$ more
explicitly.

We define a generalized divisor sum of $m$ by
\begin{align}\label{def:sig}
\sigma_m(s)=m^{(1-s)/2}\sum_{d\mid m} d^{s}\left(\chi_D(d)+\chi_D(m/d)\right).
\end{align}
It satisfies the functional equation $\sigma_m(s)=\sigma_m(-s)$.
If $p$ is a prime and $n\in \Z$, then we denote by $v_p(n)$ the additive $p$-adic valuation of $n$.

\begin{lemma}\label{lem:sigep}
i) If $m=m_0 D^\delta$ with $(m_0,D)=1$, then $\sigma_m(s)$ has the Euler product expansion
\[
\sigma_m(s)=m^{(1-s)/2}\left(1+\chi_D(m_0) D^{\delta s}\right)\prod_{\substack{\text{$p$ prime}\\ p\mid m_0}} \frac{1-\chi_D(p)^{v_p(m_0)+1} p^{(v_p(m_0)+1)s}}{1-\chi_D(p) p^{s}}.
\]
ii) If $m$ is square-free and coprime to $D$, then 
\begin{align}\label{eq:L_p}
-2\frac{\sigma_m'(-1)}{\sigma_m(-1)}=\sum_{\substack{\text{$p$ prime}\\ p\mid m} }
\frac{p-\chi_D(p)}{p+\chi_D(p)} \log(p).
\end{align}
\end{lemma}

\begin{proof}
The first formula follows from the multiplicativity of
$\sigma_m(s)$. The second (and generalizations of it) can be obtained
from the Euler product expansion in a straightforward way.
\end{proof}

\begin{lemma} 
We have
\begin{align*}
m^{s/2}\sum_{a=1}^\infty G_a(m,0)a^{-s}
%&=\frac{\zeta(s-1)}{L(s,\chi_D)}\left(1+\chi_D(m_0) D^{\delta(1-s)}\right)\prod_{\substack{\text{$p$ prime}\\ p\mid m_0}}
%\frac{1-\chi_D(p)^{v_p(m_0)+1} p^{(v_p(m_0)+1)(1-s)}}{1-\chi_D(p) p^{1-s}}
%\\
&=\frac{\zeta(s-1)}{L(s,\chi_D)}
\sigma_m(1-s).
\end{align*}
\end{lemma} 

\begin{proof} 
The exponential sum $G_a(m,0)$ is equal to
\begin{align*}
  G_a(m,0)&=\#\{\lambda\in \OK \pmod{a\frakd};\quad \norm(\lambda)\equiv -m \pmod{aD}\}\\
  &= D^{-1} \#\{\lambda\in \OK \pmod{aD};\quad \norm(\lambda)\equiv -m \pmod{aD}\}\\
  &= D^{-1} N_{aD}(-m),
\end{align*}
with
\[
N_b(n) =\#\{\lambda\in \OK/b\OK;\quad \norm(\lambda)=n \pmod{b}\}
\]
as in \cite{Za1} p.~27.  
It is easily seen that $N_b(n)$ is multiplicative in $b$.
Hence it suffices to determine $N_b(n)$ for prime powers $b=p^r$. We get the following Euler product
expansion:
\begin{align}
\nonumber
  \sum_{a=1}^\infty G_a(m,0)a^{-s}&=\frac{1}{D}\sum_{a=1}^\infty N_{aD}(-m)a^{-s}\\
  &= \frac{1}{D}\left(\sum_{r=0}^\infty N_{D^{1+r}}(-m)D^{-rs}\right)
  \prod_{\substack{\text{$p$ prime}\\p\neq D}} \left( \sum_{r=0}^\infty
    N_{p^{r}}(-m)p^{-rs}\right).
\label{eq:eulprod}
\end{align}

%We now compute the local Euler factors. 
%For simplicity we restrict ourselves to the case that $m$ is square
%free and coprime to $D$ (such that $m=m_0$). The general case can be handled in the same
%way.
%\[
%N_{D^{1+r}}(-m)=(1+\chi_D(m))D^{1+r}.
%\]
%If $q$ is a prime dividing $m$, then
%\[
%N_{q^r}(-m)=\begin{cases} 1,& r=0,\\
%q(1+\chi_D(q))-\chi_D(q),& r=1,\\
%q^{r-1}(q-1)(1+\chi_D(q)),& r>1.
%\end{cases}
%\]
%If $q$ is a prime which is coprime to $mD$, then
%\[
%N_{q^r}(-m)=\begin{cases} 1,& r=0,\\
%q^{r-1}(q-\chi_D(q)),& r>0.
%\end{cases}
%\]
%Using these formulas, we may determine the factors of the Euler
%product above. We find

The function $N_{p^r}(n)$ can be
determined explicitly by means of \cite{Za1}, Lemma 3. By a straightforward computation we find that the local Euler factors are equal to
\begin{align*}
\sum_{r=0}^\infty N_{D^{1+r}}(-m)D^{-rs} &= \frac{D}{1-D^{1-s}}\left(1+\chi_D(m_0)D^{\delta(1-s)}\right),\\
\sum_{r=0}^\infty N_{p^{r}}(-m)p^{-rs} &= \frac{1-\chi_D(p)p^{-s}}{1-p^{1-s}}
\cdot\frac{1-\chi_D(p)^{v_p(m_0)+1} p^{(v_p(m_0)+1)(1-s)}}{1-\chi_D(p) p^{1-s}},
\end{align*}
for $p$ prime with $(p,D)=1$.
Inserting this into \eqref{eq:eulprod}, we obtain the assertion by means of Lemma \ref{lem:sigep}.
%\begin{align*}
%  \sum_{a=1}^\infty G_a(m,0)a^{-s}&=\frac{\zeta(s-1)}{L(s,\chi_D)}(1+\chi_D(m)) \prod_{q\mid m}\left(1+\chi_D(q)q^{1-s}\right)\\
%  &=\frac{\zeta(s-1)}{L(s,\chi_D)}\sum_{d\mid m}
%  d^{1-s}\left(\chi_D(d)+\chi_D(m/d)\right).
%\end{align*}
%This concludes the proof of the Lemma.
\end{proof}

Hence the second term in (\ref{formel1}) is equal to
\begin{equation}\label{zwres}
 \frac{\pi \Gamma(s-1/2)^2}{\sqrt{D}\Gamma(2s)}(4/D)^s (y_1 y_2 )^{1-s}\frac{\zeta(2s-1)}{L(2s, \chi_D)}\sigma_m(2s-1).
\end{equation}
By virtue of the functional equation
\[
L(2s, \chi_D)= \left(\frac{\pi}{D}\right)^{2s-1/2} \frac{\Gamma(1/2-s)}{\Gamma(s)} L(1-2s,\chi_D)
\]
we may rewrite (\ref{zwres}) in the form 
\[
4^s\pi^{3/2-2s}\left(\frac{y_1 y_2}{D}\right)^{1-s}\frac{\Gamma(s-1/2)^2\Gamma(s)}{\Gamma(2s)\Gamma(1/2-s)} \frac{\zeta(2s-1)}{L(1-2s,\chi_D)}
\sigma_m(2s-1).
\]
Using the Legendre duplication formula, $\Gamma(s-1/2)\Gamma(s)=\sqrt{\pi} 2^{2-2s} \Gamma(2s-1)$,
we finally obtain for the second term in (\ref{formel1}):
\begin{equation}\label{formel12}
-2\left(\frac{\pi^2 y_1 y_2}{D}\right)^{1-s} \frac{\Gamma(s-1/2)}{\Gamma(3/2-s)}\frac{\zeta(2s-1)}{L(1-2s, \chi_D)}\sigma_m(2s-1).
\end{equation}

\bigskip

In the following we compute the first summand in (\ref{formel1}).

The subset
\begin{equation}\label{defsm}
S(m)=\bigcup_{\substack{\lambda\in \frakd^{-1}\\ \norm(\lambda)=-m/D}}
\{(z_1,z_2)\in\H^2;\quad \lambda y_1 +\lambda'y_2=0\}
\end{equation}
of $\H^2$ is a union of hyperplanes of real codimension $1$. It is
invariant under the stabilizer of the cusp $\infty$.  The connected
components of $\H^2-S(m)$ are called the {\em Weyl chambers} of
discriminant $m$.  For a subset $W'\subset\H^2$ and $\lambda\in K$ we
write $(\lambda,W')>0$, if $\lambda y_1+\lambda'y_2>0$ for all
$(z_1,z_2)\in W'$.

Let $W\subset\H^2$ be a fixed Weyl chamber of discriminant $m$ and
$W'\subset W$ a non-empty subset.  There are only finitely many
$\lambda\in \frakd^{-1}$ such that $\lambda>0$, $\norm(\lambda)=-m/D$,
and
\[
(\lambda,W')<0, \qquad (\eps_0^2 \lambda,W')>0.
\]
Denote the set of these $\lambda$ by $R(W',m)$.  It is easily seen that
$R(W',m)=R(W,m)$ for all non-empty subsets $W'\subset W$.  By
Dirichlet's unit theorem the set of all $\lambda\in \frakd^{-1}$ with
$\norm(\lambda)=-m/D$ is given by
\[
\{ \pm \lambda \eps_0^{2n};\quad \lambda\in R(W,m),\; n\in \Z\}.
\]
For $(z_1,z_2)\in W$,  $\lambda\in R(W,m)$, and $n\in \Z$ we have
\[
\max ( |\lambda \eps_0^{2n} y_1|, |\lambda'{\eps_0'}^{2n} y_2|)
=\begin{cases} \lambda \eps_0^{2n} y_1, & \text{if $n>0$,}\\
-\lambda'{\eps_0'}^{2n} y_2, & \text{if $n\leq 0$.}
\end{cases}
\]
On $W$ we may rewrite the first summand in (\ref{formel1}) as follows:
\begin{align*}
&
\frac{4\pi}{2s-1} \sum_{\lambda\in R(W,m)} \sum_{n\in \Z} \max(|\lambda\eps_0^{2n} y_1|,|\lambda' {\eps_0'}^{2n}y_2|)^{1-s}\min(|\lambda \eps_0^{2n}y_1|,|\lambda'{\eps_0'}^{2n} y_2|)^s\\
&=\frac{4\pi}{2s-1} \sum_{\lambda\in R(W,m)} \left(
\sum_{n\geq 1} (\lambda \eps_0^{2n} y_1)^{1-s}(-\lambda'{\eps_0'}^{2n}y_2 )^s
+\sum_{n\geq 0} (-\lambda'{\eps_0'}^{-2n}y_2)^{1-s} (\lambda  \eps_0^{-2n}y_1)^{s}\right)\\
%&=\frac{4\pi}{2s-1} \sum_{\lambda\in R(W,m)}  \left(
%(\lambda y_1)^{1-s}(-\lambda' y_2 )^s \sum_{n\geq 1} \eps_0^{2n(1-2s)}
%+(-\lambda' y_2)^{1-s} (\lambda  y_1)^{s}    \sum_{n\geq 0} \eps_0^{2n(1-2s)}
%\right)\\
&=\frac{4\pi}{2s-1} \sum_{\lambda\in R(W,m)}  \left(
(\lambda y_1)^{1-s}(-\lambda' y_2 )^s \frac{ \eps_0^{2-4s}}{1-\eps_0^{2-4s}}
+(-\lambda' y_2)^{1-s} (\lambda  y_1)^{s}    \frac{1}{1- \eps_0^{2-4s}}
\right).
\end{align*}
We summarize the above computations in the following theorem.

\begin{theorem}\label{satzformel2}
  Let $W\subset\H^2$ be a Weyl chamber of discriminant $m$. For
  $(z_1,z_2)\in W$ the constant term of the Fourier expansion of
  $\Phi_m(z_1,z_2,s)$ is given by
\begin{align}
  \nonumber 
u_0(y_1,y_2,s)&=
%  -2\left(\frac{\pi^2 m y_1 y_2}{D}\right)^{1-s}
%  \frac{\Gamma(s-1/2)}{\Gamma(3/2-s)}\frac{\zeta(2s-1)}{L(1-2s,
%    \chi_D)}\sum_{d\mid m}
%  d^{2s-1}\left(\chi_D(d)+\chi_D(m/d)\right)\\
2 \zeta(2s-1) \varphi_m(s) \left(\pi^2 y_1 y_2/D\right)^{1-s}\\
  &\phantom{=}{}+\frac{4\pi}{2s-1}\cdot \frac{1}{1- \eps_0^{2-4s}}
  \sum_{\lambda\in R(W,m)} \left(\eps_0^{2-4s} (\lambda
    y_1)^{1-s}(-\lambda' y_2 )^s +(-\lambda' y_2)^{1-s} (\lambda
    y_1)^{s} \right),
\label{formel2}
\end{align}
where 
\begin{equation}\label{defvarphi}
\varphi_m(s)=-\frac{\Gamma(s-1/2)}{\Gamma(3/2-s)}\frac{1}{L(1-2s, \chi_D)}\sigma_m(2s-1).
\end{equation}
\hfill$\square$
\end{theorem}

As a function in $s$, the Green function $\Phi_m(z_1,z_2,s)$ has a
simple pole at $s=1$ coming from the factor $\zeta(2s-1)$ in the first
term of $u_0(y_1,y_2,s)$.  However, it can be regularized at this
place by defining $\Phi_m(z_1,z_2)$ to be the constant term of the
Laurent expansion of $\Phi_m(z_1,z_2,s)$ at $s=1$ (see \cite{Br1}
p.~66).  

Using the Laurent expansion 
$\zeta(2s-1)=\frac{1}{2}(s-1)^{-1}-\Gamma'(1)+O(s-1)$ 
%\begin{align*}
%\zeta(2s-1)&=\frac{1}{2}(s-1)^{-1}-\Gamma'(1)+O(s-1),\\
%\varphi_m(s)&=\varphi_m(1) + \varphi_m'(1)(s-1) + O((s-1)^2),\\
%\left(\pi^2  y_1 y_2/D\right)^{1-s}&= 1-\log( \pi^2  y_1 y_2/D)(s-1) + O((s-1)^2),
%\end{align*}
we get at $s=1$ the expansion
\begin{align}
  \nonumber 2 \zeta(2s-1) \varphi_m(s)\left(\pi^2  y_1
    y_2/D\right)^{1-s}&=\varphi_m(1)(s-1)^{-1}- \varphi_m(1)
  \log(16\pi^2  y_1 y_2)\\
&\phantom{=}+L_m +O(s-1),
\label{constexp}
\end{align}
where
\begin{align*}
L_m & =\varphi_m'(1)-\varphi_m(1)\big( 2\,\Gamma'(1) -\log( 16 D)\big).
\end{align*}
By means of the Laurent expansion of $\Gamma(s)$ one infers that $L_m$
is more explicitely given by
\begin{align}\label{eq:L_m}
 L_m &=\varphi_m(1) \left(2\frac{L'(-1,\chi_D)}{L(-1,\chi_D)} -
 2\frac{\sigma_m'(-1)}{\sigma_m(-1)}
     +\log(D) \right).
\end{align}

In later applications it will be convenient to write the regularized
function $\Phi_m(z_1,z_2)$ as a limit.
In view of \eqref{constexp} find that
\[
\Phi_m(z_1,z_2)=\lim_{s\to 1}\left( \Phi_m(z_1,z_2,s)- \frac{\varphi_m(1)}{s-1}\right).
\]
%Notice that 
%\[
%\varphi_m(1)= -\frac{1}{L(-1, \chi_D)}\sum_{d\mid m} d\left(\chi_D(d)+\chi_D(m/d)\right).
%\]
The Fourier expansion of $\Phi_m(z_1,z_2)$ can be deduced from
(\ref{fourier1}) by virtue of Theorem \ref{satzformel2} and
(\ref{constexp}). It is given by
\[
\Phi_m(z_1,z_2)=\lim_{s\to 1}\left( u_0(y_1,y_2,s)- \frac{\varphi_m(1)}{s-1}\right)
+\sum_{\substack{\nu\in \frakd^{-1} \\ \nu\neq 0}}u_\nu(y_1,y_2,1) \e(\nu x_1+\nu'x_2).
\]
On a Weyl chamber $W$ of discriminant $m$ we get
\begin{align}\label{fourier2}
\Phi_m(z_1,z_2)&=L_m-\varphi_m(1)\log(16\pi^2  y_1 y_2)+
4\pi (\rho_W y_1+\rho_W' y_2)\\
\nonumber
&\phantom{=}{}
+\sum_{\substack{\nu\in \frakd^{-1} \\ \nu\neq 0}}u_\nu(y_1,y_2,1)
\e(\nu x_1+\nu'x_2),
\end{align}
where
\begin{align}\label{eq:wv}
\rho_W = \frac{ \eps_0}{\tr(\eps_0)}
\sum_{\lambda\in R(W,m)} \lambda
\end{align}
is the so called {\em Weyl vector} associated with $W$.

\begin{remark}
If we compare \eqref{fourier2} with the formula given in \cite{Br1} p.~67 we see
that the quantities $L$ and $q_0(m)$ in \cite{Br1} are given by
$L= L_m-\varphi_m(1)\log (16\pi^2)$ and $q_0(m) = -2  \varphi_m(1)$. Moreover, 
the Weyl vector $\rho_W$ in \cite{Br1} is equal to \eqref{eq:wv} above.
%\begin{align*}
%L&= L_m-\varphi_m(1)\log (16\pi^2),\\
%q_0(m) &= -2  \varphi_m(1),\\
%\rho_W &=\frac{ \eps_0}{\tr(\eps_0)}\sum_{\lambda\in R(W,m)} \lambda.
%\end{align*}
\end{remark}

In order to get a Green function with a ``good'' arithmetic normalization, which is compatible with our normalization of the Petersson metric, we have to renormalize as follows.  

\begin{definition}\label{def:Gm}%\label{def:L_m}
We define the normalized Green function for the divisor $T(m)$ by
\[
G_m(z_1,z_2)=\frac{1}{2}\left( \Phi_m(z_1,z_2)-L_m\right).
\]
\end{definition}

%% \begin{proposition}\label{prop:L_m} 
%% More explicitly, we have 
%% \begin{align*}
%%  L_m &=\varphi_m(1) \left(2\frac{L'(-1,\chi_D)}{L(-1,\chi_D)} -
%%  2\frac{\sigma_m'(-1)}{\sigma_m(-1)}
%%      +\log(D) \right).
%% \end{align*}
%% If $m$ is square-free and coprime to $D$, then 
%% \begin{align}\label{eq:L_p}
%% L_m&=\varphi_m(1) \Bigg(2\frac{L'(-1,\chi_D)}{L(-1,\chi_D)} +
%%  \sum_{\substack{\text{$p$ prime}\\ p\mid m} }
%% \frac{p-\chi_D(p)}{p+\chi_D(p)} \log(p)+
%%     \log(D) \Bigg). 
%% \end{align}
%% \end{proposition}
%% \begin{proof}
%% The first formula follows from the definition of $\varphi_m(s)$ by a straightforward computation. The second (and generalizations of it) can be obtained using the Euler product expansion of $\sigma_m(s)$.
%% \end{proof}
%
%In order to describe $\varphi_m'(1)$ explicitly, in view of \eqref{defvarphi} we need the following
%values of derivatives at $s=1$:
%\begin{align*}
%\frac{d}{ds}\frac{\Gamma(s-1/2)}{\Gamma(3/2-s)}\bigg|_{s=1}&=
%2\Gamma'(1)-\log 16,\\
%\frac{d}{ds}\frac{1}{L(1-2s,\chi_D)}\bigg|_{s=1}&=
%2\frac{L'(-1,\chi_D)}{L(-1,\chi_D)^2},\\
%\frac{d}{ds}\sum_{d\mid m} d^{2s-1}\left(\chi_D(d)+\chi_D(m/d)\right)\bigg|_{s=1}&=
%2\sum_{d\mid m} d\left(\chi_D(d)+\chi_D(m/d)\right)\log(d).
%\end{align*}
%Inserting in the definition of $L_m$ we get

According to \eqref{fourier2} and section 3.3 of \cite{Br1}, 
%for $y_1 y_2>m/D$ 
we have
\begin{align} 
\nonumber
G_m(z_1,z_2)&= - \frac{\varphi_m(1)}{2} \log(16 \pi^2 y_1 y_2)\\
\label{cuspest}
&\phantom{=}{}-\log \Big| \e(\rho_W z_1 + \rho'_W z_2)
  \prod_{\substack{\nu\in\frakd^{-1} \\ (\nu,W)>0\\\norm(\nu)=-m/D}} \left(1-\e(\nu z_1
    +\nu' z_2)\right) \Big|  +o(z_1, z_2).
\end{align} 
Here $W$ is a Weyl chamber of discriminant $m$ and 
$\rho_W$ the corresponding Weyl vector. Moreover, $o(z_1,z_2)$ is a $\Gamma_{K,\infty}$-invariant function on $\H^2$, which defines a smooth function in the neighborhood $V_{m/D}$ of $\infty$ and vanishes at $\infty$.
This describes the singularities of $G_m(z_1,z_2)$ near the cusp
$\infty$. Analogous expansions hold at the other cusps.

\begin{remark}\label{rem:localbp}
The function 
\[
\Psi_{m}^\infty(z_1,z_2)=\e(\rho_W z_1 + \rho'_W z_2)
\prod_{\substack{\nu\in\frakd^{-1} \\ (\nu,W)>0\\\norm(\nu)=-m/D}} 
\left(1-\e(\nu z_1
 +\nu' z_2)\right) 
\]
is a local Borcherds product at the cusp $\infty$ in the sense of \cite{BF}.
There exists a positive integer $n$ such that $n\rho_W\in \frakd^{-1}$ (by \eqref{eq:wv} one can take $n=\tr(\eps_0)$). Then $\Psi_{m}^\infty(z_1,z_2)^n$ defines a holomorphic function in a small neighborhood of $\infty$, whose divisor equals the restriction of $nT(m)$. 
Equivalently we may view it as a $\Gamma_{K,\infty}$-invariant function on $\H^2$, whose divisor on $V_{m/D}$ is equal to  $nT(m)$. In particular, this shows 
that $T(m)$ is a $\Q$-Cartier divisor near the cusps.
\end{remark}

\begin{proposition}\label{prop:Gmbasic} 
If $\Gamma\leq \GK$ is a subgroup of finite index, then 
%$\pi^*(G_m(z_1,z_2))$ 
\begin{align}\label{eq:greenobj}
\frakg(m)=(-2\partial\bar\partial \pi^* G_m, \pi^* G_m)
\end{align}
defines a pre-log-log Green object for the divisor 
$\pi^*(T(m))$ on $\widetilde{X}(\Gamma)$.
\end{proposition}

\begin{proof} 
%By \cite{Br1}, Section 3.3 and formula (38),
%\[
%\| 1_{\calO(T(m))}\| = \exp( - G_m)
%\]
We 
define a metric on the line bundle $\calO(T(m))$ on $\widetilde{X}(\Gamma)$ 
by giving the canonical section
$1_{\calO(T(m))}$ the norm 
\[
\| 1_{\calO(T(m))}\| = \exp( - \pi^* G_m).
\]
By \cite{Br1}, Section 3.3 and formula (38), this metric is smooth outside 
$\D_\Gamma$. We now 
show that it is a pre-log singular hermitian metric in the sense of Definition \ref{def:pre-log-log-hlb}. 

We have to 
consider the growth of this metric locally near points $a\in \D_\Gamma$. 
%We show that the conditions of Proposition \ref{prop:302} (iii) hold for
%$X=\widetilde{X}(\Gamma)$, $y=T(m)$, and $g_y=\pi^*(G_m)$.
%
%Using a partition of unity it suffices to check them locally in a small open neighborhood of any given point $a\in\widetilde{X}(\Gamma)$.
%We only have to do this for points belonging to the exceptional divisor $\D_\Gamma$. At the other points it follows from \cite{Br1} Theorem 4 by the usual Poincar\'e-Lelong argument.
Here we only carry out the case that $a$ lies above the cusp $\infty$. The other cusps are treated analogously, and the easier case that $a$ lies above an elliptic fixed point is left to the reader.

So let $a\in \pi^*(\infty)$, and $U\subset\widetilde{X}(\Gamma)$ be a small open neighborhood of $a$ such that $\pi(U)\subset V_C$ with $C>m/D$, as in the discussion preceding Lemma \ref{l2}.   
%After a biholomorphic change of coordinates we may
%further assume that $U=\E^2$ is a product of two unit discs, $a=(0,0)$, $\pi^* \infty=\dv(q_1^\alpha q_2^\beta)$ on $U$, and that $W$ is as in \eqref{eq:W} with $C>m/D$.
%To verify \eqref{eq:basicgreen}, we have to study the singularities of $\pi^*(G_m)$ in $U$.
In view of Remark \ref{rem:localbp}, the function $\pi^*(\Psi_{m}^\infty)$ has precisely the divisor $\pi^*(T(m))$ on $U$. 
Thus $s=1_{\calO(T(m))}/\pi^*(\Psi_{m}^\infty)$ is a trivializing section 
for $\calO(T(m))$ on $U$.
By means of \eqref{cuspest} we find that 
\[
\log\|s\|= \frac{\varphi_m(1)}{2} \pi^*(\log(y_1 y_2))
%-\log| \pi^*(\Psi_{m}^\infty)| 
+  \text{smooth function}.
\]
In the proof of Proposition \ref{prop:petersson} we already saw that $\pi^*(\log(y_1 y_2))$ is a pre-log-log form on $U$. 
Hence the assertion follows from Lemma \ref{lemm:1}.
\end{proof}

\begin{remark}\label{rem:bc}
Observe that $\pi^* T(m)$ may contain components of the exceptional divisor $\D_\Gamma$. This is actually always the case if $T(m)$ is isotropic.
If $\widetilde X(\GK)$ is the Hirzebruch desingularization of
$X(\GK)$, then $\pi^*T(m)$ is equal to the divisor $T_m^c$ considered
by Hirzebruch and Zagier in \cite{HZ}.
\end{remark}

%\begin{definition}\label{def:greenobj}
\begin{notation}
To lighten the notation we will frequently drop the $\pi^*$ from 
the notation. We write
\begin{align}\label{eq:greenobj1}
\frakg(m)=(\omega_m,G_m)=(-2\partial\bar\partial G_m, G_m).
\end{align}
\end{notation}

%\end{definition}

\begin{remark}
The Chern form $\omega_m$ is computed and studied in \cite{Br1} Theorem 7.
It turns out that 
\begin{align}\label{eq:chernform}
\omega_m= 
2\pi i \varphi_m(1) \omega + f(\eps_0 z_1, \eps_0' \bar z_2)\dd z_1 \wedge \dd \bar z_2+ f(\eps_0 z_2, \eps_0' \bar z_1)\dd z_2 \wedge \dd \bar z_1.
\end{align}
Here $f(z_1, z_2)$ is a certain Hilbert cusp form of weight $2$ for $\GK$, essentially the Doi-Naganuma lift of the $m$-th Poincar\'e series in $S_2^+(D,\chi_D)$. Green functions as $G_m$ are investigated from the point of view of the Weil representation in \cite{BrFu}. 
\end{remark}

%%% Local Variables: 
%%% mode: latex
%%% TeX-master: "h-master2"
%%% End: 

% h-star4.tex:\\
\section{Star products on Hilbert modular surfaces}

Here we compute star products on Hilbert modular surfaces related to Hirzebruch-Zagier divisors. Throughout let $\Gamma\leq \GK$ be a subgroup of finite index.

\subsection{Star products for Hirzebruch-Zagier divisors}

In general the product of a mixed growth form as $G_m$ and a pre-log-log form as $\omega^2$ need not be integrable. Therefore the following lemma, which will be crucial for Theorem \ref{thm:123star} below, is special for Hilbert modular surfaces. It seems to be related to the Koecher principle.

\begin{lemma}\label{lem:intganz}
The Green function $G_m(z_1,z_2)$
is in $L^1(\widetilde{X}(\Gamma),\omega^2)$.
%, i.e., the integral $\int_{\widetilde{X}(\Gamma)} |G_m(z_1,z_2)|\,\omega^2$ converges.
\end{lemma}

\begin{proof}
By possibly replacing $\Gamma$ by a torsion free subgroup of finite index, we may assume that $\Gamma\bs\H^2$ is regular. Since $\widetilde{X}(\Gamma)$ is compact, it suffices to show that $G_m$ is locally integrable in a neighborhood of any point of $\widetilde{X}(\Gamma)$. Outside the exceptional divisor $\D_\Gamma$ this easily follows from the fact that $G_m$ has only logarithmic singularities along $T(m)$ on $\Gamma\bs \H^2$.
Hence we only have to show local integrability at $\D_\Gamma$. For simplicity, here we only treat the points above the cusp $\infty$, for the other cusps one can argue analogously.

So let $a\in \pi^*(\infty)$ be a point above $\infty$.
Assuming the notation of the discussion preceding Lemma \ref{l2}, it suffices to prove that  
\begin{align}\label{eq:showconv}
\int\limits_{\E(\delta)^2} | \pi^*(G_m)| \,\pi^*(\omega^2)
\end{align}
converges for some $\delta>0$.

By possibly replacing $\widetilde{X}(\Gamma)$ by a desingularization of
$\widetilde{X}(\Gamma)$ with respect to $T(m)$ we may assume that
$\pi^*(\infty)$, $\pi^* T(m)$, and 
$\pi^*(\infty)\cup\pi^* T(m)$ are divisors with normal crossings supported on the coordinate axes of $\E^2$.
By virtue of \eqref{cuspest} we have
\[
\pi^*(G_m) \prec - \frac{\varphi_m(1)}{2} \pi^*(\log(y_1 y_2))
-\log|\pi^*(\Psi_{m}^\infty)|
\]
on $\E(\delta)^2$. According to Remark \ref{rem:localbp} the function $\pi^*(\Psi_{m}^\infty)$ is holomorphic on $\E^2$ with divisor $\pi^* T(m)$. 
Consequently, if we introduce polar coordinates
\[
q_j=r_j e^{i\rho_j} \qquad \text{($0\leq r_j <1$, $0\leq\rho_j < 2\pi$)}
\]
on $\E^2$, we obtain by means of \eqref{eq:omegasing1}:
\begin{align}\label{est1}
\pi^*(G_m) \prec \log (r_1 r_2).
\end{align}

We now estimate $\pi^*(\omega^2)$. As $H(q_1,q_2)$ is bounded near $(0,0)$, it follows from \eqref{eq:omegasing0} that there is a small $\delta>0$ such that
\begin{align*}
|\mu_1 \log r_1 + \mu_2 \log r_2 |&\cdot |\mu_1' \log r_1 + \mu_2' \log r_2 |\\
&\prec \pi^*(y_1 y_2) \prec |\mu_1 \log r_1 + \mu_2 \log r_2 |\cdot |\mu_1' \log r_1 + \mu_2' \log r_2 |
\end{align*}
on $\E(\delta)$. 
In view of Lemma \ref{l2}, the complex Jacobi matrix of $\pi$ at $(q_1,q_2)$ is given by
\[
J(\pi,q_1,q_2)=\frac{1}{2\pi i} \begin{pmatrix}
\mu_1/q_1 & \mu_2 / q_2\\
\mu_1'/q_1 & \mu_2' / q_2
                                \end{pmatrix}
+J(H,q_1,q_2),
\]
where $J(H,q_1,q_2)$ is smooth.

If $\mu_1$ and $\mu_2$ are both totally positive, we may infer that
$J(\pi,q_1,q_2)\prec \frac{1}{r_1 r_2}$.
Consequently, on $\E(\delta)^2$ we have
\begin{align}\nonumber
\pi^*(\omega^2)&=\pi^* \left(\frac{ 1}{8 \pi^2}
\frac{\dd x_1  \dd y_1}{y_1^2} \,\frac{\dd x_2  \dd y_2}{y_2^2}\right)\\
\label{est2}
&\prec \frac{ \dd r_1 \dd r_2 \dd\rho_1\dd \rho_2}
{r_1 r_2 (\mu_1 \log r_1 + \mu_2 \log r_2)^2 (\mu_1' \log r_1 + \mu_2' \log r_2)^2}.
\end{align}
Combining the above estimates \eqref{est1} and \eqref{est2}, we see that
\[
\int\limits_{\E(\delta)^2} | \pi^*(G_m)| \,\pi^*(\omega^2)
\prec 
\int\limits_{\E(\delta)^2} \frac{|\log(r_1 r_2)| \dd r_1 \dd r_2 \dd\rho_1\dd \rho_2}{r_1 r_2 (\mu_1 \log r_1 + \mu_2 \log r_2)^2 (\mu_1' \log r_1 + \mu_2' \log r_2)^2}.
\]
In order to prove the convergence of the latter integral, it suffices to show that
\begin{align}\label{eq:intI}
\int\limits_0^\delta \int\limits_0^\delta 
\frac{|\log r_1|  \dd r_1  \dd r_2}{r_1 r_2
(\mu_1\log r_1+\mu_2\log r_2)^2 (\mu_1'\log r_1+\mu_2'\log r_2)^2}
\end{align}
converges. If  $0\leq \eps\leq 1$, then the inequality between geometric and arithmetic mean implies that $A^{1+\eps} B^{1-\eps}\prec (A+B)^2$ uniformly for $A,B>0$.
Taking into account that $\mu_1$ and $\mu_2$ are totally positive, we derive
that
\begin{align*}
\text{\eqref{eq:intI}}&\prec \int\limits_0^\delta \int\limits_0^\delta 
\frac{|\log r_1| \dd r_1\dd r_2}{r_1 r_2
|\log r_1|^{2+2\eps} |\log r_2|^{2-2\eps}}
=  \int\limits_0^\delta \frac{\dd r_1}{r_1
|\log r_1|^{1+2\eps}} 
\int\limits_0^\delta 
\frac{\dd r_2}{r_2 |\log r_2|^{2-2\eps} }.
\end{align*}
Since the latter integrals are clearly finite for $0< \eps<1/2$, we obtain the assertion.

If $\mu_1$ and $\mu_2$ are not both totally positive, then without any restriction we may assume that $\mu_1\gg 0$ and $\mu_2 =0$. Then $(\log r_1)^2\prec \pi^*(y_1 y_2) \prec (\log r_1)^2$ and $J(\pi,q_1,q_2)\prec \frac{1}{r_1}$, and the convergence of \eqref{eq:showconv} is immediate.
\end{proof}

\begin{lemma} \label{lem:intrand}
%Let $\pi: \widetilde{X}(\Gamma)\to X(\Gamma)$ be a 
%desingularization of $X(\Gamma)$.  
Let $T(m_1)$, $T(m_2)$, $T(m_3)$ be
Hirzebruch-Zagier divisors, and $\kappa\in X(\Gamma)$ 
be a cusp or an elliptic fixed point. If $\kappa$ is a cusp, then assume in addition that $T(m_2)$ is anisotropic. 
If $B_\varepsilon(\kappa)$ denotes an $\eps$-neighborhood
of $\pi^*(\kappa) \subset \widetilde{X}(\Gamma)$,
then
\begin{align} 
&\lim_{\varepsilon \to 0}
\int\limits_{\partial \left( B_{\varepsilon}(\kappa)  \right)}
     G_{m_1}(z_1,z_2)   \dd^c G_{m_2}(z_1,z_2)\wedge \omega_{m_3}
     =0 \label{eq:intboundFG},\\
&\lim_{\varepsilon \to 0}
\int\limits_{\partial \left( B_{\varepsilon}(\kappa)  \right)}
     G_{m_2}(z_1,z_2)   \dd^c G_{m_1}(z_1,z_2)\wedge \omega_{m_3}
     =0 \label{eq:intboundGF}.
\end{align}
\end{lemma}

\begin{proof}
We only prove the assertion in the case that $\kappa$ is a cusp, leaving the easier other case to the reader.

  By possibly replacing $\widetilde{X}(\Gamma)$ by a desingularization of
  $\widetilde{X}(\Gamma)$ with respect to $T(m_1)$ we may assume that
  $\pi^*(\infty)$, $\pi^* T(m_1)$, and $\pi^*(\infty)\cup\pi^*T(m_1)$ 
  are divisors with normal crossings.
  Since $\widetilde{X}(\Gamma)$ is compact, it suffices to show that
  \eqref{eq:intboundFG} and \eqref{eq:intboundGF} hold locally. 
  We only do this for the cusp $\infty$, at the other cusps 
  one can argue analogously.
  Let $a\in \pi^*(\infty)$ be a point on the exceptional
  divisor over $\infty$, and let $U\subset \widetilde{X}(\Gamma)$ be a 
 small open neighborhood of $a$ such that
  $(\pi^* T(m_2))\cap U=\emptyset$.
  After a biholomorphic change of coordinates we may assume that
  $U=\E^2$, $a=(0,0)$, and $\pi^* T(m_1)=\dv(q_1^\alpha q_2^\beta)$ on
  $U$.  We assume the notation of the proof of Lemma \ref{lem:intganz}.
  
Without loss of generality it suffices to show that for some $1\geq \delta>0$:
\begin{align}\label{eq:sts1} 
&\lim_{\varepsilon \to 0}
\int\limits_{\E(\delta)\times\partial \E(\eps)}
     |\pi^*(G_{m_1})   \dd^c \pi^*( G_{m_2})\wedge \pi^*(\omega_{m_3})|
     =0 ,\\
\label{eq:sts2} 
&\lim_{\varepsilon \to 0}
\int\limits_{\E(\delta)\times\partial \E(\eps)}
     |\pi^*(G_{m_2})  \dd^c \pi^*( G_{m_1})\wedge \pi^*(\omega_{m_3})|
     =0 .
\end{align}
We use the local expansions of $G_{m_1}$ and $G_{m_2}$ at the cusp 
$\infty$ given in 
\eqref{cuspest}. On $\E^2$ we have
\begin{align*}
\pi^*(G_{m_1}) &= - \frac{\varphi_{m_1}(1)}{2} \pi^*(\log(y_1 y_2))
-\log|q_1^\alpha q_2^\beta|+\text{smooth function},\\
\pi^*(G_{m_2}) &= - \frac{\varphi_{m_2}(1)}{2} \pi^*(\log(y_1 y_2))+\text{smooth function}.
\end{align*}
Here we have used that $T(m_2)$ is anisotropic, which implies that $\Psi_{m_2}^\infty=1$.
Furthermore, by means of \eqref{eq:chernform} we find
%For $\omega_{m_3}$ one easily infers from \cite{Br1} Theorem 7 that 
\begin{align*}
\pi^*(\omega_{m_3})=2\pi i \varphi_{m_3}(1)\pi^*(\omega)+\text{smooth differential form}.
\end{align*}

We now estimate the quantities occurring in the above boundary integrals in polar coordinates. We recall from \eqref{eq:omegasing1} that
\begin{align*}%\label{eq:omegasing11}
\pi^*(\log(y_1 y_2)) &= \log\left( -\tfrac{1}{4\pi}g_1(q_1,q_2)\right)
+\log\left(-\tfrac{1}{4\pi} g_2(q_1,q_2)\right),
\end{align*}
where 
\begin{align*}
g_1(q_1,q_2)= \mu_1\log |q_1|^2+\mu_2\log |q_2|^2 
-4\pi \Im  H_1(q_1,q_2),\\ 
g_2(q_1,q_2)= \mu_1'\log |q_1|^2+\mu_2'\log |q_2|^2 
-4\pi\Im  H_2(q_1,q_2).
\end{align*}
There is a small $1>\delta>0$ such that
\begin{align*}
\mu_1 \log r_1 + \mu_2 \log r_2 \prec g_j(q_1,q_2) \prec \mu_1 \log r_1 + \mu_2 \log r_2,
\end{align*}
for $j=1,2$, 
%and $\pi^*(\log(y_1 y_2))\prec \log|\log r_1| + \log|\log r_2|$
on $\E(\delta^2)$.
It follows that
\begin{align*}
\pi^*(G_{m_1}) &\prec \log r_1 + \log r_2,\\
\pi^*(G_{m_2}) &\prec \log|\log r_1| + \log|\log r_2|.
\end{align*}
Moreover,
\begin{align*}
\dd^c \pi^*(G_{m_2}) &\prec r_1 \dd \rho_1 + r_2 \dd \rho_2+\dd r_1 + \dd r_2+ \dd^c\pi^*(\log(y_1 y_2)),\\
&\prec r_1 \dd \rho_1 + r_2 \dd \rho_2+\dd r_1 + \dd r_2 + \frac{\dd^c g_1}{g_1}+ \frac{\dd^c g_2}{g_2}\\
&\prec r_1 \dd \rho_1 + r_2 \dd \rho_2+\dd r_1 + \dd r_2 + \frac{\mu_1 \dd \rho_1+\mu_2 \dd \rho_2+ \dd r_1 + \dd r_2}{|\mu_1 \log r_1 + \mu_2 \log r_2|}.
\end{align*}
Here ``$\prec$'' is understood component-wise in the present coordinates. In the same way we get
\begin{align*}
\dd^c \pi^*(G_{m_1}) &\prec \dd \rho_1 + \dd \rho_2+\dd r_1 + \dd r_2+ \dd^c\pi^*(\log(y_1 y_2)),\\
&\prec \dd \rho_1 + \dd \rho_2+\dd r_1 + \dd r_2 +\frac{\mu_1 \dd \rho_1+\mu_2 \dd \rho_2+ \dd r_1 + \dd r_2}{|\mu_1 \log r_1 + \mu_2 \log r_2|}.
\end{align*}
Finally, we find 
\begin{align*}
\pi^*(\omega) &=-\ddc \pi^*(\log(y_1 y_2))\\
&=\frac{\dd g_1 \dd^c g_1}{g_1^2} + \frac{\dd g_2 \dd^c g_2}{g_2^2}
 +\frac{4\pi\ddc \Im H_1}{g_1}    +\frac{4\pi \ddc \Im H_2}{g_2}\\
&\prec \frac{ ( \mu_1 \tfrac{\dd r_1}{r_1} + \mu_2   \tfrac{\dd r_2}{r_2} + r_1 \dd \rho_1 + r_2 \dd \rho_2)  
( \mu_1 \dd \rho_1 + \mu_2 \dd \rho_2 + \dd r_1 +\dd r_2 )  }{(\mu_1 \log r_1 + \mu_2 \log r_2)^2}  \\
&\phantom{\prec}{}+
\frac{(r_1 \dd \rho_1 + r_2 \dd \rho_2+ \dd r_1 + \dd r_2)^2}{|\mu_1 \log r_1 + \mu_2 \log r_2|} . 
\end{align*}

We now estimate the integrals in \eqref{eq:sts1} and \eqref{eq:sts2}. We only consider the case that both, $\mu_1$ and $\mu_2$, are non-zero, leaving the easier case that one of them vanishes to the reader. Only the $\dd r_1 \dd \rho_1 \dd \rho_2$ component of the integrand gives a non-zero contribution.
After a calculation, we find that in \eqref{eq:sts1} this component is bounded by
\[
\frac{\dd r_1 \dd \rho_1 \dd \rho_2}{r_1 |\log r_1 + \log r_2|^2 } + \frac{r_1+r_2}{|\log r_1 + \log r_2|} \dd r_1 \dd \rho_1 \dd \rho_2.
\]
In \eqref{eq:sts2} it is bounded by the same quantity times  $\log|\log r_1| + \log|\log r_2|$. In both cases this implies that the integrals are $O(1/\sqrt{\log(\eps)})$ as $\eps \to 0$.  
\end{proof}

\begin{theorem} 
\label{thm:123star}
Let $T(m_1)$, $T(m_2)$, and $T(m_3)$ be Hirzebruch-Zagier divisors, such that all possible intersections on $X(\Gamma)$ among them are proper, and such that $T(m_2)$ is anisotropic.
Then 
%If $\mathfrak{g}(m_j)= (\omega_{m_j}, G_{m_j})$ is the Green object
%  associated with $\widetilde{T}(m_j)$, then
\begin{align} 
&  \frac{1}{(2\pi i)^2}\int\limits_{\widetilde{X}(\Gamma)} 
\frakg(m_1)*\frakg(m_2)*\frakg(m_3) \notag\\
&= \frac{1}{(2\pi i)^2} \int\limits_{\Gamma \setminus \H^2} 
G_{m_1}  \, \omega_{m_2} \land \omega_{m_3} +
\frac{1}{2\pi i}  \int\limits_{T(m_1)'} G_{m_2} \,\omega_{m_3}
+ \int\limits_{T(m_1)' \cap T(m_2)'} G_{m_3}.
\label{starformula}
\end{align}
Here $T(m_j)'$ denotes the strict transform of the divisor $T(m_j)$ on $X(\Gamma)$.
\end{theorem}

\begin{proof}  Let $\D_1$ be the sum of the connected components of 
$\D_\Gamma$, which have
  non empty intersection with $T(m_1)$ on $\widetilde X(\Gamma)$.  
To ease notation we put $
  \frakg_{2,3} = \frakg(m_2)*\frakg(m_3)$.  We note that
  $\frakg_{2,3}=(-2 \partial \bar \partial g_{2,3} ,g_{2,3})$ is a
  log-log Green object for the zero cycle $ T (m_2) \cap T(m_3)$,
  which has empty intersection with $\D_1$.  Thus the assumptions of
  Theorem \ref{STAR-PRODUCT} are satisfied.  Consequently,
\begin{align*}
  &\frac{1}{(2 \pi i)^2}   \int\limits_{\widetilde{X}(\Gamma)} \frakg(m_1) * \frakg_{2,3}\\
  &=\lim_{\varepsilon \to 0} \left(\frac{1}{(2\pi i)^2}
    \int\limits_{\widetilde{X}(\Gamma) \setminus B_{\varepsilon}(\D_\Gamma)} G_{m_1}
    \wedge -2 \partial \bar \partial g_{2,3} - \frac{2}{2\pi i}
    \int\limits_{\partial \left( B_{\varepsilon}(\D_1) \right)}
    g_{2,3} \wedge \dd^c G_{m_1}
    - G_{m_1} \wedge \dd^c g_{2,3} \right)\\
  &\phantom{=}{}+\frac{1}{ 2 \pi i} \int\limits_{ \overline{
      T(m_1) \setminus ( T(m_1)\cap \D_1)} } g_{2,3}.
\end{align*}
Since $-2 \partial \bar \partial g_{2,3} = \omega_{m_2} \land
\omega_{m_3}$ and $\omega_{m_2} \land \omega_{m_3} \prec \omega^2$, we obtain by Lemma \ref{lem:intganz} for the first
integral in question
\begin{align} \label{prop:sum1}
\lim_{\varepsilon \to 0}
\left(
 \int_{\widetilde{X}(\Gamma) \setminus B_{\varepsilon}(\D_\Gamma)} 
    G_{m_1} \wedge -2 \partial \bar \partial   g_{2,3} \right) 
 = \int_{\Gamma \setminus \ohh^2}
G_{m_1}  \omega_{m_2} \land \omega_{m_2}.  
\end{align}

We now show that the integrals along the boundary vanish in the limit.
Recall that by \eqref{eq:starprodformel} a representative of 
$\mathfrak{g}_{2,3}$
%\in\widehat{H}^4_{\mathfrak{D},Z}(\widetilde{X}(\Gamma),2) 
is given by
the pair
\begin{align*}
 \big( \omega_{m_2} \land \omega_{m_3},  \, 
\sigma_{_{3,2}} G_{m_2} 
\cdot  \omega_{m_3}  -2 \partial \bar \partial 
  \left(\sigma_{_{2,3}} G_{m_2}\right) 
\cdot G_{m_3}\big);
\end{align*} 
here $\sigma_{_{2,3}}:=\sigma_{T(m_2),T(m_3)}$ and $\sigma_{_{3,2}}$ form a partition of unity
as in \eqref{eq:starprodformel1}. Because $T(m_2)\cap T(m_3)\cap E_j=\emptyset$ for all irreducible components $E_j$ of $\D_1$, we may choose $\sigma_{_{2,3}}$ so that
\begin{align*}
g_{2,3}\big|_{\partial (B_\epsilon(E_j))} = 
\begin{cases}  
   G_{m_3}  \omega_{m_2}, & \text{if $T(m_2) \cap E_j
  \neq \emptyset$,} \\
G_{m_2} \omega_{m_3}, & \text{if $T(m_2) \cap E_j = \emptyset$.}
\end{cases}
\end{align*} 
In particular, by our assumption on $T(m_2)$, if $E_j$ lies above a cusp, then the second case applies.
Because $\omega_{m_j} \prec
\omega$ for $j=2,3$, we find  by Lemma \ref{lem:intrand} :
\begin{align} \label{prop:sum3}
  \lim_{\varepsilon \to 0} \left( 
\int\limits_{\partial \left( B_{\varepsilon}(\D_1) \right)}
    g_{2,3} \wedge {\rm d}^c G_{m_1} -G_{m_1}
    \wedge {\rm d}^c g_{2,3} \right)=0.
\end{align}
Applying Theorem \ref{STAR-PRODUCT} again to the remaining integral, we
obtain
\begin{align}\label{eq:littlestar}
\frac{1}{2 \pi i}  \int\limits_{
\overline{T(m_1) 
\setminus ( T(m_1)\cap \D_1)} } g_{2,3}
=\frac{1}{2 \pi i} 
 \int\limits_{T(m_1)'} G_{m_2} \,\omega_{m_3}
+ \int\limits_{T(m_1)' \cap T(m_2)'} G_{m_3};
\end{align}
here the vanishing of the boundary terms follows from the
compatibility
of this star product and the formulas for the generalized arithmetic
intersection number at the infinite places in \cite{Kue2}. 
Putting together all the 
terms yields the claim.
\end{proof}

\begin{remark}
i) The formula for the star product in
  the above theorem is formally the same as the usual formula for the
  star product of Green currents considered by Gillet and Soul\'e.
  The main difference is that in the above star product the
  components of the desingularization of  the divisor
  $T(m_1)$  do not contribute. 
In other words this star product is independent of the choice we
  made in the desingularization and only  depends on the Baily-Borel
  compactification $X(\Gamma)$.
  
ii) For other choices of $T(m_2)$ the contribution of \eqref{eq:littlestar}
  has to be calculated by different formulas, which then essentially contain
  the generalized arithmetic intersection number at the
  infinite places given in \cite{Kue2}.
\end{remark}

\begin{remark} \label{rem:subgroupstar}
Since $T(m)$ and $\frakg(m)$ are invariant under the full Hilbert modular group $\GK$, we clearly have
\begin{align*} 
%\frac{1}{(2\pi i)^2}
\int\limits_{\widetilde{X}(\Gamma)} 
\frakg(m_1)*\frakg(m_2)*\frakg(m_3)= [\GK:\Gamma]\int\limits_{\widetilde{X}(\GK)} 
\frakg(m_1)*\frakg(m_2)*\frakg(m_3).
\end{align*}
\end{remark}

\subsection{Integrals of Green functions}

The purpose of this section is to compute the first integral in
(\ref{starformula}) in the case that $\omega_{m_2}=\omega_{m_3}=2\pi i \omega$. 
%We begin by computing an analogous 
%integral over $\Phi_m(z_1,z_2,s)$ in Theorem \ref{intphif} below.
Recall our normalization (\ref{defomega}) of the
invariant K\"ahler form $\omega$. In
addition we put
\[
\eta_1 =\frac{1}{4\pi} \frac{\dd x_1 \dd y_1}{y_1^2},\qquad 
\eta_2 =\frac{1}{4\pi} \frac{\dd x_2 \dd y_2}{y_2^2}.
\]

%First, we describe the Hirzebruch-Zagier divisor $T(m)$ on
%$\GK\bs\H^2$ more precisely (see \cite{Za1} \S1, \cite{Ge} Chapter V.1).

We consider the quadratic space $V$ and the lattices $L=L(\OK)$,
$L'=L'(\OK)$ defined at the beginning of Section \ref{sect:2.3}. We
put $W=\kzxz{0}{-1}{1}{0}$.
%\cite{Za2} \S6).  
The Hilbert
modular group $\GK$ acts on $L'$ by $\gamma.A= \gamma A {\gamma'}^t$ for $\gamma\in\GK$.
For every $A=\kzxz{a}{\nu}{\nu'}{b}\in L'$
the graph
$\{ (WAz,z);\; z\in \H\}$ defines a divisor in $\H^2$ given by the equation 
$az_1z_2+\nu z_1+\nu' z_2+b=0$.  
The stabilizer
$\Gamma_{K,A}$ of $A$ under the action of $\GK$ can be viewed as an
arithmetic subgroup of $\Sl_2(\R)$ with finite covolume (see \cite{Za1} \S1, \cite{Ge} Chapter V.1). 
By reduction theory, the subset
\[
L'_m=\{ A\in L';\quad \det A = m/D\}
\]
of elements with norm $m/D$ decomposes into finitely many
$\Gamma_K$-orbits. The divisor $T(m)$ on $\GK\bs\H^2$ is given by
\begin{equation}\label{tmisom}
T(m)\cong\bigcup_{A\in \GK\bs L'_m/\{\pm 1\}} \Gamma_{K,A} \bs \H.
\end{equation}
Here $\{\pm 1\}$ acts on $L'_m$ by scalar multiplication.

We may rewrite $\Phi_m(z_1,z_2,s)$ using this splitting.  
For $A=\kzxz{a}{\nu}{\nu'}{b}\in L'_m$ we define
\begin{equation}\label{defga}
d_A(z_1,z_2)=1+\frac{|z_1-WAz_2|^2}{2\Im( z_1) \Im (WAz_2)}=1+\frac{ |az_1 z_2 +\nu z_1+\nu' z_2 +b|^2}{ 2y_1 y_2 m /D }.
\end{equation}
If $\gamma\in \GK$, we have
\[
d_{\gamma.A}(z_1,z_2)=d_A(\gamma^t z_1,\gamma'{}^t z_2).
\]
This follows from the fact that
$\frac{|z_1-z_2|^2}{2y_1 y_2}$ is a point-pair invariant, i.e., only
depends on the hyperbolic distance of $z_1$ and $z_2$.  
%
%Let
%\[
%\calA=\{ A\in M_2(\OK); \quad A^*=A'\}
%\]
%be the lattice of $2\times2$-matrices $A=\kabcd$, whose adjoint
%$A^*=\kzxz{d}{-b}{-c}{a}$ equals the conjugate $A'$.  The Hilbert
%modular group $\GK$ acts on $\calA$ from the right via $M\circ A= M^*
%AM'$ for $M\in \GK$.  For every $A\in \calA$ with $\det A>0$ the graph
%$\{ (Az,z);\; z\in \H\}$ defines a curve in $\H^2$. The stabilizer
%$\Gamma_{K,A}$ of $A$ under the action of $\GK$ can be viewed as an
%arithmetic subgroup of $\Sl_2(\R)$ with finite covolume.  The subset
%\[
%\calA_m=\{ A\in \calA;\quad \det A = m\}
%\]
%of elements with determinant $m$ decomposes into finitely many
%$\Gamma_K$-orbits. The divisor $T(m)$ on $\GK\bs\H^2$ is given by
%\begin{equation}\label{tmisom}
%T(m)\cong\bigcup_{A\in \GK\bs\calA_m/\{\pm 1\}} \Gamma_{K,A} \bs \H.
%\end{equation}
%Here $\{\pm 1\}$ acts on $\calA_m$ by scalar multiplication.
%
%We may rewrite $\Phi_m(z_1,z_2,s)$ using this splitting.  A typical
%matrix in $\calA_m$ has the form
%$A=\kzxz{\theta'}{-b\sqrt{D}}{a\sqrt{D}}{\theta}$ with $a,b\in \Z$ and
%$\theta\in \OK$.  It is easily checked that the function
%\begin{equation}\label{defga}
%d_A(z_1,z_2)=1+\frac{|z_1-Az_2|^2}{2\Im z_1 \Im Az_2}
%\end{equation}
%is equal to
%\[
%1+\frac{ |az_1 z_2 +\lambda z_1+\lambda' z_2 +b|^2}{ 2y_1 y_2 m /D },
%\]
%with $\lambda=\theta/\sqrt{D}$.
%Moreover, $d_A(z_1,z_2)$ satisfies the identity
%\[
%d_A(Mz_1,M'z_2)=d_{M^*AM'}(z_1,z_2)
%\]
%for $M\in \GK$. This follows from the fact that
%$\frac{|z_1-z_2|^2}{2y_1 y_2}$ is a point-pair invariant, i.e., only
%depends on the hyperbolic distance of $z_1$ and $z_2$.  
%
Consequently,
\begin{equation}\label{phirew}
\Phi_m(z_1,z_2,s) = \sum_{A\in L'_m} Q_{s-1}\left( d_A(z_1,z_2)\right)
=\sum_{A\in \GK\bs L'_m}\sum_{\gamma\in \Gamma_{K,A}\bs\GK}
Q_{s-1}\left( d_A(\gamma z_1, \gamma' z_2)\right).
\end{equation}

\begin{proposition}\label{intprop}
The integral
\[
\int_{\GK\bs\H^2} | \Phi_m(z_1,z_2,s) | \, \omega^2
\]
converges for all $s\in \C$ with $\Re(s) >1$.
\end{proposition}

\begin{proof}
  According to (\ref{phirew}) we have the formal identity
\begin{align}
\nonumber
\int\limits_{\GK\bs\H^2} | \Phi_m(z_1,z_2,s) | \, \omega^2
&= \sum_{A\in \GK\bs L'_m\;}\int\limits_{\GK\bs\H^2}
\sum_{\gamma\in \Gamma_{K,A}\bs\GK} | Q_{s-1}\left( d_A(\gamma z_1,\gamma' z_2)\right) |\,\omega^2\\
\nonumber
&= \sum_{A\in \GK\bs L'_m\;}
\int\limits_{\Gamma_{K,A}\bs\H^2} | Q_{s-1}\left( d_A( z_1, z_2)\right) |\,\omega^2\\
&= 2 \sum_{A\in \GK\bs L'_m\;}\int\limits_{z_2\in \Gamma_{K,A}\bs\H}\,
\int\limits_{z_1\in \H}
 | Q_{s-1}\left( d_A( z_1, z_2)\right) |\, \eta_1\,\eta_2.
\label{intprop1}
\end{align}
By a standard Fubini type lemma on integrals over Poincar\'e series
(see for instance \cite{Fr} AII.7), the integral on the left hand side
converges (and equals the right hand side), if the latter integral on
the right hand side converges.  Thus it suffices to prove that
\begin{equation}\label{intprop2}
\int\limits_{z_2\in \Gamma_{K,A}\bs\H}\,\int\limits_{z_1\in \H}
 | Q_{s-1}\left( d_A( z_1, z_2)\right) |\, \eta_1\,\eta_2 
\end{equation}
converges.  We notice that the inner integral actually does not depend
on $z_2$ and $A$. Using the fact that $\frac{|z_1-z_2|^2}{2y_1 y_2}$
is a point-pair invariant and the invariance of $\eta_1$, we find
that (\ref{intprop2}) is equal to
\[
\int\limits_{z_2\in \Gamma_{K,A}\bs\H}\,\int\limits_{z_1\in \H}
 \left| Q_{s-1}\left( 1+\frac{| z_1- i|^2}{2\Im z_1 \Im i}\right) \right|\, \eta_1\,\eta_2. 
\]
This integral is obviously bounded by
\[
\vol( \Gamma_{K,A}\bs\H) \int\limits_{z_1\in \H}
 \left| Q_{s-1}\left( 1+| z_1- i|^2/2y_1\right) \right|\, \eta_1.
\]
That the latter integral is finite for $\Re(s)>1$ is a well know fact.
See for instance \cite{La} Chapter 14 \S3.
\end{proof}

\begin{lemma}\label{greenh}
  Let $h:\H\to \C$ be a bounded eigenfunction of the hyperbolic
  Laplacian $\Delta_1$ with eigenvalue $\lambda$. Then for $s\in \C$
  with $\Re(s)>1$ we have 
\[
\int\limits_{\H} Q_{s-1}\left( 1+\frac{| z_1- z_2|^2}{2y_1 y_2}\right) h(z_1) 
\,\eta_1 = \frac{1/2}{s(s-1)-\lambda} h(z_2).
\]
\end{lemma}

\begin{proof} This statement is well known. It can be proved using 
the Green formula. (See for instance \cite{Iw} Chapter 1.9. Notice the
different normalization there.)
\end{proof}

\begin{theorem}\label{intphif}
  Let $f:\GK\bs\H^2\to \C$ be a bounded eigenfunction of the Laplacian 
  $\Delta_1$ (or $\Delta_2$) with
  eigenvalue $\lambda$.  Then for $s\in \C$ with $\Re(s)>1$ we have
\[
\int\limits_{\GK\bs \H^2} \Phi_m(z_1,z_2,s) f(z_1,z_2)\, \omega^2 = 
\frac{1}{s(s-1)-\lambda} \int\limits_{T(m)} f(z_1,z_2) \,\omega.
\]
\end{theorem}

\begin{proof}
  Since $\Phi_m(z_1,z_2,s)$ is invariant under $(z_1,z_2)\mapsto
  (z_2,z_1)$, we may assume that $f$ is an eigenfunction of $\Delta_1$
  with eigenvalue $\lambda$.
  
  First, we notice that the integral $I$ on the left hand side
  converges by Proposition \ref{intprop}. Similarly as in the proof of
  Proposition \ref{intprop} we rewrite it as follows:
\begin{align*}
I &= 
\int\limits_{\GK\bs \H^2} 
\sum_{A\in \GK\bs L'_m}\sum_{\gamma\in \Gamma_{K,A}\bs\GK}
Q_{s-1}\left( d_A(\gamma z_1,\gamma' z_2)\right)f(z_1,z_2)\, \omega^2\\
&= 2 \sum_{A\in \GK\bs L'_m}
\int\limits_{z_2\in \Gamma_{K,A}\bs\H}\,\int\limits_{z_1\in \H}
Q_{s-1}\left( 1+\frac{| z_1- WA z_2|^2}{2\Im( z_1) \Im(W Az_2)} \right) f(z_1,z_2)\, \eta_1\, \eta_2.
\end{align*}
Here the inner integral can be computed by means of Lemma \ref{greenh}. We obtain
\begin{align*}
I &= 
\frac{2}{s(s-1)-\lambda}\sum_{A\in \GK\bs L'_m/\{\pm 1\}}
\int\limits_{z_2\in \Gamma_{K,A}\bs\H} f(WA z_2,z_2)\, \eta_2\\
&=\frac{1}{s(s-1)-\lambda}
\int\limits_{T(m)} f(z_1,z_2)\, \omega.
\end{align*}
This concludes the proof of the theorem.
\end{proof}

\begin{corollary}\label{cor:greenint}
i) If $s\in \C$ with $\Re(s)>1$, then
\[
\int\limits_{\GK\bs \H^2} \Phi_m(z_1,z_2,s) \, \omega^2 = 
\frac{2}{s(s-1)}\vol(T(m)). 
\]

ii) The volume of $T(m)$ with respect to the invariant volume form $\frac{\dd x
  \dd y}{4\pi y^2}$ on $\H$  equals
\begin{align}\label{eq:voltm}
\vol(T(m))=\frac{1}{2} \zeta_K(-1) \varphi_m(1)=\frac{1}{24} \sigma_m(-1).
\end{align}

iii) We have 
\begin{equation} \label{intphi1}
\int\limits_{\GK\bs \H^2} G_m(z_1,z_2) \, \omega^2 
= -\vol(T(m))\left(2\frac{L'(-1,\chi_D)}{L(-1,\chi_D)}-
 2\frac{\sigma_m'(-1)}{\sigma_m(-1)}
     +1 +\log(D) \right).
\end{equation}
\end{corollary}

\begin{proof}
If we use Theorem \ref{intphif} with $f=1$, we get
\begin{align*}
\int\limits_{\GK\bs \H^2} \Phi_m(z_1,z_2,s) \, \omega^2 &= 
\frac{1}{s(s-1)} \int\limits_{T(m)}  \,\omega
=\frac{2}{s(s-1)} \vol(T(m))\\
&=\frac{2 \vol(T(m))}{s-1} - 2 \vol(T(m)) + O(s-1)  .
\end{align*}
By \eqref{eq:ZETAVOLUMEN} we have in addition
\begin{align*}
\int\limits_{\GK\bs \H^2} \frac{\varphi_m(1)}{s-1} \, \omega^2 
=\frac{1}{s-1} \zeta_K(-1) \varphi_m(1).
\end{align*}
Since $\Phi_m(z_1,z_2) = \lim_{s \to 1}\left( \Phi_m(z_1,z_2,s) -
  \frac{\varphi_m(1)}{s-1} \right)$ is regular and integrable at $s=1$, we derive the second claim by comparing residues in the latter two equalities.
The third claim follows from Definition \ref{def:Gm} and \eqref{eq:L_m} by comparing the constant terms. 
\end{proof}

Notice that our normalization of $\vol T(m)$ equals $-1/2$ times the
normalization of \cite{HZ} and \cite{Ge} Chapter V.5. In particular,
here the volume of $T(1)\cong \Gamma(1)\bs\H$ is $1/12$.

%\begin{remark} In \cite{Br2} analogous Green
%  functions are considered in a more general setting. They are associated with Heegner divisors on hermitian quotients of the orthogonal group
%$\Orth(2,n)$. Using a similar argument as above, one can
%calculate the corresponding integrals, and should thereby obtain a generalization of recent results
%of Kudla \cite{Ku4}. 
%We hope to come back to this in a subsequent paper. 
%\end{remark} 

\subsection{Star products on isotropic Hirzebruch-Zagier divisors}
\label{sect:4.3}

For the rest of this section we assume that $p$ is a prime which is
split in $\OK$ or $p=1$.  Let $\frakp$ be a prime ideal of $\OK$ above
$p$. There is a fractional ideal $\frakc$ and a totally positive
$\lambda\in K$ such that $\frakp=\lambda\frakc^2$.  We fix a matrix
$M\in \kzxz{ \frakc^{-1}}{ \frakc^{-1}}{\frakc}{ \frakc}\cap
\Sl_2(K)$.

It is well-known that the anisotropic Hirzebruch-Zagier divisor
$T(p)\subset \GK \setminus \H^2$ is irreducible. It may have points of
self intersection, and its normalization is isomorphic to the
non-compact modular curve $Y_0(p)=\Gamma_0(p) \setminus \H$. The
normalization of the closure of $T(p)$ in $X(\GK)$ is isomorphic to
the compact modular curve $X_0(p)$, the standard compactification of
$Y_0(p)$ (cf.~\cite{Ge} Chapter V.1 and \cite{HZ}).  For basic facts
on integral models of $X_0(p)$, the line bundle of modular forms on
it, and the normalization of the corresponding Petersson metric, we
refer to \cite{Kue2}. We now describe how $T(p)$ can be parametrized.

On the Hilbert modular surface $X(\Gamma(\frakp))$ the Hirzebruch-Zagier divisor 
$T_\frakp(p)$ is simply given by the diagonal. 
More precisely, the assignment $\tau\mapsto (\tau,\tau)$ induces a morphism of degree $1$
\begin{align}\label{eq:phidiag}
X_0(p)\longrightarrow X(\Gamma(\frakp)), 
\end{align}
whose image is $T_\frakp(p)$. In fact, the 
vector $\kzxz{0}{p/\sqrt{D}}{-p/\sqrt{D}}{0}\in L'(\frakp)$ has determinant $p^2/D$. Since $T_\frakp(p)$ is irreducible, the image of the diagonal in $X(\Gamma(\frakp))$ is $T_\frakp(p)$. Moreover, it is easily checked that the stabilizer in $\Gamma(\frakp)$ of the diagonal in $\H^2$ is equal to $\Gamma_0(p)$.
The cusp $\infty$ (respectively $0$) of $X_0(p)$ is mapped to the cusp  $\infty$ (respectively $0$) of $X(\Gamma(\frakp))$.

\begin{proposition} \label{prepullerback}
The map $\H\to \H^2$, $\tau \mapsto M^{-1}\kzxz{\lambda}{0}{0}{1} (\tau,\tau)$, induces a morphism
\begin{align}\label{eq:defvarphi}
\varphi: X_0(p) \longrightarrow  X(\GK),
\end{align}
which is generically one-to-one and whose image equals $T(p)$.
The cusp $\infty$ (respectively $0$) of $X_0(p)$ is mapped to the cusp  $\frakc$ (respectively $\frakc^{-1}$) of $X(\GK)$.
\end{proposition}

\begin{proof}
The map $\varphi$ is given by the commutative diagram
\begin{align*}
\xymatrix{
X_0(p)\ar[r] \ar[dr]_{\varphi} & X(\Gamma(\frakp))\\
& X(\GK)\ar[u]
},
\end{align*}
where the horizontal arrow is given by \eqref{eq:defvarphi} and the vertical arrow (which is an isomorphism) by \eqref{eq:aiso}. The properties of the latter two maps imply the assertion.
%
%We consider the curve $C=\{ (\theta z, \theta'z);\;
%z\in \H\}$ in $\H^2$. 
%Its stabilizer in $\GK$ is given by
%\[
%\left\{ \zxz{a}{b\theta}{c/\theta}{d};\quad \abcd\in \Gamma_0(p)\right\}.
%\]
%This implies that $\varphi$ is well-defined. 
%The image of $\GK C$ in $\GK\bs\H^2$ equals $T(p)$, because $C
%\subseteq T(p)$, and $T(p)\subset \GK \setminus \H^2$ is irreducible.
%\end{comment}
\end{proof}

%It is easily checked that $\varphi$  extends to a morphism   $ X_0(p)\to X(\GK)$
%of the compactifications, which maps the two cusps $0,\infty$ of $X_0(p)$ to the cusp $\infty$ of $X(\GK)$. It follows from Chow's lemma that $\varphi$ is algebraic.  

%
It is easily seen that the pull-back of the K\"ahler form $\omega$ equals
 $2 \frac{\dd x\dd y}{4 \pi  y^2}$.

For the next proposition we recall that the Fricke involution
$W_p$ on the space of modular forms of weight $k$ for $\Gamma_0(p)$ 
is defined by $ f(z) \mapsto (f\mid W_p )(z)=p^{k/2} z^{-k} f(-1/pz)$.
Moreover, we define the Petersson slash operator for Hilbert modular forms in
weight $k$ by
\[
(F\mid_k \gamma) (z_1, z_2) = \norm(\det(\gamma))^{k/2}(cz_1+d)^{-k}(c'z_2+d')^{-k}
F(\gamma z_1, \gamma'z_2),
\]
for $\gamma=\kabcd\in \Gl_2^+(K)$.

\begin{proposition} \label{prop:pullerback}
i) 
If $F$ is a Hilbert modular form of weight $k$ for $\GK$,
then its pull-back 
\[
(\varphi^* F)(\tau)=\left(F\mid_k M^{-1}\kzxz{\lambda}{0}{0}{1}\right)(\tau,\tau)
\]
is a modular form of weight $2k$ for the group
$\Gamma_0(p)$.

ii) 
If $F$ is in addition holomorphic and has the Fourier expansions 
\begin{align*}
F\mid_k M^{-1}\kzxz{\lambda}{0}{0}{1}&=
\sum_{\nu\in \frakp\frakd^{-1}}a_\nu\, \e(\nu z_1+\nu' z_2),\\
F\mid_k {M'}^{-1}\kzxz{\lambda'}{0}{0}{1}&=
\sum_{\nu\in \frakp'\frakd^{-1}}b_\nu\, \e(\nu z_1+\nu' z_2),
\end{align*}
at the cusps $\frakc$ and $\frakc'$, respectively, then  
the Fourier expansions of $\varphi^* F$ at $\infty$ and $0$ 
are given by
\begin{align}\label{eq:fourierinf}
\varphi^*F& =\sum_{n=0}^\infty\sum_{\substack{\nu\in \frakp\frakd^{-1}\\\tr(\nu)=n}} a_\nu \,\e(n\tau),\\
(\varphi^*F)\mid W_p & =\sum_{n=0}^\infty\sum_{\substack{\nu\in \frakp'\frakd^{-1}\\\tr(\nu)=n}} b_\nu \,\e(n\tau).
\end{align}

%iii) The modular form $\varphi^* F$ is
%invariant under the Fricke involution $W_p$,
%$W_p=\kzxz{0}{-1}{p}{0}$ on$X_0(p)$ 
%whenever $F$ is symmetric, i.e., $F(z_1,z_2)=F(z_2,z_1)$.

iii) For the Petersson metric we have:
\begin{equation}
\label{eq:metricchange}
 \varphi^*\left(\log\|  F\|_{\Pet}\right)
%=\log\left(p^k\|\varphi^*( F)\|^2_{\Pet}\right).
= \log \left( \|\varphi^* F\|_{\Pet}\right).
\end{equation}
\hfill$\square$
\end{proposition}

\begin{remark}\label{rem:pullerback}
The above proposition implies in particular that Hilbert modular forms with rational Fourier coefficients are
mapped to modular forms for $\Gamma_0(p)$ with rational Fourier
coefficients. 
This shows that $\varphi$ is actually defined over $\Q$ (cf.~Proposition \ref{prop:t-rational}).
\end{remark}

%% An explicit isomorphism is easily obtained as follows. It induces a morphism 
%% \begin{equation}\label{defi}
%% \iota : X_0(p)\longrightarrow T(p)'\subset \widetilde{X(\GK)}
%% \end{equation}
%% onto the strict transform $T(p)'$ of $\widetilde{T(p)}=\pi^*T(p)$.
%% Recall that, if it is clear from the context we 
%% write $T(p)$ instead of $T(p)'$.

%% The restriction of the K\"ahler form $\omega$ is $2 \frac{dxdy}{4 \pi
%%   y^2}$. If $F \in M_k(\GK)$ is a Hilbert modular form of weight $k$,
%% then its pull-back to $T(p)$, denoted by $\iota^* F(z)=F( \theta z,
%% \theta' z)$, is a modular form of weight $2k$ for the group
%% $\Gamma_0(p)$.  In terms of the Fourier expansion of $F$ as in
%% (\ref{FFourier}) we have
%% \begin{align} \label{eq:iotaFinfty}
%% \varphi^* F(z)&=  a_0 +\sum_{\substack{\nu \in {\frakd}^{-1} \\ \nu \gg 0}} 
%% a_\nu \, e\left (\tr(\nu \theta) z\right). 
%% \end{align}
%%   Moreover, it is easily checked that $\varphi^* F$ is
%% invariant under the Fricke involution $W_p=\kzxz{0}{-1}{p}{0}$ on
%% $X_0(p)$ 
%% %, if
%% whenever $F$ is symmetric, i.e.~$F(z_1,z_2)=F(z_2,z_1)$.  Finally,
%% the identity $ \varphi^*(16 \pi^2 y_1 y_2) = p (4 \pi y)^2$ implies
%% that
%% \begin{equation}
%% \label{eq:metricchange}
%%  \varphi^*\left(\log\|  F\|^2_{\Pet}\right)
%% %=\log\left(p^k\|\varphi_p^*( F)\|^2_{\Pet}\right).
%% =k \log( p )+ \log \left( \|\varphi^* F\|^2_{\Pet}\right).
%% \end{equation}

%\subsection{Star products on isotropic Hirzebruch-Zagier  divisors}

We now compute star products of pull-backs of Hilbert
modular forms via $\varphi$.  
%After a brief
%discussion of the pull-back map we give a first formula for the star
%products in question. 
%We also consider the interplay of
%bad reduction and the (modular) extension of divisors to the minimal
%regular model $\mathcal{X}_0(p)$ of $X_0(p)$ over $\Spec \Z$.

\begin{theorem} \label{thm:modulstar} 
Let $F$, $G$ be Hilbert modular forms of weight $k$ with 
%integral 
rational Fourier coefficients. Assume that all possible intersections on $X(\GK)$ of  $T(p)$,  $\dv(F)$,  $\dv(G)$ are proper, and that $F$ does not vanish at the cusps $\frakc$ and $\frakc^{-1}$  of $\GK$.
Then
\begin{align}
  \nonumber &\int_{X_0(p)} \varphi^* \left( \log\|F\|_{\Pet} \right) 2 k
  \frac{\dd x\dd y}{4 \pi y^2}
  +\int_{\dv(\varphi^*F)} \varphi^* \left(\log\|G\|_{\Pet}\right)\\
  \label{modulstar1}
&= (2 k)^2 \vol(T(p))
\left( \frac{\zeta'(-1)}{\zeta(-1)} +\frac{1}{2} 
%+\frac{1}{2} \log(p)
\right)
+\big( \dv(\varphi^*F), \dv(\varphi^* G)
  \big)_{\mathcal{X}_0(p),\fin}.
\end{align}
Here $ \left( \dv(\varphi^* F),\dv( \varphi^* G)
\right)_{\mathcal{X}_0(p),\fin}$ denotes the intersection number at
the finite places on the minimal regular model
$\mathcal{X}_0(p)$ of $X_0(p)$  of the divisors associated with 
the sections of the line
bundle of modular forms corresponding to $\varphi^* F$ and 
$\varphi^* G$ (cf.~\cite{Kue2}).
\end{theorem}

\begin{proof} %We assume that the reader is familiar with the theory of 
 % modular curves and modular forms over $\ZZ$ as used in \cite{Kue2}.
  By assumption $\varphi^* F, \varphi^* G$ are modular forms for
  $\Gamma_0(p)$ of weight $2 k$ with rational Fourier expansions, and
  $\varphi^* F$ does not vanish at the cusps.  They determine sections
  %also denoted by $\varphi^* F$, $\varphi^* G$, 
of the line bundle of
  modular forms of weight $2 k$ on $\mathcal{X}_0(p)$.  Notice that on
  $X_0(p)$ we have $\dv (\varphi^*F) \cap \dv( \varphi^*G)=
  \emptyset$.  Because of \eqref{eq:metricchange} the left hand side
  of \eqref{modulstar1} equals the negative of the formula for the
  generalized arithmetic intersection number $\langle \varphi^* F,
  \varphi^* G \rangle_\infty$ at the infinite place given in Lemma~3.9
  of \cite{Kue2}. Hence we have
\begin{align*}
&\int_{X_0(p)}  \varphi^* \left( \log\|F\|_{\Pet} \right) 
2 k \frac{\dd x\dd y}{4 \pi y^2}
+\int_{\dv(\varphi^* F)} \varphi^* \left(\log\|G\|_{\Pet}\right)
=-\langle \varphi^* F, \varphi^*
G \rangle_\infty.
%+ \vol X_0(p) \cdot  2 k^2      \log(p). 
\end{align*}
%Here the extra $\log p$-term comes from the change of metrics
%\eqref{eq:metricchange} and the fact that the degree of $\dv
%\varphi^*F$ is equal to $2 k \vol X_0(p)$.  
Now the claim follows from
Corollary 6.2 of \cite{Kue2}, i.e.,
\[
(2 k)^2 (p+1)
\left(\frac{1}{2}\zeta(-1)+\zeta'(-1)\right)=
\big(\dv(\varphi^* F),
\dv(\varphi^* G)\big)_{\mathcal{X}_0(p),\fin} +\langle \varphi^* F,
\varphi^* G \rangle_\infty,
\]
using the identity $\vol T(p)= \vol X_0(p)=- \zeta(-1)[\Gamma(1):\Gamma_0(p)]=
(p+1)/12$.
\end{proof}

\subsection{Star products on Hilbert modular surfaces}

We combine the results of the previous sections to compute star
products on Hilbert modular surfaces.

\begin{theorem}
\label{thm:hilbertstar}
Let $p$ be a prime which is split in $\OK$ or $p=1$, and let $F$, $G$ be Hilbert modular forms of weight $k$
with 
%integral 
rational Fourier coefficients. Assume that all possible intersections on $X(\GK)$ of  $T(p)$,  $\dv(F)$,  $\dv(G)$ are proper, and that $F$ does not vanish at any cusp of $\GK$.
Then
%Let $F$ and $G$
%be  Hilbert modular forms of weight $k$ with coprime integral Fourier
%coefficients such that $T(p) \cap \dv(F) \cap \dv(G) = \emptyset$
%and  $F(\kappa)=1$ for all cusps $\kappa$ of $X(\GK)$. Then
\begin{align*} 
&  \frac{1}{(2\pi i)^2} \int\limits_{\widetilde{X}(\Gamma_K)}
  \mathfrak{g}(p)*\mathfrak{g}(F)*\mathfrak{g}(G) \\
&=
- k^2 \vol(T(p)) \left(2\frac{ \zeta_K'(-1)}{ \zeta_K(-1)}
    + 2 \frac{\zeta'(-1)}{\zeta(-1)} +3  
+\log( D) \right)\\
&\phantom{=}{} - k^2 \vol(T(p)) \frac{p-1}{p+1} \log(p) 
 -\big(\div(\varphi^*F), \dv(\varphi^*G)
  \big)_{\mathcal{X}_0(p),\fin}. 
\end{align*}
\end{theorem}

\begin{proof} The logarithm of the Petersson norm
  of a Hilbert modular form satisfies along $\D_\Gamma$ the same
  bounds as the Green functions $G_m$.  Hence we may calculate the star
  product in question by means of the formula of Theorem \ref{thm:123star}:
\begin{align*} 
& \frac{1}{(2\pi i)^2} \int\limits_{\widetilde{X}(\Gamma_K)} 
\mathfrak{g}(p)*\mathfrak{g}(F)*\mathfrak{g}(G) \notag\\
&= k^2 \int\limits_{\GK \setminus \H^2} 
G_p \, \omega^2  
+ k \int\limits_{T(p)'} -\log\|F\|_{\Pet} \,\omega
+ \int\limits_{T(p)' \cap \dv(F)'} -\log\|G\|_{\Pet}.
\end{align*}
By Corollary \ref{cor:greenint} (iii) and \eqref{eq:L_p} the first integral is given by 
\begin{align*}
  k^2 \int\limits_{\GK \setminus \H^2} G_p \, \omega^2 &= -k^2
  \vol(T(p)) \left(2 \frac{L'(-1,\chi_D)}{L(-1,\chi_D)} +
    \frac{p-1}{p+1}\log(p) +1 + \log(D) \right) .
\end{align*}
Here we have used that  $\chi_D(p)=1$.    For the remaining
integrals we use the morphism $\varphi$ defined by
\eqref{eq:defvarphi} to infer
\begin{align*}
&k  \int\limits_{T(p)'} -\log\|F\|_{\Pet} \,\omega
+ \int\limits_{T(p)' \cap \dv(F)'} -\log\|G\|_{\Pet}\\
&= - 2 k  \int\limits_{X_0(p)}\varphi^*( \log\|F\|_{\Pet})
 \,\frac{\dd x\dd y}{4\pi y^2}
- \int\limits_{ \dv(\varphi^*F) } \varphi^* \log\|G\|_{\Pet}\\
&=-  k^2 \vol(T(p))  \left( 4
    \frac{\zeta'(-1)}{\zeta(-1)} +2 
%+2 \log(p)
  \right)
 - \big(\dv(\varphi^*F),\dv(\varphi^*G)
  \big)_{\mathcal{X}_0(p),\fin} ,
\end{align*} 
where the last equality was derived by means
of Theorem \ref{thm:modulstar}. 
 Adding the above expressions, the claim follows by the
identity
$\zeta_K(s) = \zeta(s) L(s,\chi_D)$.
\end{proof}

%\begin{remark}
%If the forms $F$ and $G$  in Theorem \ref{thm:hilbertstar} are in
%addition symmetric and have coprime integral Fourier coefficients, then one can use the explicit description of the
%minimal regular model of $X_0(p)$ due to Deligne and Rapoport
%\cite{DeRa}
%to show that 
%\begin{align*}
%&\big(\div(\varphi^*F), \dv(\varphi^*G)
%  \big)_{\mathcal{X}_0(p),\fin} \\
%&= 
%\big(\div(\varphi^*F)^{hor}, \dv(\varphi^*G)^{hor}
%  \big)_{\mathcal{X}_0(p),\fin} -
%k^2 \vol(T(p)) \frac{p-1}{p+1} \log(p).
%\end{align*}
%Here $\div(\varphi^*F)^{hor}$ and $\div(\varphi^*G)^{hor}$ denote the
%Zariski closure of the divisors restricted to the generic fiber.
%\end{remark}

%%% Local Variables: 
%%% mode: latex
%%% TeX-master: "h-master2"
%%% TeX-region: "h-star3-region"
%%% End: 

% h-borcherds3.tex:\\
\section{Borcherds products on Hilbert modular surfaces}

It was shown in section 4 of \cite{Br1} that for certain integral
linear combinations of the $G_m(z_1,z_2)$ all Fourier coefficients,
whose index $\nu\in \frakd^{-1}$ has negative norm, vanish.  Such a
linear combination is then the logarithm of the Petersson metric of a
Hilbert modular form, which has a Borcherds product expansion.  We now
explain this in more detail.

\subsection{Basic properties of Borcherds products}

Recall our assumption that $D$ be a prime.  Let $k$ be an integer. We
denote by $A_k(D,\chi_D)$ the space of {\em weakly holomorphic}
modular forms of weight $k$ with character $\chi_D$ for the group
$\Gamma_0(D)$. These are holomorphic functions $f:\H\to \C$, which
satisfy the transformation law
\[
f(M \tau)=\chi_D(d)(c\tau + d)^k f(\tau)
\]
for all $M=\kabcd\in \Gamma_0(D)$, and are meromorphic at the cusps of
$\Gamma_0(D)$.  If $f=\sum_{n\in \Z} c(n) q^n\in A_k(D,\chi_D)$, then the
Fourier polynomial
\[
\sum_{\substack{n<0}} c(n) q^n
\]
is called the {\em principal part} of $f$.  
Here
$q=e^{2\pi i \tau}$ as usual.
We write $M_k(D,\chi_D)$ (respectively $S_k(D,\chi_D)$) for
the subspace of holomorphic modular forms (respectively cusp forms).

For $\epsilon\in \{\pm 1\}$ we let $A^\epsilon_k(D,\chi_D)$ be the
subspace of all $f=\sum_{n\in \Z} c(n) q^n$ in $A_k(D,\chi_D)$, for
which $c(n)=0$ if $\chi_D(n)=-\epsilon$ (cf.~\cite{BB}). 
A classical Lemma due to Hecke implies
that
\begin{equation}\label{split}
A_k(D,\chi_D)=A^+_k(D,\chi_D)\oplus A^-_k(D,\chi_D).
\end{equation}
Finally, we define the
spaces $M_k^\epsilon(D,\chi_D)$ and $S^\epsilon_k(D,\chi_D)$
analogously.

Here we mainly consider $M_2^+(D,\chi_D)$ and $A^+_0(D,\chi_D)$. The
Eisenstein series
\begin{equation}\label{eis}
E(\tau) = 1+\sum_{n\geq 1} B_D(n)q^n=  1+\frac{2}{L(-1,\chi_D)} \sum_{n\geq 1} 
%\sum_{d\mid n} d \left( \chi_D(d) +  \chi_D(n/d)\right) 
\sigma_n(-1)q^n
\end{equation}
is a special element of $M_{2}^{+}(D,\chi_D)$.  Note that by \eqref{eq:voltm}
\begin{align}\label{eq:frankehausmann}
B_D(n)=-2\varphi_n(1)=-\frac{4}{\zeta_K(-1)}\vol(T(n)).
\end{align}  
The space $M_2^+(D,\chi_D)$ is the
orthogonal sum of $\C E$ and the subspace of cusp forms $
S_{2}^{+}(D,\chi_D)$.  

The existence of weakly holomorphic modular forms in $A^+_0(D,\chi_D)$
with prescribed principal part is dictated by $S_2^+(D,\chi_D)$.
Before making this more precise, it is convenient to introduce the following notation. If $\sum_{n\in\Z}c(n)q^n\in \C\{q\}$ is a formal Laurent series, we put
\begin{equation}\label{ctilde}
\tilde c(n)=\begin{cases} 2c(n),& \text{if $n \equiv 0\pmod{D}$},\\
c(n),& \text{if $n \not\equiv 0\pmod{D}$}.
\end{cases}
\end{equation}
We now recall a theorem from \cite{BB}, 
which is a reformulation of Theorem 3.1 in \cite{Bo3}:

\begin{theorem}\label{crit}
  There exists a weakly holomorphic modular form $f\in
  A_0^+(D,\chi_D)$ with prescribed principal part $\sum_{n<0}c(n)q^n$
  (where $c(n)=0$ if $\chi_D(n)=-1$), if and only if
\begin{equation}\label{cond}
\sum_{n< 0} \tilde c(n)b(-n)=0
\end{equation}
for every cusp form $g=\sum_{m>0} b(m)q^m$ in $S_2^{+}(D,\chi_D)$.
Then the constant term $c(0)$ of $f$ is given by the coefficients of
the Eisenstein series $E$:
\begin{align*}
c(0)&=- \frac{1}{2}\sum_{n<0} \tilde c(n) B_D(-n)=\frac{2}{\zeta_K(-1)}\sum_{n<0} \tilde c(n) \vol(T(-n)).
%\\&=  \sum_{\substack{n<0 \\ n \not\equiv 0 \;(D)}} c(n) \varphi_{-n}(1) + 
% 2\sum_{\substack{n<0 \\ n \equiv 0 \;(D)}} c(n) \varphi_{-n}(1).
\end{align*}
\hfill $\square$
\end{theorem}

Using orthogonality relations for the non-degenerate bilinear pairing  
between Fourier polynomials in $\C[q^{-1}]$ and formal
 power series in $\C[[q]]$ given by
\begin{equation*}
\left\{\sum_{n\le 0} c(n)q^n,\sum_{m \ge0} b(m)q^m\right\}= \sum_{n\leq 0} \tilde c(n)b(-n),
\end{equation*}
one can deduce:

%By extending the bilinear pairing defined by \eqref{cond} to a non-degenerate pairing between Fourier polynomials in $\C[q^{-1}]$ and formal power series in $\C[[q]]$, and using orthogonality relations, one easily deduces:

\begin{corollary}\label{critcor}
A formal power series $\sum_{n\geq 0} b(n)q^n \in\C[[q]]$ (where $b(n)=0$ if $\chi_D(n)=-1$) is a modular form in $M_2^{+}(D,\chi_D)$, if and only if 
\[
\sum_{n\leq 0}\tilde c(n) b(-n)=0
\]
for every $f=\sum_{n\gg-\infty} c(n)q^n$ in $A_0^{+}(D,\chi_D)$.
\hfill $\square$
\end{corollary}

By Borcherds' theory \cite{Bo2} there is a lift from weakly
holomorphic modular forms in $A_0^+(D,\chi_D)$ to Hilbert modular
forms for the group $\GK$, whose divisors are linear combinations of
Hirzebruch-Zagier divisors. Since this result is vital for us, we
state it in detail.

\begin{theorem}\label{borcherdsprod}
  (See \cite{Bo2} Theorem 13.3, \cite{Br1} Theorem 5, \cite{BB}
  Theorem 9.)  Let $f=\sum_{n\in \Z}c(n)q^n\in A_0^+(D,\chi_D)$ and
  assume that $\tilde c(n)\in \Z$ for all $n<0$. Then there is a
  meromorphic function $F(z_1,z_2)$ on $\H^2$ with the
  following properties:
  
  i) $F$ is a meromorphic modular form for $\Gamma_K$ with some
  multiplier system of finite order. The weight of $F$ is equal to
  the constant coefficient $c(0)$ of $f$. It can also be computed
  using Theorem \ref{crit}.
  
  ii) The divisor of $F$ is determined by the principal part of $f$.
  It equals \[\sum_{n<0} \tilde c(n) T(-n).\]
  
  iii) Let $W\subset\H^2$ be a Weyl chamber attached to $f$,
  i.e., a connected component of
\[
\H^2-\bigcup_{\substack{n<0\\ c(n)\neq 0}}S(-n);
\]
and define the ``Weyl vector'' $\rho_W\in K$ for $W$ by
\begin{equation}\label{defweyl}
\rho_W= \frac{ \eps_0}{\tr(\eps_0)}\sum_{n<0} \tilde c(n) 
\sum_{\lambda\in R(W,-n)} \lambda.
\end{equation}
The function $F$ has the Borcherds product expansion
\[
F(z_1,z_2)=\e(\rho_W z_1 + \rho'_W z_2) \prod_{\substack{\nu\in\frakd^{-1} \\ (\nu,W)>0}} \left(1-\e(\nu z_1 +\nu' z_2)\right)^{\tilde c(D\nu\nu')}.
\]
The product converges normally for all $(z_1,z_2)$ with $y_1 y_2 >
|\min \{n;\; c(n)\neq 0\}|/D$ outside the set of poles.

%4) There exists a positive integer $c$ such that $F^c$ has integral
%rational Fourier coefficients with greatest common divisor $1$ (at all cusps of $\GK$).

iv) The Petersson metric of $F$ is given by
\begin{equation}\label{petpsi}
-\log\|F\|_{\Pet}
=\sum_{n<0} \tilde c(n)G_{-n}(z_1,z_2).
%-\frac{1}{2}\left( -\Phi_{-n}(z_1,z_2)+L_{-n}\right).
\end{equation}

v) We have 
\begin{align*}
&\int\limits_{\GK\bs \H^2} \log\|F(z_1,z_2)\|_{\Pet} \, \omega^2 \\ 
&=  \sum_{n<0} \tilde c(n)
 \vol(T(-n))\left(2\frac{L'(-1,\chi_D)}{L(-1,\chi_D)}-
 2\frac{\sigma_{-n}'(-1)}{\sigma_{-n}(-1)}
     +1 +\log(D) \right).
\end{align*}
%where $L_{-n}$ is the constant of Definition \ref{def:L_m}.
\end{theorem}

\begin{proof}[Proof of Theorem \ref{borcherdsprod}]
The statements (i), (ii), and (iii) are proved in \cite{BB} using Theorem 13.3 of \cite{Bo2}.
Therefore we only have to verify (iv) and (v). 

By Theorem \ref{crit} the existence of $f\in A_0^+(D,\chi_D)$ implies
that condition (\ref{cond}) is fulfilled for all cusp forms $g\in
S_2^{+}(D,\chi_D)$.  Using Poincar\'e series it is easily checked that
(\ref{cond}) actually holds for all $g\in S_2(D,\chi_D)$.  Thus, by
Theorem 5 of \cite{Br1} the right hand side of (\ref{petpsi}) is equal
to the logarithm of the Petersson metric of a Hilbert modular form
$F'$ with the same divisor as $F$.  Hence the quotient
$F'/F$ is a Hilbert modular form without any zeros and poles on
$\H^2$ and thereby constant.  This shows that (\ref{petpsi}) holds up
to an additive constant.  By comparing the constant terms in the
Fourier expansions of both sides one finds that this constant equals
$0$.  Here the Fourier expansion of the right hand side is given by (\ref{cuspest}). This proves (iv). 
The last assertion follows from (iv) and \eqref{intphi1} in 
Corollary~\ref{cor:greenint}.
\end{proof}

Notice that in \cite{Br1} the assertions (i) and (ii) are deduced from (iv).
However, there a slightly different product expansion is obtained,
which involves Fourier coefficients of weakly holomorphic Poincar\'e
series of weight $2$. Similarly as in \cite{Br2} Chapter 1, these could
be related to the coefficients of weakly holomorphic modular forms of
weight $0$. In that way a more direct proof of the above theorem could
be given.  Another direct proof could be obtained by completely
arguing as in \cite{Br2}. There the Green functions
$\Phi_m(z_1,z_2,s)$ are constructed as regularized theta lifts of
non-holomorphic Hejhal-Poincar\'e series of weight $0$.  Here, for
brevity, we have preferred to argue as above. Observe that (v) also follows from \cite{Ku5}.

\begin{definition}\label{def:intbp}
Hilbert modular forms that arise as lifts via Theorem \ref{borcherdsprod} are called {\em Borcherds products}. A holomorphic Borcherds product is called {\em integral}, if it has trivial multiplier system and integral coprime Fourier coefficients.
A meromorphic Borcherds product is called integral, if it is the quotient of two holomorphic integral Borcherds products.
\end{definition}

\begin{proposition}\label{prop:intbp}
For any Borcherds product $F$ there exists a positive integer $N$ such that $F^N$ is integral.
\end{proposition}

\begin{proof}
Let $f=\sum_{n\in \Z}c(n)q^n\in A_0^+(D,\chi_D)$ as in Theorem  \ref{borcherdsprod} be the pre-image of $F$ under the Borcherds lift. It is explained in \cite{BB} that the condition $\tilde c(n)\in \Z$ for $n<0$ automatically implies that all coefficients $c(n)$ of $f$ are rational with bounded denominators. 

Thus, if $F$ is holomorphic, the Borcherds product expansion of $F$ implies that a suitable power of $F$ has integral coprime Fourier coefficients. Since the multiplier system of $F$ has finite order, we obtain the assertion in that case.

It remains to show that any meromorphic Borcherds product is the quotient of two holomorphic Borcherds products. In view of Theorem \ref{borcherdsprod} (ii) it suffices to show that there exist two weakly holomorphic modular forms $f_j=\sum_{n\in \Z}c_j(n)q^n\in A_0^+(D,\chi_D)$ such that $\tilde c_j(n)\in \Z_{\geq 0}$ for all $n<0$ (where $j=1,2$) and $f=f_1-f_2$. Then $F$ is the quotient of the Borcherds lifts of $f_1$ and $f_2$. We will now construct such  forms explicitly.

Let 
\[
E_{12}^+(\tau) =  1+\frac{2}{L(-11,\chi_D)} \sum_{n\geq 1} 
%\sum_{d\mid n} d \left( \chi_D(d) +  \chi_D(n/d)\right) 
n^5\sigma_n(-11)q^n
\]
be the normalized Eisenstein series in $M_{12}^+(D,\chi_D)$. It follows from the functional equation and the Euler product expansion of $L(s,\chi_D)$ that $L(-11,\chi_D)>0$. Moreover, $\sigma_n(-11)>0$ if $\chi_D(n)\neq-1$, and $\sigma_n(-11)=0$ if $\chi_D(n)= -1$.
Let $\Delta=q\prod_{n\geq 1}(1-q^n)^{24}$ be the Delta function, $E_4$ the normalized Eisenstein series of weight $4$ for $\Sl_2(\Z)$, and $j=E_4^3/\Delta$ the $j$-function.
The partition theoretic interpretation of $1/\Delta$ shows that all Fourier coefficients of $1/\Delta$ and $j$ with index $\geq -1$ are positive integers.
Consequently, if $c$ is a positve integer, then 
\[
g(\tau)=\frac{E_{12}^+(\tau)}{\Delta(D\tau)}j(D\tau)^{c-1}\in A_0^+(D,\chi_D),
\]
and the Fourier coefficients $b(n)$ of $g$ satisfy:
\begin{align*}
b(n)=0,&\qquad \text{if $n<-Dc$,}\\
b(n)=0,&\qquad \text{if $n\geq -Dc$ and $\chi_D(n)= -1$,}\\
b(n)>0,&\qquad \text{if $n\geq -Dc$ and $\chi_D(n)\neq -1$.}
\end{align*}
Thus, if we chose $c$ large enough, then there is a positive integer $c'$ such that
$f_1:=f+c'g$ and $f_2:=c'g$ are elements of $A_0^+(D,\chi_D)$ with the required properties.
\end{proof}

\begin{remark}
According to a result of Hecke the dimension of
$M_2^+(D,\chi_D)$ is given by $[\frac{D+19}{24}]$. 
% For $S_2^+(D,\chi_D)$ it is $[\frac{D-5}{24}]$.
In particular,
$S_2^+(D,\chi_D)=\{0\}$ for the primes $D=5,13,17$. In this case for any Hirzebruch-Zagier divisor $T(m)$ there exists a Borcherds product of weight $\varphi_m(1)$ with divisor $T(m)$. For explicit examples and calculations we refer to \cite{BB}.
\end{remark}

\subsection{Density of Borcherds products}
\label{sect:density}

\begin{definition}
Intersection points  $z\in \GK\bs\H^2$ of  Hirzebruch-Zagier divisors
are called {\em special points} (see \cite{Ge} Chapter V.6).  
%A point $z\in \GK\bs\H^2$ is called {\em special\/} if it is in an  
%intersection point of two Hirzebruch-Zagier divisors. 
\end{definition}

%In this section we show that there are ``many'' Borcherds products, whose divisors are disjoint to a given finite set of special points, and satisfy some additional technical conditions.

\begin{theorem}\label{density}
If $S\subset \GK\bs\H^2$ is a finite set of
special points, then there exist
infinitely many meromorphic Borcherds products of non-zero weight, whose divisors are disjoint to $S$ and given by linear combinations of isotropic Hirzebruch-Zagier divisors $T(p)$ with $p$ prime and coprime to $D$. 
\end{theorem}

Here and in the following, by ``infinitely many Borcherds products'' we
understand infinitely many Borcherds products whose divisors have
pairwise proper intersection.

To prove this theorem, it is convenient to view $\Sl_2(\OK)$ as an orthogonal group. We briefly recall some facts on the identification of $(\Sl_2(\R)\times\Sl_2(\R))/\{\pm (1,1)\} $ with the group $\operatorname{SO}^0(2,2)/\{\pm 1\}$ (for more details see \cite{Ge} Chapter V.4 and \cite{Bo3} Example 5.5).

We consider the quadratic space $V$ and the lattices $L=L(\OK)$, $L'=L'(\OK)$ defined at the beginning of Section \ref{sect:2.3}.

%even lattice $L$ of matrices $\lambda=\kzxz{a}{\nu}{\nu'}{b}$ with $a,b\in \Z$ and $\nu\in \OK$, with the quadratic form $q(\lambda)=- \det(\lambda)$. The signature of $L$ equals $(2,2)$, and the dual lattice $L^\vee$ 
%is given by matrices $\lambda$ as before, but with $\nu\in \frakd^{-1}$.
%The group $\GK$ acts on $L$ by $\lambda\mapsto M\lambda {M'}^t$ for 
%$M\in \GK$, the quadratic form being preserved. In that way one gets an injective homomorphism $\GK/\{\pm 1\}\to \Orth(L)$ into the orthogonal group of $L$. 

The upper half plane $\H^2$ can be identified with the Grassmannian 
\[
\Gr(L)=\{ v\subset L\otimes_\Z\R;\quad \dim(v)=2, \; q|_v<0\}
\]
of $2$-dimensional negative definite subspaces of $L\otimes_\Z\R$. The action of $\Sl_2(\OK)$ on $\H^2$ corresponds to the linear action of $\Orth(L)$ 
on $\Gr(L)$.
In terms of $\Gr(L)$ the Hirzebruch-Zagier divisor $T(m)$ is given by
\[
\bigcup_{\substack{\lambda\in L'\\ q(\lambda)=m/D}} \lambda^\perp,
\]
where $\lambda^\perp$ means the orthogonal complement of $\lambda$ in $\Gr(L)$.

For $v\in \Gr(L)$ we denote by $L_v$ the lattice $L'\cap v^\perp$ with the integral quadratic form $q_v=D \cdot q|_{L_v}$. If $g\in \Orth(L)$, then the quadratic modules $(L_v,q_v)$ and $(L_{gv},q_{gv})$ are equivalent. Therefore the equivalence classes of these quadratic forms can be viewed as invariants of the points of $\GK\bs\H^2$.
The following lemma is well known.

\begin{lemma}\label{speciallem}
Let $z\in \GK\bs\H^2$, and assume that $z$ corresponds to $v\in \Gr(L)$. Then $z\in T(m)$, if and only if the quadratic form $q_v$ on $L_v$ represents $m$.
\hfill$\square$
\end{lemma}

It is easily checked that $L_v$ has rank $2$, if and only if $v$ corresponds to a special point $z\in \H^2$.
In this case $(L_v,q_v)$ is a positive definite integral binary quadratic form.
If we write $\Delta_v$ for its discriminant, then $\Delta_v<0$ and $\Delta_v\equiv 0,1 \pmod{4}$.

\begin{lemma}\label{fffexist}
If $S\subset \GK\bs\H^2$ is a finite set of special points, then there
are infinitely many primes $p$ coprime to $D$ such that $T(p)$ is non-empty, 
isotropic, and 
$T(p)\cap S=\emptyset$. 
%
%Similarly, there are infinitely many non-trivial anisotropic Hirzebruch-Zagier
%divisors disjoint from $S$.
\end{lemma}

\begin{proof}
Let $Q_1,\dots,Q_r$ be the positive definite integral binary quadratic forms corresponding to the special points in $S$.
In view of the above discussion it suffices to show that there exist infinitely many primes $p$, which are not represented by $Q_1,\dots,Q_r$, and such that $\chi_D(p)=1$. 

Let $\Delta_j$ be the discriminant of $Q_j$. Then $\Delta_j<0$ and $\Delta_j\equiv 0,1 \pmod{4}$. 
Let  $n$ be a non-zero integer coprime to $\Delta_j$.
It is well known that $R(\Delta_j,n)$, the {\em total} representation number of $n$ by positive definite integral binary quadratic forms of discriminant 
$\Delta_j$,
is given by
\[
R(\Delta_j,n)=\sum_{d\mid n } \chi_{\Delta_j}(d)
\]
(see \cite{Za4} \S8).
Hence, if $R(\Delta_j,n)=0$ for $n$ coprime to $\Delta_j$, then $Q_j$ does not represent $n$. In particular, any prime $p$ with $ \chi_{\Delta_j}(p)=-1$ is not represented by $Q_j$.

Thus it suffices to show that there are infinitely many primes $p$ with
\begin{align*}
\chi_{\Delta_j}(p)&=-1\quad (j=1,\dots,r),\\
\chi_{D}(p)&=1.
\end{align*}
Since the $\Delta_j$ are negative and $D$ is positive, this is clearly true. In fact, even a positive proportion of primes has these properties.
\end{proof}

For the rest of this section we temporarily abbreviate $M:=M_2^+(D,\chi_D)$ and $S:=S_2^+(D,\chi_D)$. We denote the dual $\C$-vector spaces by $M^\vee$ and $S^\vee$, respectively. For any positive integer $r$, the functional
\[
a_r:M\to\C,\quad f=\sum_n b(n)q^n\mapsto a_r(f):=b(r)
\]
is a special element of $M^\vee$, and $M^\vee$ is generated by the family $(a_r)_{r\in \N}$ as a vector space over $\C$.
We denote by
$M_\Z^\vee$ the $\Z$-submodule of $M^\vee$ generated by the $a_r$ ($r\in \N$). The fact that $M$ has a basis of modular forms with integral coefficients implies that the rank of $M_\Z^\vee$ equals the dimension of $M$, and that $M_\Z^\vee\otimes_\Z\C=M^\vee$.
We write 
\[
\bar{}:M^\vee\to S^\vee,\quad a\mapsto \bar{a},
\]
for the natural map given by the restriction of a functional.

\begin{lemma}\label{speciallem2}
Let $I$ be an infinite set of positive integers $m$ with $\chi_D(m)\neq -1$, and let $A^\vee$ be the $\Z$-submodule of $M_\Z^\vee$ generated by the $a_m$ with $m\in I$. Then there is a non-zero $a\in A^\vee$ with the property that $\bar{a}=0$ in $S^\vee$. 
\end{lemma}

\begin{proof}
We consider the image $\bar{A}^\vee$ of $A^\vee$ in $S^\vee$. It is a free
$\Z$-module of rank $d\leq \dim S$.
There exist $n_1,\dots,n_d\in I$ such that $\bar{a}_{n_1},\dots,\bar{a}_{n_d}$ are linearly independent in $\bar{A}^\vee$. 
Then for any $m\in I$ there is a linear relation 
\begin{align}\label{linh1}
r_0(m)\bar{a}_m+r_1(m)\bar{a}_{n_1}+\dots+ r_d(m)\bar{a}_{n_d}=0
\end{align}
in $\bar{A}^\vee$ with integral coefficients $r_j(m)$ and $r_0(m)\neq 0$.  

If the corresponding linear combination 
\begin{align}\label{linh2}
r_0(m) B_D(m)+r_1(m) B_D(n_1)+\dots+ r_d(m) B_D(n_d)
\end{align}
of the coefficients of the Eisenstein series $E$ does not vanish, then 
$a=r_0(m)a_m+r_1(m)a_{n_1}+\dots+ r_d(m)a_{n_d}$ is a non-zero element of $A^\vee$ with the claimed property, and we are done.

We now assume that the linear combination \eqref{linh2} vanishes for all $m\in I$ and derive a contradiction. If the vector $r(m)=(r_1(m),\dots,r_d(m))$ is equal to $0$ for some $m$, then the vanishing of \eqref{linh2} implies that $r_0(m)=0$, contradicting our assumption on the $r_j(m)$. Therefore we may further assume that $r(m)\neq 0$ for all $m\in I$.

Equation \eqref{linh1} and the vanishing of \eqref{linh2} imply
\begin{align}\label{linh3}
B_D(m)\big( r_1(m)\bar{a}_{n_1}+\dots+ r_d(m)\bar{a}_{n_d} \big)
= \big( r_1(m)B_D(n_1)+\dots+ r_d(m)B_D(n_d) \big) \bar{a}_m
\end{align} 
for all $m\in I$. 
We write $\|r\|$ for the Euclidean norm of a vector
$r=(r_1, \dots, r_d)\in \C^d$ and also denote by $\|\cdot\|$ a norm on $S^\vee$, say the operator norm.
Since $\bar{a}_{n_1},\dots,\bar{a}_{n_d}$ are linearly independent, there exists an $\eps>0$ such that
\[
\| r_1\bar{a}_{n_1}+\dots+ r_d \bar{a}_{n_d} \| \geq \eps \|r \|
\]
for all $r\in \C^d$.
Moreover, there exists a $C>0$ such that
\[
|r_1 B_D(n_1)+\dots+ r_d B_D(n_d)| \leq  C   \|r \|
\]
for all $r\in \C^d$.
If we insert these estimates into the norm of \eqref{linh3}, we obtain
\[
\eps | B_D(m)| \cdot\| r(m)\|\leq C\| r(m)\|  \cdot\| \bar{a}_m\|.
\]
Since $\|r(m)\|\neq 0$, we find that
\[
| B_D(m)| \leq C/\eps   \cdot\| \bar{a}_m\|,
\]
for all $m\in I$.

By \eqref{eis}, for any $\delta>0$ the coefficients $ B_D(m)$ with $\chi_D(m)\neq 1$ satisfy $| B_D(m)|\geq C'm^{1-\delta}$ as $m\to \infty$ for some positive constant $C'$. But the Deligne bound for the growth of the coefficients of cusp forms in $S$ implies that $\| \bar{a}_m\|=O_\delta(m^{1/2+\delta})$, contradicting the above inequality.
\end{proof}

\begin{proof}[Proof of Theorem \ref{density}]
Let $I$ be the infinite set of primes coprime to $D$ that we get by Lemma \ref{fffexist}. According to Lemma \ref{speciallem2} there exists a non-zero integral finite linear combination
\[
a=\sum_{p\in I} \tilde c(p) a_p \in M_\Z^\vee
\]
with $\bar{a}=0$ in $S^\vee$. 
But then Theorem \ref{crit} implies that there is a weakly holomorphic modular form $f\in A_0^+(D,\chi_D)$ with principal part $\sum_{p\in I} c(p) q^{-p}$ and
{\em non-zero} constant coefficient $c(0)= -\frac{1}{2}\sum_{p\in I} \tilde c(p) B_D(p)$. 
If we apply Theorem \ref{borcherdsprod} to this $f$, we get a Borcherds product with the claimed properties.

We may remove the primes $p$ with $c(p)\neq 0$ from the set $I$ and repeat the above construction to get another Borcherds product. By induction we get infinitely many Borcherds products with the claimed properties.
\end{proof}

\begin{theorem}\label{densitytriple}
Let $C\in \Div(X(\GK))$ be a linear combination of Hirzebruch-Zagier divisors. Then there exist infinitely many meromorphic integral Borcherds products $F_2$ and $F_1$ of non-zero weights such that
\begin{enumerate}
\item[i)] all possible intersections of $\dv(F_1)$, $\dv(F_2)$, and $C$ are proper,
%$\dv(F_1) \cap \dv(F_2) \cap C = \emptyset$, 
\item[ii)] $F_2(\kappa)=1$ at all cusps $\kappa$ of $X(\GK)$,
%$\dv(F_2)$ is a linear combination of anisotropic Hirzebruch-Zagier divisors,
\item[iii)] $\dv(F_1)$ is a linear combination of isotropic
  Hirzebruch-Zagier divisors of prime discriminant $p$ coprime to $D$ (thus $\chi_D(p)=1$).
\end{enumerate}
\end{theorem}

\begin{proof}
Since there are infinitely many anisotropic Hirzebruch-Zagier divisors, 
Theorem \ref{crit} and Lemma \ref{speciallem2} imply that there are infinitely many Borcherds products $F_2$ of non-zero weight, whose divisor consists of anisotropic Hirzebruch-Zagier divisors, and such that $C$ and $\dv(F_2)$ intersect properly. Choose such an $F_2$. Since $F_2$ has anisotropic divisor, the Weyl vectors in the Borcherds product expansion of $F_2$ at the different cusps of $X(\GK)$ equal $0$. Consequently $F_2$ is holomorphic at all cusps with value $1$.

Let $S\subset \GK\bs \H^2$ be the finite set of intersection points 
$\dv(F_2)\cap C$. By Theorem \ref{density} there exist infinitely many Borcherds products $F_1$ of non-zero weight such that $\dv(F_1)$ is disjoint to $S$ and such that properties (i) and (iii) are fulfilled.

By possibly replacing $F_2$, $F_1$ by sufficiently large powers, we may assume that these are integral Borcherds products.
\end{proof}

In the rest of this section we essentially show that the subspace of $\Pic(X(\Gamma_K))\otimes_\Z\Q$ spanned by all Hirzebruch-Zagier divisors, is already generated by the $T(p)$ with prime index $p$ and $\chi_D(p)=1$. 

\begin{proposition}\label{density2}
  Let $I$ be a set of primes containing almost all primes $p$ with
  $\chi_D(p)=1$. Then the functionals $a_p$ with $p\in I$ generate a
  $\Z$-submodule of finite index in $M^\vee_\Z$.
\end{proposition}

\begin{proof}
Since $M$ has a basis of modular forms with integral coefficients it suffices to show that the $a_p$ with $p\in I$ generate $M^\vee$ as a $\C$-vector space.
In view of Lemma \ref{speciallem2} it suffices to show that the $\bar{a}_p$ generate
$S^\vee$. Therefore the assertion is a consequence of the following Lemma \ref{densitylem}.
\end{proof}

\begin{lemma}\label{densitylem}
Let $I$ be as in Proposition \ref{density2}.
If $f\in S$ is a cusp form which is annihilated by all $\bar{a}_p$ with $p\in I$, then $f=0$.
\end{lemma}

\begin{proof}
Since $S$ has a basis of modular forms with rational Fourier coefficients we may assume, without loss of generality, that the Fourier coefficients of $f$ are algebraic. 
By the hypothesis, and because $f\in S=S^+_2(D,\chi_D)$,  we have $a_p(f)=0$ for almost all primes $p$.

The assertion follows from the properties of the $\ell$-adic Galois representations associated to a basis of normalized newforms of $S$ using the main lemma of \cite{OS} (by a similar argument as on p.~461). 
Notice that $S$ does not contain any eigenforms with complex multiplication, since $D$ is a prime $\equiv 1\pmod{4}$.
\end{proof}

\begin{theorem}\label{densityshim}
Let $T(m)$ be any Hirzebruch-Zagier divisor. Then there exist infinitely many meromorphic integral Borcherds products $F$ of non-zero weight such that
\[
\dv(F)=\tilde c (m)T(m)+\sum_{\substack{\text{$p$ prime}\\ \chi_D(p)=1}} \tilde c(p) T(p)
\] 
with suitable integral coefficients $\tilde c(p)$ and $\tilde c(m)\neq 0$.
\end{theorem}

\begin{proof}
Let $I$ be the set of all primes $p$ with $\chi_D(p)=1$.
By Proposition \ref{density2} there exists a non-zero finite 
linear combination 
\[
a=\tilde c(m)a_m + \sum_{p\in I} \tilde c(p) a_p\in M^\vee_\Z
\]
with integral coefficients such that $\tilde c(m)\neq 0$ and $\bar{a}=0\in S^\vee$.
Therefore, in view of Theorem \ref{crit} there exists a weakly holomorphic modular form $f\in A_0^+(D,\chi_D)$ with principal part
\[
c(m)q^{-m} + \sum_{p\in I} c(p) q^{-p}
\]
and non-vanishing constant term $c(0)$. The Borcherds lift of $f$ is a Borcherds product with divisor of the required type. By taking a sufficiently large power, we may assume that it is integral.

We may now remove the primes $p$ occurring with $c(p)\neq 0$ in the above sum from the set $I$ and repeat the argument. Inductively we find infinitely many Borcherds products of the required type.
\end{proof}

\begin{remark}\label{rem:densityshim}
Let $T(m)$ be any Hirzebruch-Zagier divisor. Then there also exist infinitely many meromorphic integral Borcherds products $F$ of weight $0$ such that the same conclusion as in Theorem \ref{densityshim} holds.
\end{remark}

%%% Local Variables: 
%%% mode: latex
%%% TeX-master: "h-master"
%%% TeX-master: "h-master2"
%%% End: 

% h-rmavlev9.tex, \\

\section{Arithmetic theory of Hilbert modular surfaces}

Throughout this section we keep the assumptions of the previous sections. 
In particular $D$ is a prime congruent $1$ modulo $4$ and $K= \Q(\sqrt{D})$.
%Thus the fundamental unit  $\eps_0\in \OK$ has norm $-1$, 
%and wide equivalence of ideal classes coincides with narrow
%equivalence. 

\subsection{Moduli spaces of abelian schemes with real multiplication}  
Suppose that $A \to S$ is an abelian scheme.  
Then there exists a dual abelian scheme
$A^\vee \to S$ and a natural isomorphism $A \cong
(A^\vee)^\vee$. If $\phi: A \to B$ is a homomorphism of
abelian schemes, then there is a dual morphism $\phi^\vee: B^\vee
\to A^\vee$. A homomorphism $\mu \in \Hom(A,A^{\vee})$ is
called symmetric, if $\mu$ equals the composition $A \cong
(A^\vee)^\vee \stackrel{\mu^\vee}{\longrightarrow} A^\vee$. In this case we write
$\mu = \mu^\vee$. We denote by
$\Hom(A,A^{\vee})^{\sym}$ the space of symmetric homomorphisms.

%If $\lambda: A \stackrel{\sim}{\longrightarrow} A^{\vee}$ is a
%principal polarization, then $\mu \mapsto \eta=\lambda^{-1} \mu$
%induces an isomorphism $ \operatorname{Hom}(A,A^{\vee})
%\stackrel{\sim}{\longrightarrow} \operatorname{End}(A)$. 
%The Rosati
%involution (with respect to $\lambda$) is defined by $\eta^\dag =
%\lambda^{-1} \eta^\vee \lambda$, where $\eta$ is any endomorphism of
%$A$.  If $\eta= \lambda^{-1} \mu$ is as above, then $\mu$ is symmetric
%(with respect to $\dag$), if $\eta^\dag=(\lambda^{-1} \mu)^\dag
%=\lambda^{-1} \mu=\eta$.  We write $\End(A)^{\sym,\dag}$ for the space
%of symmetric endomorphisms.

An abelian scheme $A$ of 
relative dimension $2$ together with a ring homomorphism
$$
\iota: \OK \hookrightarrow \operatorname{End}(A)
$$
will be called an abelian surface with multiplication by $\OK$ and
denoted by the pair $(A,\iota)$.  Via $\alpha \mapsto
\iota(\alpha)^\vee$ we obtain real multiplication on the dual
abelian surface. If $\fraka\subset \OK$ is an ideal, we write $A[\fraka]$ for the $\fraka$-torsion on $A$.

Suppose from now on that
$(A,\iota)$ is an abelian surface with multiplication by $\OK$. 
 An element $\mu \in\Hom(A,A^\vee)$ is called $\mathcal{O}_K$-linear if
$\mu \iota(\alpha) = \iota(\alpha)^\vee \mu$ for all $\alpha \in \OK$.
%We denote by $\Hom_\OK(A,A^\vee)$ the space of $\OK$-linear morphisms.
%If $\lambda, \mu \in \Hom_\OK(A,A^\vee)$ and $\lambda$ is a
%principal polarization, then for any $\alpha \in \OK$ we have 
%$\eta=\lambda^{-1} \mu \in \End(A)$ and $\eta \iota(\alpha)=\lambda^{-1} \mu
%\iota(\alpha) = (\lambda^{-1} \iota(\alpha)^\vee \lambda) \lambda^{-1}
%\mu = \iota(\alpha)^\dag \eta$.  Given $\eta \in \End(A)$ we say that
%$\eta$ is $\OK$-linear, if $\eta \iota(\alpha)= \iota(\alpha)^\dag
%\eta$ and denote the space of $\OK$-linear endomorphisms by
%$\End_\OK(A)^\dag$.
%If in addition  $\eta$ is invertible in $\operatorname{End}(A)\otimes \Q$,
% then $\eta$ is 
%$\mathcal{O}_k$-linear iff  $\eta \alpha \eta^{-1} = \alpha^*$ iff 
%there exist $\eta_0$ such that $\eta= \eta_0 \iota(x)$ for some $x \in K$.   
%We denote by $\operatorname{Hom}_\OK(A,A^{\vee})^{\sym}$, 
%resp.~$\operatorname{End}_\OK(A)^{\sym,\dag}$, the space of 
%symmetric $\OK$-linear
%homomorphisms, resp. symmetric $\OK$-linear  endomorphisms.
We denote by 
\begin{align}\label{def:P(A)}
P(A)=\operatorname{Hom}_\OK(A,A^{\vee})^{\sym}
\end{align}
the space of 
symmetric $\OK$-linear
homomorphisms and write $P(A)^+$ for the $\OK$-linear symmetric polarizations.

Let $\frakl$ be a fractional ideal of $K$ and $\frakl^+$ be 
the subset of totally positive elements. An 
$\frakl$-polarization
on $(A,\iota)$ 
is  a homomorphism of $\OK$-modules
$$
\psi: \frakl \longrightarrow  
P(A)
$$
that takes $\frakl^+$ to $P(A)^+$.
We say that $\psi$
satisfies the Deligne-Pappas condition (DP), if
\begin{align*}\label{eq:cond-DP}
\psi: A\otimes_{\OK} \frakl \longrightarrow A^\vee \tag{DP},\qquad
x \otimes \lambda \mapsto \psi(\lambda)(x) \notag
\end{align*}
is an isomorphism.

The Deligne-Pappas condition \eqref{eq:cond-DP} is over schemes of
characteristic prime to $D$ equivalent to the Rapoport condition (R)
that $\Omega^1_{A/S}$ be a locally free $\mathcal{O}_S \otimes_\Z
\OK$-module. Moreover, in characteristic $0$ it holds automatically 
(see e.g.~\cite{Goren}, p.~99). 

Let $N$ be a positive integer. Suppose that $S$ is scheme over $\Spec \Z[1/N]$.
A full level-$N$ structure on an abelian surface $A$ over $S$ 
with real multiplication
by $\OK$ is an $\OK$-linear isomorphism
\begin{align*}
\left( \OK / N\OK\right)_S^2 \longrightarrow A[N]
\end{align*}
between the constant group scheme defined by 
$(\OK / N\OK)^2$ and the  $N$-torsion on $A$.

With the formulation of the next theorem that summarizes some properties
of the moduli spaces of abelian surfaces needed below, we follow
\cite{Pa}, p.~47 Theorem 2.1.2 and Remark 2.1.3; see also \cite{Ch},
\cite{Ra}, \cite{Goren}, in particular Theorem 2.17 on p.~57, Lemma
5.5 on p.~99. Furthermore, we fix a primitive $N$-th root of unity $\zeta_N$.

\begin{theorem} \label{thm:depa-stack}
  The moduli problem ``Abelian surfaces over $S$ with real
  multiplication by $\calO_K$ and $\frakl$-polarization with $(DP)$
  and full level-$N$ structure'' is represented by a regular
  algebraic stack $\mathcal{H}^\frakl(N) $, which is flat and of
  relative dimension two over $\Spec \ZZ[\zeta_N,1/N]$. It is smooth
  over $\Spec \Z[\zeta_N, 1/N D]$ and the fiber of
  $\mathcal{H}^\frakl(N)$ above $D$ is smooth outside a closed subset
  of codimension $2$. Moreover, if $N\ge 3$, then
$\mathcal{H}(N)$ is a scheme. 
\hfill $\square$
\end{theorem}

It is well known that $\Gamma(\fraka)\bs \H^2$
 can be identified with $\mathcal{H}^{\fraka\frakd^{-1}}(1)(\C)$. 
The isomorphism can be described as follows. 
(See e.~g.~\cite{Goren} Chapter II.2 for a detailed
discussion.) 

To $z=(z_1,z_2) \in \H^2$
we associate the lattice
$$
\Lambda_z = \OK z + \fraka^{-1} 
= \begin{pmatrix} \OK &\fraka^{-1}\end{pmatrix}
 \begin{pmatrix} z\\1\end{pmatrix}
%= \left\{ \alpha z + \beta ;\quad\alpha,\beta\in \OK\right\} 
= \left\{
\begin{pmatrix} \alpha z_1 +  \beta \\
\alpha' z_2 +  \beta'
\end{pmatrix} 
\in \CC^2;\quad \alpha\in \OK,\,\beta\in \fraka^{-1}\right\}.
$$
If $\gamma=\kabcd\in\Gl_2(K)$ with totally positive determinant, then 
\begin{align}\label{eq:laction}
\Lambda_{\gamma z}=\frac{1}{cz+d}\begin{pmatrix} \OK &\fraka^{-1}\end{pmatrix}
\abcd \begin{pmatrix} z\\1\end{pmatrix}.
\end{align}
In particular, if $\gamma\in \Gamma(\fraka)$, then  $\Lambda_{\gamma z}=\frac{1}{cz+d}\Lambda_z$.
For any $r \in \fraka\frakd^{-1}$ we define a 
hermitian form on $\C^2$ by
$$
H_{r,z}(u,v) =r \frac{ u_1 \bar{v}_1}{ y_1}+ r' \frac{ u_2 \bar{
    v}_2}{y_2}, \qquad (u,v \in \CC^2).$$
For $\alpha z+ \beta
,\gamma z+ \delta \in \Lambda_z \subset \C^2$ we have
$$
E_{r,z}(\alpha z+ \beta , \gamma z+ \delta ) = \Im H_{r,z}(\alpha +
\beta z , \gamma + \delta z)= \tr \left( r (\alpha \delta - \beta
  \gamma)\right) \in \Z. $$
The hermitian form $H_{r,z}$ is positive
definite (and therefore defines a polarization), if and only if $r \in
(\fraka\frakd^{-1})^+$.  

We see that $A_z= K \otimes_\QQ \CC /
\Lambda_z$ is an abelian surface.  On $A_z$ we have a natural
$\OK$-multiplication $\iota: \OK \hookrightarrow \End(A_z)$, where
$\nu \in \OK$ acts via $\iota(\nu)=\kzxz{\nu}{0}{0}{\nu'}$ on $\CC^2$.
We write $\lambda_r$ for the $\OK$-linear homomorphism 
in $P(A_z)$ 
given by $x \mapsto H_{r,z}(x, \cdot)$.
%One easily checks that $\lambda_r$ has degree $\norm( r \frakd)^2$,
%because the dual lattice of $\Lambda_z$ with respect to $E_{r,z}$ is
%isomorphic to $(r \frakd)^{-1} \oplus (r \frakd)^{-1}$.  
The
assignment $r\mapsto \lambda_r$ defines an $\OK$-linear isomorphism
$\psi:\fraka\frakd^{-1} \to P(A_z)$,
which maps the totally positive elements to polarizations. 
%The dual lattice of $\Lambda_z$ with respect to $E_{r,z}$ is
%isomorphic to $\fraka (r \frakd)^{-1}z+  (r \frakd)^{-1}$.  

%If $\eps\in \OK$ is a totally positive unit and $r_0=\eps /\sqrt{D}$, 
%then $\lambda_{r_0}$ is a principal polarization. 
%Therefore $\psi$ satisfies the Deligne-Pappas condition.

%% The last datum we need is a full level $N$-structure. We define it as
%% the $\OK$-linear isomorphism $\alpha:\left( \OK / N\OK\right)^2 \to
%% A_z[N]$ such that $\alpha(1,0)=1/N+\Lambda_z$ and
%% $\alpha(0,1)=z/N+\Lambda_z$.  The assignment $z\mapsto
%% (A_z,\psi,\alpha)$ induces the desired isomorphism $\GK(N) \bs \H^2
%% \to \mathcal{H}(N)(\C)$.

Observe that the conditions which involve $\frakl$ do in fact only
depend on the ideal class.
% $\cl(\frakl)^+$.  
To lighten notation, we 
will frequently omit the superscript $\frakl$ whenever $\frakl\cong \frakd^{-1}$ and e.~g.~simply write $\calH(N)$ for $\calH^{\frakl}(N)$.
Moreover, if $N=1$, we will 
frequently omit the argument $N$ and e.~g.~simply write 
$\calH$ instead of $\calH(1)$.
By abuse of notation, we denote the (coarse) moduli
scheme associated with $\mathcal{H}(N)$ by $\mathcal{H}(N)$, too.

%For each embedding $\sigma \in \Hom(\Q(\zeta_N),\C)$ the
%complex space $\mathcal{H}_\sigma(N) (\C)$ is isomorphic to 
%$\Gamma_K(N) \setminus \H^2$.

For the following results we refer to \cite{Ch} and \cite{Ra}.

\begin{theorem}  \label{torcomp}
There is a toroidal compactification
$h_N:\widetilde{\mathcal{H}}(N)\to \Spec \ZZ[\zeta_N,1/N] $ of $\mathcal{H}(N)$ that 
is smooth at infinity, and such that forgetting the level induces a morphism
$\pi_N:\widetilde{\mathcal{H}}(N)\to \widetilde{\mathcal{H}}(1)$
which is a Galois cover.
The complement $\widetilde{\mathcal{H}}(N) \setminus
\mathcal{H}(N)$
is a relative divisor with normal crossings. 
\hfill $\square$
\end{theorem}

%We also refer to the fact that $\widetilde{\mathcal{H}}^\frakl(N)$ is
%a Galois cover of $\widetilde{\mathcal{H}}^\frakl(1)$ with Galois
%group $G= \Sl_2( \OK/N\OK)$, and $\pi_N$ equals the natural projection
%$ \widetilde{\mathcal{H}}^\frakl(N) \to
%G \bs \widetilde{\mathcal{H}}^\frakl(N)
%=\widetilde{\mathcal{H}}^\frakl(1)$.

\begin{theorem}[$q$-expansion principle]
  There is a positive integer $n_0$ (depending on $K$ and $N$) such
  that in all weights $k$ divisible by $n_0$ there exists a line
  bundle $\mathcal{M}_k(\Gamma_K(N))$ on $\widetilde{\mathcal{H}}(N)$,
  whose global sections correspond to holomorphic Hilbert modular forms of weight
  $k$ for $\Gamma_K(N)$ with Fourier coefficients in $\ZZ[\zeta_N,
  1/N]$.  \hfill
  $\square$
\end{theorem}

We call $\mathcal{M}_k(\Gamma_K(N))$ the line bundle of Hilbert
modular forms. It is given by the $k$-th power of the Hodge bundle,
i.e., the pull-back along the zero section of the determinant of the
relative cotangent bundle of the universal family over
$\mathcal{H}(N)$.
%If $N=1,2$, then  $\mathcal{M}_k(\Gamma_K(N))$ can be defined as the $\Gamma_K(N)/\Gamma_K(N')$-invariants of $\mathcal{M}_k(\Gamma_K(N'))$ on some cover $\widetilde{\mathcal{H}}^\frakl(N')$ of sufficiently large level $N'$.

According to Proposition \ref{prop:intbp} any integral Borcherds product of weight $k$ (divisible by $n_0$) defines a rational section of $\mathcal{M}_k(\Gamma_K(N))$.

%By the Koecher principle, any section of $\mathcal{M}_k(\Gamma_K(N))$ over $\calH(N)$ extends to a section over $\widetilde\calH(N)$, that is
%$H^0(\calH(N),  \,\mathcal{M}_k(\Gamma_K(N)))
%=H^0(\widetilde\calH(N),  \,\mathcal{M}_k(\Gamma_K(N)))$ (see \cite{Ch}, Main Theorem (i)). 

\begin{theorem} The minimal
compactification of the (coarse) moduli scheme $\mathcal{H}(N)$ is given by
\begin{align}\label{eq:minhilbert}
\overline{\calH}(N) = 
\Proj  \Bigg( \bigoplus_{k\geq 0,\; n_0|k} 
H^0\left(\widetilde\calH(N),  \,\mathcal{M}_k(\Gamma_K(N))\right)\Bigg).
\end{align}
The scheme $\overline{\calH}(N)$ is  normal, projective,
and flat over $\Spec\ZZ[\zeta_N,1/N]$ (see \cite{Ch}, p.~549). 
Furthermore, for $N=1$ its fibers over $\Spec \Z$ are 
irreducible (see \cite{DePa}, p.~65). \hfill
  $\square$
\end{theorem}

For any embedding $\sigma \in
\Hom(\Q(\zeta_N),\C)$ the complex variety $\overline{\calH}(N)_\sigma(\CC)$ is isomorphic to  $X(\GK(N))$.

\begin{proposition}\label{prop:t-rational}
  For any $m$ the divisor $T(m) \subset X(\GK)=\overline{\calH}(\C)$ is defined over
  $\QQ$.
\end{proposition}
\begin{proof}  
This fact is well known see e.g.~\cite{HLR, Ge} 
for related formulations.
We sketch a different proof 
using the theory of Borcherds products. 

Let $p$  be a prime which is split in $\OK$. 
According tho Proposition \ref{prepullerback}, $T(p)$ is the image of the morphism $\varphi$ defined in \eqref{eq:defvarphi}. By Proposition \ref{prop:pullerback} the  pullback via $\varphi$ of a  Hilbert modular with rational Fourier coefficients is a modular form for $\Gamma_0(p)$ with rational coefficients.
Thus the ideal of holomorphic Hilbert modular forms vanishing along $T(p)$ in the graded ring of holomorphic Hilbert modular forms is generated by Hilbert modular forms with rational coefficients. 
Since $X(\GK)$ is equal to 
$\Proj \bigoplus_k H^0(\widetilde\calH,\calM_k(\GK))(\C)$,
we obtain the claim for $T(p)$. (See also Proposition \ref{prop:tp} below.)

If $T(m)$ is any Hirzebruch-Zagier divisor, then by Theorem
\ref{densityshim} there exists an integral Borcherds product 
whose divisor is a linear combination of $T(m)$
and divisors $T(p)$ with prime index $p$ as above. Now the claim
follows by linearity.
\end{proof}

\begin{definition} 
  We define the Hirzebruch-Zagier divisor $T_N(m)$ on the generic
  fiber $\widetilde{\calH}(N)_{\Q(\zeta_N)}$ as the pullback of $T(m)$
  on $\overline{\calH}_\Q$.  Moreover, we define the
  Hirzebruch-Zagier divisor $\mathcal{T}_N(m)$ on
  $\widetilde{\calH}(N)$ as the Zariski closure of $T_N(m)$.
\end{definition}

If $N=1$, we frequently write $\mathcal{T}(m)$ for $\mathcal{T}_1(m)$. Observe that $\mathcal{T}_N(m)=\pi_N^*\mathcal{T}(m)$.

\begin{proposition}\label{prop:moduli-divisor}
  Let $F$ be an integral Borcherds product of weight $k$, and write
  $\div_N(F)$ for the divisor on $\widetilde{\calH}(N)$ of the
  corresponding global section of $\mathcal{M}_k(\GK(N))$.  Then
  $\div_N(F)$ is equal to the Zariski closure of the induced divisor on
  the generic fiber.
%, i.e.,
%\begin{align*}
%\dv(F)= \overline{ \dv(F)\times_{\Z[\zeta_N,1/N]} \QQ(\zeta_N)}^{\text{Zar}}
%\end{align*}
In particular, if $ \dv(F)(\C)= \sum_m a(m) T(m)$ on $X(\GK)$,
then we have on $\widetilde{\calH}(N)$: 
 \begin{align*}
 \dv_N(F)= \sum_m a(m) \,\mathcal{T}_N(m).
\end{align*}
\end{proposition}

\begin{proof}
Without loss of generality we may assume that $F$ is holomorphic. 
  It suffices to show that $\dv_N(F)$ is a horizontal divisor. 
  Since $F$ is a modular form for the
  full group $\GK$, the divisor $\div_N(F)$ is the pull-back of a
  divisor on $\widetilde{\calH}$.
  Because $\widetilde{\calH}$ is
  geometrically irreducible at all primes (see \cite{DePa}), the vertical part of $\dv_N(F)$ can only contain multiples of full fibers of $\widetilde{\calH}(N)$. But 
the Borcherds product expansion implies that the Fourier coefficients of $F$ 
  are coprime. %(for all cusps).  
  Therefore, by the $q$-expansion
  principle \cite{Ch}, $\dv_N(F)$ does not contain a full fiber of
  $\widetilde{\calH}(N)$ above $\Spec \ZZ[\zeta_N,1/N]$.  
This concludes the proof of the proposition.
\end{proof}

\subsection{Modular morphisms}

We extend the morphism $\varphi:Y_0(p) \to \GK\bs \H^2$ of Section \ref{sect:4.3} to integral models by giving a modular interpretation.   
We basically follow the descriptions in \cite{Goren} and \cite{Lan} (see also \cite{KR} for related results). Moreover, we extend it to compactifications by 
means of the  $q$-expansion principle.

For the rest of this section we assume that $p$ is a prime which is
split in $\OK$ or $p=1$.  Let $\frakp$ be a prime ideal of $\OK$ above
$p$. There is a fractional ideal $\frakc$ and a totally positive
$\lambda\in K$ such that $\frakp=\lambda\frakc^2$. We may assume that $\norm(\lambda)$ is a power of $p$ (for instance if we take $\frakc=\frakp^{\frac{h+1}{2}}$, where $h$ is the class number of $K$). We fix a matrix
$M\in \kzxz{ \frakc^{-1}}{ \frakc^{-1}}{\frakc}{ \frakc}\cap
\Sl_2(K)$.

The following two lemmas were communicated to us by T.~Wedhorn.

\begin{lemma}\label{thorsten0}
Let $A/S$ and $B/S$ be abelian surfaces with
$\OK$-multiplication, such that $P(A)$ and $P(B)$ are locally constant sheaves
of projective $\OK$-modules of rank 1. 

i) Let $\pi: A \to B$ be an $\OK$-linear isogeny
whose degree is invertible on $S$ and denote by $\pi^*:P(B)\hookrightarrow P(A)$, $g\mapsto \pi^\vee g\pi$ the canonical map. 
Then the index of $\pi^*P(B)$ in $P(A)$ is
equal to the degree of $\pi$.

ii) Let $\fraka$ be an ideal of $\OK$, and let $\pi: A \to B := A/A[\fraka]$
be the canonical projection. Then
$
\fraka^2 P(A) \subset \pi^*P(B).
$
If the norm of $\fraka$ is invertible on $S$, we have equality.
\end{lemma}

\begin{proof}
i)  We write $(\ )_p$ for $(\ ) \otimes_{\Z} \Z_p$. It suffices to show
that the length of the cokernel of $\pi^*_p: P(B)_p \hookrightarrow P(A)_p$ is equal to the
$p$-adic valuation of $\deg(\pi)$ for any prime number $p$ which divides the
degree of $\pi$.
We have a commutative diagram
\begin{align*}
\xymatrix{
P(A)_p \ar[r]^-{\sim} \ar[d]_{\pi^*_p} & 
\left(\bigwedge^2_{\OK\otimes\Z_p} T_p(A)\right)^*\ar[d]^{T_p(\pi)^*}\\
P(B)_p \ar[r]^-{\sim}  & 
\left(\bigwedge^2_{\OK\otimes\Z_p} T_p(B)\right)^*
},
\end{align*}
where $T_p(A)$ and $T_p(B)$ denote the Tate modules corresponding to $A$ and $B$, and for an alternating form of $\OK\otimes\Z_p$-modules $\beta: \wedge^2_{\OK\otimes\Z_p} T_p(A) \to \OK\otimes\Z_p$ we put
$T_p(\xi)^*(\beta)=\beta\circ \wedge^2  T_p(\xi)$ (see e.g.~\cite{Vo} (1.3.4)).
Since
$$
v_p(\deg(\pi)) = \operatorname{length}(\Coker(T_p(\pi))) = \operatorname{length}(\Coker(T_p(\pi)^*))=\operatorname{length}(\Coker(\pi_p^*))
$$
% For the first equality: Mumford, Abelian varieties, \S19, Thm. 4.
we obtain the first claim.

ii)
Let $C$ be any abelian scheme over $S$ with $\OK$-multiplication.
We first claim that $\fraka\Hom_{\OK}(A,C) \subset \Hom_{\OK}(B,C)\pi$. Let
$x \in \fraka$ and $\alpha \in \Hom_{\OK}(A,C)$. We have to show that $x\alpha$
annihilates $A[\fraka]$, but this is obvious as $\alpha$ is $\OK$-linear.

Next we claim that $\fraka\Hom_{\OK}(C,A^\vee) \subset
\pi^\vee\Hom_{\OK}(C,B^\vee)$. Note that if we dualize the canonical
projection $\pi:A \to B$, $\pi^\vee$ is the canonical projection
$B^\vee \to B^\vee/B^\vee[\fraka]$. This follows from the fact that the dual
of the multiplication with an element $x \in \OK$ on any abelian scheme $A$
with $\OK$-multiplication is just the multiplication with $x$ on $A^\vee$.
Hence, the second claim follows from the first one by dualizing.

Altogether, these two claims imply that $\fraka^2\Hom_{\OK}(A,A^\vee) \subset
\pi^*\Hom_{\OK}(B,B^\vee)$ and hence $\fraka^2P(A) \subset \pi^*P(B)$.

It follows that $\fraka^2P(A) \subset \pi^*P(B) \subset P(A)$ and as $P(A)$ is a
locally constant sheaf of projective $\OK$-modules of rank 1, this is also
true for $P(B)$. Moreover, the degree of $\pi$ is $\norm(\fraka)^2$. Hence the assertion follows from (i).
\end{proof}

\begin{remark}
The hypothesis in (i) that $\deg(\pi)$ is invertible on $S$ is in
fact superfluous: As $P(A)$ and $P(B)$ are locally constant, we can assume
that $S = \Spec(k)$ where $k$ is an algebraically closed field. If $l =
\operatorname{char}(k)$ divides the degree of $\pi$ one uses the Dieudonn\'e module instead
of the Tate module and works with $P(A) \otimes_{\Z} W(k)$ where $W(k)$ is the
Witt ring of $k$. Otherwise the reasoning is the same.
In the same way we actually always have equality in (ii).
\end{remark}

\begin{lemma}\label{thorsten}
Let $S$ be a scheme over $\Z[1/p]$ and $A/S$ an abelian surface with $\OK$-multiplication, such that $P(A)$ is a locally constant sheaf of projective $\OK$-modules of rank 1.
Let $0\subsetneq H\subsetneq A[\frakp]$ be an $\OK$-invariant locally free subgroup scheme and write $\pi: A\to A/H$ for the canonical projection.
Then the image of the canonical map
$\pi^*:P(A/H)\to
%\stackrel{\pi^*}{\longrightarrow}
P(A)
%,\quad g\mapsto \pi^\vee g \pi
$
is equal to $\frakp P(A)$.
\end{lemma}

\begin{proof}
%The injectivity is clear. 
If we briefly write $B=A/H$, 
we have a commutative diagram
\begin{align*}
\xymatrix{
A \ar[r]^{\pi} \ar[dr]_{\pr} & B \ar[d]^{\xi}\\
& A/A[\frakp]
},
\end{align*}
where $\pr$ denotes the natural projection. The morphism $\xi$ exists, because $H\subset A[\frakp]$. The condition $0\subsetneq H\subsetneq A[\frakp]$ implies that both $\pi$ and $\xi$ have degree $>1$. Since the degree of $\pr$ is equal to $p^2$, we have $\deg(\pi)=\deg(\xi)=p$.
In view of Lemma \ref{thorsten0} the above diagram induces the following inclusions 
for the polarization modules:
\begin{align*}
\frakp^2 P(A)=\pr^* P(A/A[\frakp])\subsetneq \pi^*P(B)\subsetneq P(A).
\end{align*}
%\begin{align*}
%\xymatrix{
%P(A) \ar[r]^{\xi^*} \ar[dr]_{\cdot \theta^2} & P(B) \ar[d]^{\pi^*}\\
%& P(A)
%},
%\xymatrix{
%P(A)  & P(B) \ar[l]_{\pi^*}\\
%& P(A/A[\frakp])\ar[u]_{\xi^*}\ar[ul]^{\pr^*}
%},
%\end{align*}
%where all arrows are injective.
%In particular $P(B)$ is again a sheaf of projective $\OK$-modules of rank 1.
%Consequently, $\pr^* P(A)\subset \pi^*P(B)\subset P(A)$. 
%Since $\norm(\frakp)=p$, Lemma \ref{thorsten0} (i) implies that the indices for both inclusions are $p$. Moreover, by Lemma \ref{thorsten0} (ii) we know that 
%$\pr^* P(A/A[\frakp])=\frakp^2 P(A)$.
Therefore $\pi^*P(B)=\frakp P(A)$.
\end{proof}

Let $\frakl$ be a fractional ideal of $K$ and $E/S$ an elliptic curve over a scheme $S$. We consider the abelian surface 
\begin{align}\label{eq:efrakl}
E\otimes_\Z \frakl \cong E\times_S E
\end{align}
with the canonical $\OK$-action, the isomorphism being obtained by chosing a $\Z$-basis to $\frakl$.
Then the natural $\OK$-action on $\frakl$ is given by a ring homomorphism $\OK\hookrightarrow \operatorname{M}_2(\Z)$, $u\mapsto R_u$. The corresponding $\OK$-action on $E\times_S E$ is 
given by the inclusion $\operatorname{M}_2(\Z) \subset \End(E\times_S E)$.
One easily checks that for any $u\in \OK$ the conjugate acts by $R_{u'}=R_u^*{}^t= \kzxz{0}{-1}{1}{0}R_u\kzxz{0}{1}{-1}{0}$.

%the $\OK$-action is given by a ring homomorphism $\OK\hookrightarrow \operatorname{M}_2(\Z)\subset\End( E\times_S E)$ with $u\mapsto R(u)\in \operatorname{M}_2(\Z)$.

\begin{lemma}\label{inken}
The dual of the abelian surface with $\OK$-multiplication $E\otimes_\Z \frakl$ is given by
$E\otimes_\Z (\frakl\frakd)^{-1}$. Moreover, $P(E\otimes_\Z \frakl)\cong \frakl^{-2}\frakd^{-1}$, the isomorphism preserves the positivity, and the Deligne-Pappas condition (DP) holds.
\end{lemma}

\begin{proof}
If we choose a $\Z$-basis to $\frakl$, the natural $\OK$-action on $\frakl$ is given by a ring homomorphism $\OK\hookrightarrow \operatorname{M}_2(\Z)$, $x\mapsto R_x$. 
%The corresponding $\OK$-action on $E\otimes_\Z \frakl$ is 
%given by the inclusion $\operatorname{M}_2(\Z) \subset \End(E\times_S E)$.
The dual of  $\frakl$ with respect to the trace form on $K$ is equal to
$(\frakl\frakd)^{-1}$, and the dual $\OK$-action with respect to the dual basis is given by $x\mapsto R_x^t$.

On the other hand we may identify $(E\times_S E)^\vee$ canonically with 
$E\times_S E$. This identifies the dual $R^\vee$ of a morphism
$R\in \operatorname{M}_2(\Z) \subset \End(E\times_S E) $ with $R^t$.
This yields the first assertion.

We obtain an $\OK$-linear monomorphism 
$
\frakl^{-2}\frakd^{-1}\to P(E\otimes_\Z \frakl)$  by the assignment 
$x\mapsto (e\otimes l \mapsto e\otimes xl)$. One can check that totally positive elements are mapped to polarizations.
The natural composite morphism
\[
(E\otimes_\Z \frakl)\otimes_\OK (\frakl^{-2}\frakd^{-1})\longrightarrow 
(E\otimes_\Z \frakl)\otimes_\OK  P(E\otimes_\Z \frakl) \longrightarrow (E\otimes_\Z \frakl)^\vee
\]
is an isomorphism. Now the assertion follows from \cite{Vo} Proposition~3.3.
\end{proof}

\begin{proposition} \label{prop:tp} 
  There
  exists a morphism from the moduli stack $\calY_0(p)$ over
  $\Sch/\Z[1/p]$ of elliptic curves together with a cyclic subgroup of
  order $p$ to the moduli stack $\mathcal{H}$ over
  $\Sch/\Z[1/p]$ of abelian schemes with $\OK$-multiplication and
  $\OK$-polarization. Moreover, the associated morphism
  $\calY_0(p)(\C)\to\mathcal{H}(\C) $ is induced by the morphism of 
Section \ref{sect:4.3},
  $\H\to\H^2$, $\tau \mapsto M^{-1}\kzxz{\lambda}{0}{0}{1} (\tau,\tau)$.
\end{proposition}

\begin{proof} 
Let $S$ be a scheme over $\Z[1/p]$, $E$ an elliptic curve over $S$,
and $C$ be a finite locally free subgroup scheme of $E$, whose
geometric fibers are cyclic groups of order $p$.  We consider
$A=E\otimes_\Z\frakc$ as in \eqref{eq:efrakl}.  
The subgroup $C
\otimes_\Z\frakc\subset A$ of order $p^2$ is invariant under the
action of $\OK$. 
It is also invariant under the 
isomorphism $\sigma=\kzxz{0}{-1}{1}{0} \in\Sl_2(\Z)\subset 
\operatorname{End}(A)$.
%involution $\sigma$

The $\frakp$-torsion $A[\frakp]\subset A$, 
is a subgroup of order $p^2$, which is also invariant
under the action of $\OK$, but $A[\frakp]^\sigma = A[\frakp']$.  
Consequently $H = (C
\otimes_\Z\frakc) \cap
  A[\frakp]$ is a subgroup of $A$ of order $p$, which is
  $\OK$-invariant. The fact that the orders of $C
\otimes_\Z\frakc$ and $
  A[\frakp]$ are invertible on the base implies that they are both
  \'etale over $S$. Hence $H$ is also \'etale over $S$.  We have an
  exact sequence of abelian schemes with $\OK$-multiplication
  $$
  0 \longrightarrow H \longrightarrow A \longrightarrow B
  \longrightarrow 0.
  $$
  Since $H$ satisfies the conditions of
  \cite{AnGo}, 
%  [Andreatta-Goren, Geometry of Hilbert modular varieties over...,
%  IMRN, to appear], 
Corollary 3.2, the abelian scheme $B$ satisfies
  (DP) for $P(B)$, 
which is equal to 
$\frakp\frakc^{-2}\frakd^{-1}\cong\frakd^{-1}$ by Lemma \ref{thorsten} and Lemma \ref{inken}. The assignment $(E,C)\mapsto B$ is functorial and defines a morphism $\calY_0(p)\to \calH$.

  We now trace this morphism on the complex points.  Recall that
  points on $\Gamma_0(p) \bs \H$ correspond to isomorphism classes of
  elliptic curves over $\C$ together with a subgroup of order $p$ via
  the assignment $\tau\mapsto =(E_\tau,C_\tau):=
  (\C/\Lambda_\tau,\langle 1/p\rangle)$, where $\Lambda_\tau = \Z\tau+\Z$, and $\langle 1/p\rangle$ denotes the subgroup of $E_\tau$ generated by the point $1/p+\Lambda_\tau$.  The abelian surface $A$ is
  defined by the lattice $\Lambda_{\tau,\tau}=\frakc(\tau,\tau)+\frakc$.  
The lattice for the quotient
  $A/(C_{\tau}\otimes \frakc)$ is given by
  $\frakc (\tau,\tau)+\frac{1}{p}\frakc\subset\C^2$ and that for
  $A/A[\frakp]$ by
  $\frakc\frakp^{-1} (\tau,\tau)+\frakc\frakp^{-1}$.  Therefore
  the lattice for $B$ is equal to
  $\Lambda_B=\frakc (\tau,\tau)+\frakc\frakp^{-1}$.  
On the other hand, the abelian surface corresponding to the point
$M^{-1}\kzxz{\lambda}{0}{0}{1} z\in \GK\bs\H^2$ is given by the lattice
$
\tilde{\Lambda}_z=\OK M^{-1}\kzxz{\lambda}{0}{0}{1} z+\OK
$.
If we write $M=\kabcd$, and use \eqref{eq:laction}, we see that
\begin{align*}
\tilde{\Lambda}_z&=\frac{1}{-cz+a}\begin{pmatrix} \OK &\OK\end{pmatrix}
M^{-1}\zxz{\lambda}{0}{0}{1} \begin{pmatrix} z\\1\end{pmatrix}\\
&=\frac{ \lambda}{-cz+a}\begin{pmatrix}\frakc &(\lambda
\frakc)^{-1}\end{pmatrix}
\begin{pmatrix} z\\1\end{pmatrix}\\
&=\frac{ \lambda}{-cz+a}\left(\frakc z +\frakc \frakp^{-1}\right).
\end{align*}
Therefore $\Lambda_B\cong\tilde\Lambda_{(\tau,\tau)}$ and 
the associated morphism
  $\calY_0(p)(\C)\to\mathcal{H}(\C) $ is induced by 
$\tau \mapsto M^{-1}\kzxz{\lambda}{0}{0}{1} (\tau,\tau)$.
\end{proof}

\begin{remark}
i) If we consider
$A=E\otimes_\Z\OK$ in the above construction instead of 
$A=E\otimes_\Z\frakc$, we obtain a moduli description of the morphism of 
\eqref{eq:phidiag}.

ii) If we simply consider 
$E\mapsto A=E\otimes_\Z\OK$ (and do not take a quotient), we obtain a moduli description of the diagonal embedding $\Sl_2(\Z)\bs\H\to\GK$ (see \cite{Goren} Chapter~2.5).
\end{remark}

%\begin{proposition} \label{prop:t1-level}
%  There exists a morphism from the moduli stack $\mathfrak{M}(N)$ over
%  $\Sch/\Z[\zeta_N][1/N]$ of elliptic curves with full level-$N$
%  strucuture to the moduli stack $\mathcal{H}(N)$ over $\Sch/\Z$ of
%  abelian surfaces with $\OK$-multiplication and $\OK$-polarization and
%  full level-$N$ structure such that the associated morphism
%  $\mathfrak{M}_1(\C)\to\mathcal{H}(N)(\C) $ is induced by
%  $\H\to\H^2$, $\tau\mapsto (\tau,\tau)$.
%\end{proposition}  
%\begin{proof} 
%\texttt{extend proof of \ref{prop:t1}}
%\end{proof}

%\begin{proposition} \label{prop:tp-level} 
%  Let $p$ be a prime which is split in $\OK$ or $p=1$ and write
%  $p=\theta\theta'$ with $\theta\in \OK$ totally positive.  Then there
%  exist a morphism from the moduli stack $\mathfrak{M}_0(p,N)$ over
%  $\Sch/\Z[\zeta_N,1/Np]$ of elliptic curves together with a cyclic
%  subgroup scheme of order $p$ and full level-$N$ strucuture to the moduli
%  stack $\mathcal{H}(N)$ over $\Sch/\Z[\zeta_N,1/Np]$ of abelian surfaces
%  with $\OK$-multiplication and $\OK$-polarization and full level-$N$
%  structure. Moreover, the associated morphism
%  $\mathfrak{M}_0(p,N)(\C)\to\mathcal{H}(N)(\C) $ is induced by
%  $\H\to\H^2$, $\tau\mapsto (\theta\tau,\theta'\tau)$.
%\end{proposition} 
%\begin{proof} 
%\texttt{extend proof of \ref{prop:tp}}
%\end{proof}

By abuse of notation, we denote the coarse moduli scheme associated with 
$\calY_0(p)$ by $\calY_0(p)$, too. 
In the following proposition we extend the morphism of moduli schemes 
$\calY_0(p)\to \mathcal{H}$ given by Proposition 
\ref{prop:tp} 
to the minimal compactifications.
Recall that the minimal compactification $\calX_0(p)$ of $\calY_0(p)$
can be described analogously to \eqref{eq:minhilbert}
as
\begin{align}\label{eq:minmodular}
\calX_0(p) = 
\Proj \Bigg( \bigoplus_{k\geq 0,\; n_1|k} 
H^0\left(\calX_0(p),  \,\mathcal{M}_k(\Gamma_0(p))\right)\Bigg).
\end{align}
Here $n_1$ is a suitable positive integer (depending on $p$), and (for $k$ divisible by $n_1$) $\mathcal{M}_k(\Gamma_0(p))$ denotes the line bundle of modular forms on $\calX_0(p)$. 

%If $N\geq 3$, then $\mathcal{M}_k(\Gamma_0(p,N))$ is the $k$-th power of the pull-back of the relative cotangent bundle of the universal family over  
%$\calX_0(p,N)$.

%\texttt{Properties of $\calX_0(p)$, geometric description of the bundle of modular forms.}

Moreover, by the $q$-expansion principle on modular curves, the global
  sections of the bundle $\mathcal{M}_k(\Gamma_0(p))$ correspond to modular forms of weight $k$ for $\Gamma_0(p)$ with
 Fourier coefficients in $\ZZ$. 

\begin{proposition} \label{prop:tpmin} 
%Let $p$ be a prime which is split in $\OK$ or $p=1$.
There
exists a unique proper morphism  
$\bar\varphi:\calX_0(p)\to \overline{\calH}$
of schemes over $\Z[1/p]$ such that the diagram
\begin{align*}
\xymatrix{
\calX_0(p) \ar[r]^{\bar\varphi} & \overline{\calH}\\
\calY_0(p) \ar[r]^{\varphi} \ar@{^(->}[u] & \calH \ar@{^(->}[u] }
\end{align*}
commutes. Here $\varphi$ denotes the morphism given by Proposition 
\ref{prop:tp}. Moreover, $\bar\varphi$ maps the cusps of $\calX_0(p)$ to cusps of $\overline{\calH}$. The image of $\bar\varphi$ is the Zariski closure of the Hirzebruch-Zagier divisor $T(p)$ on the generic fiber $\overline{\calH}_\Q$. 
\end{proposition}

\begin{proof}
We have to prove the existence of $\bar\varphi$. The uniqueness follows from the separatedness of $\overline{\calH}$.  
%For simplicity we carry out the proof only for $N=1$. The straightforward generalization to the general case is left to the reader.

We may assume that the integer $n_1$ of \eqref{eq:minmodular} is equal to $2 n_0$.
The morphism $\varphi$ induces an isomorphism of line bundles on $\calY_0(p)$:
\begin{align}\label{bundelisom}
\varphi^*(\mathcal{M}_k(\Gamma_K))\stackrel{\sim}{\longrightarrow} \calM_{2k}(\Gamma_0(p))|_{\calY_0(p)}.
\end{align}
This can be seen 
by considering $\varphi$ on the universal objects over the moduli stacks $\calY_0(p)$ and $\calH$ and using the definition of the line bundles of modular forms in terms of relative cotangent bundles. The isomorphism \eqref{bundelisom} induces a morphism on the sections over any open subset $U\subset\calH$:
\[
\varphi^*:\mathcal{M}_k(\Gamma_K))(U)\longrightarrow \calM_{2k}(\Gamma_0(p))(\varphi^{-1} U),\quad s\mapsto s\circ\varphi.
\]
Consequently, $\varphi$ gives rise to a graded morphism of graded rings
\[
S:=\bigoplus_{k\geq 0,\; n_0 | k} 
H^0 \big( \widetilde\calH, \,\mathcal{M}_k(\Gamma_K)\big)
\stackrel{\varphi^*}{\longrightarrow} 
\bigoplus_{k\geq 0,\; n_0 | k} 
H^0\big( \calY_0(p),  \,\mathcal{M}_{2k}(\Gamma_0(p))\big).
\]
The induced morphism over $\C$ is 
(up to a scalar power of $p$ factor)
given by mapping a Hilbert modular form $F$ of weight $k$ to the modular form 
$\varphi^* (F)$
%(\tau)=\left(F\mid_k M^{-1}\kzxz{\lambda}{0}{0}{1}\right)(\tau,\tau)
%$F(\theta \tau, \theta'\tau)$ 
of weight $2k$ for the group $\Gamma_0(p)$ as in Proposition \ref{prop:pullerback} (which a priori might have singularities at the cusps of $\Gamma_0(p)$). 
%Here $\theta\in \OK$ is a totally positive integer with
%$p=\theta\theta'$. 

We first show that
$\varphi^*(F)\in  H^0\left( \calX_0(p),  \mathcal{M}_{2k}(\Gamma_0(p))\right)$
for any $F\in 
H^0 ( \widetilde\calH,  \,\mathcal{M}_k(\Gamma_K))$.
This implies that $\varphi^*$ factors into a graded homomorphism of graded rings
\[
\bigoplus_{n_0 | k} H^0\big( \widetilde\calH,  \mathcal{M}_k(\Gamma_K)\big)\stackrel{\bar\varphi^*}{\longrightarrow}
\bigoplus_{n_0 | k} H^0\big( \calX_0(p),  \mathcal{M}_{2k}(\Gamma_0(p))\big)
\]
and the homomorphism induced by the inclusion $\calY_0(p)\to \calX_0(p)$.

%It suffices to show that we have 
%\begin{align*}
%\xymatrix{
%\bigoplus_{n_0 | k} 
%H^0\left( \overline{\calH},  \mathcal{M}_k(\Gamma_K)\right)
%\ar[r]^{\varphi^*}\ar[rd]_{\bar\varphi^*} & 
%\bigoplus_{n_0 | k} 
%H^0\left( \calY_0(p),  \mathcal{M}_{2k}(\Gamma_0(p))\right)\\
% & \bigoplus_{n_0 | k} H^0\left( \calX_0(p),  \mathcal{M}_{2k}(\Gamma_0(p))\right)
%\ar[u]
%}.
%\end{align*}
%Here the right arrow is induced by the inclusion $\calY_0(p)\to \calX_0(p)$.
%It suffices to show that we have 
%$\varphi^*(F)\in  H^0\left( \calX_0(p),  \mathcal{M}_{2k}(\Gamma_0(p))\right)$
%for any $F\in 
%H^0\left( \overline{\calH},  \,\mathcal{M}_k(\Gamma_K)\right)$.

Over the complex points this follows
by looking at the Fourier
expansions.   
%If $F\in H^0( \widetilde\calH(\C),
%  \,\mathcal{M}_k(\Gamma_K))$, then $F$ can be viewed as a
%holomorphic Hilbert modular form with Fourier expansion
%$a_0+\sum_{\nu\gg 0} a_\nu \e(\nu z_1 + \nu'z_2)$ at the cusp
%$\infty$.  The Fourier expansion of $\varphi^*(F)$ at the cusp
%$\infty$ of $\Gamma_0(p)$ is given by
%\[
%\varphi^*(F)=F(\theta z_1, \theta' z_2)=a_0+\sum_{n\geq 1}
%\sum_{\substack{\nu\gg 0\\ \tr(\theta\nu)=n}} a_\nu \, q^n.
%\]
%Up to a power of $p$ the Fourier expansion of $\varphi^*(F)$ at the
%cusp $0$ is equal to the Fourier expansion of $\varphi^*(F)\mid W_p$
%at the cusp $\infty$. It equals
%\[
%\varphi^*(F)\mid W_p=a_0+ \sum_{n\geq 1} \sum_{\substack{\nu\gg 0\\ 
%    \tr(\theta\nu)=n}} a_{\nu'}\, q^n .
%\]
By Proposition \ref{prop:pullerback} (ii)
we see that $\varphi^*(F)\in H^0\left( \calX_0(p)(\C),
  \mathcal{M}_{2k}(\Gamma_0(p))\right)$.  Moreover, Hilbert modular
forms with integral Fourier coefficients are mapped to modular forms
for $\Gamma_0(p)$ with coefficients in $\Z[1/p]$.  But, since taking
$q$-expansions commutes with base change, and since the algebraic
$q$-expansions over $\C$ equal the holomorphic $q$-expansions, the
same map on $q$-expansions defines the desired homomorphism
$\bar\varphi^*$ over $\Z[1/p]$.

By \cite{Ha}, Ex.~II.2.14, $\bar\varphi^*$ induces a morphism
\[
U\stackrel{\bar\varphi}{\longrightarrow} \Proj S=\overline{\calH},
\]
defined on some open subset $U\subset\Proj \bigoplus_{n_0 | k}
H^0\left( \calX_0(p), \mathcal{M}_{2k}(\Gamma_0(p))\right)$. By
construction $U$ contains $\calY_0(p)$.  We now prove that $U=
\calX_0(p)$.  It suffices to show that there is an $F\in S_+$ such
that $\bar\varphi^*(F)$ does not vanish on the cusps $0$ and $\infty$ as
sections over $\Z[1/p]$.  By the $q$-expansion principle it suffices
to show that there is a holomorphic Hilbert modular $F$ of positive
weight with integral Fourier coefficients such that the constant
coefficient of $\bar\varphi^*(F)$ at $\infty$ (respectively $0$) is a unit in 
$\Z[1/p]$.  Since the constant coefficient of $\bar\varphi^*(F)$ at $\infty$
is up to a power of $p$ equal to the constant coefficient of $F$ at $\frakc$
(respectively the constant coefficient of $\bar\varphi^*(F)$ at $0$ is
up to a power of $p$
equal to the constant coefficient of $F$ at $\frakc^{-1}$), 
it suffices to show that there is a
holomorphic Hilbert modular form $F$ of positive weight with integral
Fourier coefficients whose constant coefficient at $\frakc$ 
(respectively $\frakc^{-1}$)
is a unit in $\Z[1/p]$.
But the existence of such an $F$ follows from \cite{Ch}, Proposition 4.5.
(Alternatively we could use an integral Borcherds product whose divisor consists of
anisotropic Hirzebruch-Zagier divisors.)

The properness of $\bar\varphi$ follows from \cite{Ha}, Corollary
II.4.8(e).  By Proposition \ref{prop:tp} we know that the
Hirzebruch-Zagier divisor $T(p)$ on the generic fiber
$\overline{\calH}_\Q$ is equal to the image of
$\bar\varphi\otimes\Q$. Since $\bar\varphi$ is proper, it follows that
the Zariski closure of $T(p)$ is contained in the image of
$\bar\varphi$. The assertion follows from the irreducibility of
$\calX_0(p)$.
\end{proof}

\begin{remark} \label{rem:mfpb}
The above argument in particular shows that for a Hilbert modular form $F$ of weight $k$ with rational coefficients the pullback to $\calX_0(p)$ of the corresponding section of $\mathcal{M}_k(\GK)$ over $\Z[1/p]$ is equal to the section of $\mathcal{M}_{2k}(\Gamma_0(p))$ over $\Z[1/p]$ corresponding to the pullback $\bar\varphi^*(F)$
%=F(\theta z_1, \theta' z_2)$ 
over $\C$.
%The above argument in particular shows that for an integral Borcherds product $F$ of weight $k$ the divisor $\dv (\varphi^* F)$
%of the section of $\mathcal{M}_{2k}(\Gamma_0(p))$ induced by
%$\varphi^*(F)=F(\theta
%z_1, \theta' z_2)$ equals the divisor $\dv (\bar \varphi^* F)$ of
%the pull-back of the section of $\mathcal{M}_k(\GK)$ associated with $F$.
\end{remark}

\begin{proposition}\label{finiteexplicit}
%Let $p$ be a prime which is split in $\OK$ or $p=1$.
Let $F$, $G$ be Hilbert modular forms for $\GK$ with rational coefficients whose divisors intersect $T(p)$ properly on $X(\GK)$. Let $\dv_N(F)$, $\dv_N(G)$ be the divisors on $\widetilde\calH(N)$ 
of the rational sections of $\mathcal{M}_{k}(\GK(N))$ associated with $F$ and $G$.
Moreover, let $S=\Spec\Z[\zeta_N, 1/Np]$ and $g:\calX_0(p)\to S$ 
be given by the structure morphism.
Then we have in $\CH^1(S)$:
\begin{align*} 
 (h_N)_* \big(\mathcal{T}_N(p) \cdot \dv_N(F) \cdot \dv_N(G)
 \big)= \deg(\pi_N) \cdot
 g_* \big(  \dv(\bar\varphi^*F) \cdot  \dv( \bar\varphi^*
 G)\big).
\end{align*}
\end{proposition}

%Let $F, G$  be integral Borcherds products whose divisors intersect $T(p)$ 
%properly on $X(\GK)$,  
%and assume that $F(\kappa)=1$ at all cusps $\kappa$ of $\GK$. 
%such that
%$$
%T(p)(\CC) \cap \dv (F)(\CC) \cap \dv (G)(\CC) = \emptyset,
%$$ 
%Let $\dv_N(F)$, $\dv_N(G)$ be the divisors on $\widetilde\calH(N)$ 
%of the sections of $\mathcal{M}_{k}(\GK(N))$ associated with $F$ and $G$.
%and let $\dv(\varphi^*F)$, $\dv(\varphi^*G)$
%be the divisors on $\calX_0(p)$ associated with the modular forms for $\Gamma_0(p)$ given by the pull-back of $F$ and $G$.
%Moreover, let $S=\Spec\Z[\zeta_N, 1/Np]$ and $g:\calX_0(p)\to S$ 
%be given by the structure morphism.
%Then we have in $\CH^1(S)_\QQ$:
%\begin{align*} 
% (h_N)_* \big(\mathcal{T}_N(p) \cdot \dv_N(F) \cdot \dv_N(G)
% \big)= \deg(\pi_N) \cdot
% g_* \big(  \dv(\varphi^*F) \cdot  \dv( \varphi^*
% G)\big).
%\end{align*}
%\end{theorem}

\begin{proof} 
In the following we consider all schemes as schemes over $S$. 
We use the diagram
%\begin{align*}
%\xymatrix{
% & \widetilde \calH(N)\ar[d]^{\pi_N}\\
%\widetilde{\mathcal{X}}_0(p)  \ar[r]^{\tilde \varphi} \ar[rd]^{\widetilde{g}} 
%&\widetilde \calH \ar[d]^h\\
% & S.
%}
%\end{align*}
\begin{align}\label{finexpldiag}
\xymatrix{
  \widetilde \calH(N)\ar[r]^{\pi_N} & \widetilde \calH \ar[r]^u & \overline \calH \ar[r]^{\bar h}  & S  \\
& \calT(p)   \ar@{^(->}[u] \ar[ru]& \mathcal{X}_0(p) \ar[u]^{\bar\varphi}\ar[ur]^g.
} 
\end{align}
We may view $F$ and $G$ as sections of the line bundle of modular forms on $\widetilde \calH(N)$, $\widetilde \calH$, and $\overline \calH$. Throughout the proof, we temporarily denote the corresponding Cartier divisors, by $\dv_N (F)$, $\dv_1 (F)$, and  $\dv (F)$, respectively (and analogously for $G$). So $\dv_N (F)=\pi_N^*\dv_1(F)$ 
and $\dv_1 (F)=u^*\dv(F)$.
We may use intersection theory for the intersection of a Cartier divisor and a cycle as described in \cite{Fu} Chapters 1, 2 on the normal schemes $\widetilde{ \calH}$ and $\overline \calH$ as intermediate steps (see loc.~cit.~Chapter 20.1).

By means of the projection formula we obtain in
$\CH^2(\widetilde \calH)$:
\begin{align*} 
(\pi_N)_* \big(\mathcal{T}_N(p) \cdot \dv_N(F) \big)
&=(\pi_N)_* \big(\pi_N^* \mathcal{T}(p) \cdot \pi_N^*\dv_1(F) \big)  \\
&= (\pi_N)_* \pi_N^* \mathcal{T}(p) \cdot \dv_1(F)\\
&=\deg(\pi_N)\cdot \mathcal{T}(p) \cdot \dv_1(F).
\end{align*}
For the latter equality we have used that $\pi_N$ is a flat morphism.
Therefore we have in $\CH^3(\widetilde \calH)$:
\begin{align}
\nonumber
  (\pi_N)_* \big(\mathcal{T}_N(p) \cdot \dv_N(F) \cdot \dv_N(G)\big)
  &=(\pi_N)_* \big(\pi_N^* \mathcal{T}(p) \cdot \pi_N^*\dv_1(F) \cdot \pi_N^*\dv_1(G)\big)  \\
\nonumber
  &= (\pi_N)_* \big(\pi_N^* \mathcal{T}(p) \cdot \pi_N^*\dv_1(F)\big) 
\cdot \dv_1(G)\\
  &=\deg(\pi_N)\cdot \mathcal{T}(p) \cdot \dv_1(F)\cdot \dv_1(G).
\label{eq:flatpb}
\end{align}
Proposition \ref{prop:tpmin} implies that $u_*(\mathcal{T}(p))=\bar \varphi_* ( \calX_0(p))$. Thus, by the projection formula we get  in
$\CH^3 (\overline{\calH})$ the equality
\begin{align*}
u_*\big(  \mathcal{T}(p) \cdot \dv_1(F) \cdot \dv_1(G) \big)&= 
u_*(  \mathcal{T}(p)) \cdot \dv(F) \cdot \dv(G) \\
&=\bar\varphi_* ( \calX_0(p)) \cdot \dv(F) \cdot \dv(G).
\end{align*}
Since $\bar\varphi$ is proper, again using the projection 
formula, we find 
\begin{align*}
  \bar \varphi_* (\calX_0(p)) \cdot \dv(F) &=
  \bar \varphi_* \big( \calX_0(p) \cdot
  \bar \varphi^* \dv(F) \big),
\end{align*}
which implies in $\CH^3(\overline\calH)$ the equality
\begin{align*}
  \bar \varphi_* (\calX_0(p)) \cdot \dv(F) \cdot
  \dv(G)&= \varphi_* \big( \calX_0(p) \cdot
  \bar\varphi^* \dv(F)  \big) \cdot \dv(G) \\
  &=\bar\varphi_* \big( \calX_0(p) \cdot \bar \varphi^*\dv(F)
  \cdot \bar \varphi^* \dv(G) \big)\\
  &=\bar \varphi_* \big( \bar\varphi^*\dv(F) \cdot \bar\varphi^*
  \dv(G) \big).
\end{align*}
%In view of diagram \eqref{finexpldiag} we have $(h_N)_* = h_* \circ
%(\pi_N)_*$ and 
Since $ g_*= \bar h_* \circ \bar \varphi_*$, we find in $\CH^1(S)$:
\begin{align*}
  \bar h_* u_* \big(\calT(p) \cdot \dv_1(F)
  \cdot \dv_1(G) \big) &= g_*\left(\bar\varphi^* \dv(F)
    \cdot \bar\varphi^*\dv(G) \right).
\end{align*}
Combining this with \eqref{eq:flatpb} and Remark \ref{rem:mfpb}, we obtain the assertion.
%
%
%To analyze the quantity on the right hand side, we consider the
%commutative diagram
%\begin{align*}
%\xymatrix{ 
%\widetilde \calX_0(p) \ar[r]^{\tilde \varphi}\ar[d]  &
% \widetilde{\calH} \ar[d]\\ 
%\calX_0(p)\ar[r]^{\bar{\varphi}} & \overline \calH
%}.
%\end{align*}
%As $\dv(F)$ is disjoint to the cusps,  we have
%\[
%\tilde g_* \big(\tilde \varphi^* \dv(F) \cdot \tilde \varphi^*
%  \dv(G)\big)=
% g_* \big(\bar\varphi^* \dv(F) \cdot \bar\varphi^*
%  \dv(G)\big).
%\]
%Finally, since $\calX_0(p)$ is smooth over $S$ and since
%$\dv(\bar\varphi^*G)$ and $\dv(\bar\varphi^*F)$ are
%horizontal divisors, we find that
% $\bar\varphi^* \dv(F)=\dv(\varphi^*F)$ and 
%$\bar\varphi^* \dv(G)=\dv(\varphi^*G)$.
%\begin{align*}
% g_* \big(\bar\varphi^* \dv(F) \cdot \bar\varphi^*
%  \dv(G)\big)
%=  g_* \big(  \dv(\bar\varphi^*F) \cdot  \dv( \bar\varphi^*
%  G)\big).
%\end{align*}
%This concludes the proof of the proposition.
\end{proof}

\section{Arithmetic intersection theory on Hilbert modular varieties}

We keep our assumption that $K$ be a real quadratic field 
with prime discriminant $D$. Moreover, let $N$ be an integer  $\geq 3$.
Throughout this section we mainly work over the arithmetic ring $\Z[\zeta_N,1/N]$. 
We put $S_N= \Spec\Z[\zeta_N,1/N]$ and 
\begin{align*}
\R_N = \R \Big/ \Big< \sum_{p|N} \Q \cdot\log(p) \Big>.
\end{align*} 
Consequently, we have an  arithmetic degree map 
\begin{align*}
 \dega: \cha^1(S_N) \longrightarrow \R_N.
\end{align*}
We consider  a toroidal compactification
$h_N: \widetilde{\calH}(N) \to S_N$ 
of the Hilbert modular variety $\calH(N)$
for the principal congruence subgroup $\GK(N)$ as in Theorem \ref{torcomp}.
It is an arithmetic variety of
dimension $2$ over $\Z[\zeta_N,1/N]$ (see
\cite{BKK}, Definition 4.3). 
If $\Sigma= \Hom(\Q(\zeta_N),\CC)$ denotes the
set of all embeddings from $\Q(\zeta_N)$ into $\C$, then the associated
analytic space is
\begin{align*}
  \widetilde{\calH}(N)_\infty=\widetilde{\calH}(N)_\Sigma(\C)=
  \coprod_{\sigma \in \Sigma}\widetilde{\calH}(N)_\sigma(\CC).
\end{align*}
Analogously we let $D_\infty= \coprod_{\sigma \in \Sigma}
\big(\widetilde{\calH}(N) \setminus \calH(N)\big)_\sigma (\C)$.  By
construction $D_\infty$ is a normal crossing divisor on
$\widetilde{\calH}(N)_\infty$ which is stable under $F_\infty$. 
We write $\mathcal{D}_\pre$ for the Deligne algebra with 
pre-log-log forms along $D_\infty$. These data give rise to
arithmetic Chow groups with pre-log-log forms denoted by
 $\cha^*(\widetilde{\calH}(N),\mathcal{D}_\pre)$.

We write 
\begin{align}\label{def:d_N}
d_N =[\Q(\zeta_N):\Q] [\GK:\Gamma_K(N)] 
\end{align}
for the total degree
of the morphism $\widetilde{\calH}(N) \to\widetilde{\calH}$ of schemes over 
$\Z[1/N]$.

 From Proposition \ref{prop:petersson} we deduce that the
 line bundle of modular forms of weight $k$ on $\widetilde{\calH}(N)$
 equipped with the Petersson metric is a pre-log singular hermitian line
 bundle, which we denote by $\overline{\mathcal{M}}_k(\GK(N))$. Its
 first arithmetic Chern class defines a class 
\begin{align*}
\cca\left(\overline{\mathcal{M}}_k(\GK(N))\right) \in
 \cha^1(\widetilde{\calH}(N),\mathcal{D}_\pre)
\end{align*}
 that one may represent
 by $\big(\dv(F), \frakg_N(F) \big)$, where $F$
 is a Hilbert modular form of weight $k$ with Fourier
 coefficients in $\Q(\zeta_N)$ and $\frakg_N(F)=
(2 \pi i k \cdot \omega, -\log\|F^\sigma\|)_{\sigma \in \Sigma}$.
The first arithmetic Chern class 
$\cca\left(\overline{\mathcal{M}}_k(\GK(N))\right)$ 
is linear in the weight $k$.

\subsection{Arithmetic Hirzebruch-Zagier divisors}
%\subsection{Arithmetic generating series}

If we consider $G_m(z_1,z_2)$ as a Green function on $\GK(N)\bs \H^2$, 
it immediately follows from Proposition \ref{prop:Gmbasic} 
that 
%the function $G_m(z_1,z_2)$,
%considered as an function on $\GK(N)\bs \H^2$, is a basic
%Green function for the cycle $T_N(m)$.  Therefore 
\[
\frakg_N(m) =
\big(-2 \partial \bar \partial G_m(z_1,z_2), G_m(z_1,z_2)\big)_{\sigma
  \in \Sigma}
\] 
is a Green object for the divisor
$T_N(m)_\infty$ on $\widetilde\calH(N)_\infty$.  
We obtain the
following arithmetic Hirzebruch-Zagier divisors:
\begin{align}\label{def:tmhat}
  \widehat{\mathcal{T}}_N(m) = ( \mathcal{T}_N(m),\frakg_N(m) ) \in
  \cha^1(\widetilde{\calH}(N),\mathcal{D}_\pre).
\end{align}
Moreover, we define for $k$ sufficiently divisible:
 \begin{align*}
\cca(\overline{\mathcal{M}}_{1/2}^\vee)= 
-\frac{1}{2k} \cca(\overline{\mathcal{M}}_k(\GK(N))) \in 
\cha^1(\widetilde{\calH}(N),\mathcal{D}_\pre)_\QQ.
\end{align*}

\begin{theorem}\label{thm:chowfinite}
  The subspace of $\cha^1(\widetilde{\calH}(N),\mathcal{D}_\pre)_\QQ$
  spanned by the $\widehat{\mathcal{T}}_N(m)$ has dimension $\dim_\C
  M_2^+(D,\chi_D)=\left[\frac{D+19}{24}\right]$. It is already generated by 
$\{ \widehat{\mathcal{T}}_N(p);\; p\in I\}$, where $I$ is any set consisting of all but finitely many primes which split in $\OK$. 
\end{theorem}

\begin{proof}
We choose for $\widehat{\mathcal{T}}_N(-n)$ the representatives given in \eqref{def:tmhat}.
  That the subspace of $\cha^1(\widetilde{\calH}(N),\mathcal{D}_\pre)_\QQ$
  spanned by the $\widehat{\mathcal{T}}_N(m)$ has dimension $\leq\dim_\C
  M_2^+(D,\chi_D)$ follows by Theorem \ref{crit} arguing as in Lemma 4.4 of \cite{Bo3}. Here one also needs Proposition \ref{prop:moduli-divisor} and its analytical counterpart Theorem \ref{borcherdsprod} (iv).
  
  On the other hand the dimsion is $\geq\dim_\C
  M_2^+(D,\chi_D)$ because of  
  Theorem 9 of \cite{Br1}, which says that any rational function on
  $X(\GK)$ whose divisor is supported on Hirzebruch-Zagier divisors is a
  Borcherds product.

Let $I$ be a set of primes as above and $T(m)$ any Hirzebruch-Zagier divisor.
According to Remark \ref{rem:densityshim}, there exists an integral Borcherds product of weight $0$, and thereby a rational function on $\widetilde{\calH}(N)$, with divisor  $\tilde{c}(m) T(m) + \sum_{p\in I} \tilde{c}(p) T(p)$ on $X(\GK)$ and $\tilde{c}(m)\neq 0$.   
We may conclude by Proposition
 \ref{prop:moduli-divisor} and Theorem \ref{borcherdsprod} (iv) that
 \[
\widehat{\mathcal{T}}_N(m) =-\frac{1}{\tilde{c}(m)}\sum_{p\in I} \tilde{c}(p) \widehat{\mathcal{T}}_N(p)
\in \cha^1(\widetilde{\calH}(N),\mathcal{D}_\pre)_\QQ.
\]
This proves the theorem.
\end{proof}

\begin{theorem}\label{thm:genser}
The arithmetic generating series
\begin{align}\label{def:genser}
\widehat{A}_{N}(\tau)=\cca(\overline{\mathcal{M}}_{1/2}^\vee) + 
\sum_{m>0} \widehat{\mathcal{T}}_N(m) q^m
\end{align}
%where $q=\e(\tau)$ with $\tau\in \H$, 
is a modular form in
$M^+_2(D,\chi_D)$ with values in
$\cha^1(\widetilde{\calH}(N),\mathcal{D}_\pre)_\QQ$, i.e., an element
of
$M^+_2(D,\chi_D)\otimes_\Q\cha^1(\widetilde{\calH}(N),\mathcal{D}_\pre)_\QQ$.
\end{theorem}

\begin{proof} By Theorem \ref{thm:chowfinite} and in view of Corollary \ref{critcor} 
%Theorem 3.1 in \cite{Bo3} and its reformulation 
%given in \texttt{Kudla/dmv}, 
  it suffices to show that for any weakly holomorphic modular form $f=
  \sum_n c(n)q^n$ in $A_0^+(D,\chi_D)$ we have the relation
\begin{align}\label{eq:gs-condition}
\tilde c(0) \cca(\overline{\mathcal{M}}_{1/2}^\vee) + 
\sum_{n<0} \tilde c(n)\widehat{\mathcal{T}}_N(-n) =0 
\end{align}
in $\cha^1(\widetilde{\calH}(N),\mathcal{D}_\pre)_\Q$.  Since
$A_0^+(D,\chi_D)$ has a basis of modular forms with rational
coefficients, it suffices to check \eqref{eq:gs-condition} for those
$f$ with $\tilde c(n)\in \Z$ for $n<0$.  Then, by Theorem
\ref{borcherdsprod}, there exists a Borcherds product $F$ of weight
$c(0)$ and divisor $\sum_{n<0} \tilde{c}(n) T(-n)$ on $X(\GK)$.  We
may assume that $F$ is an integral Borcherds product and therefore defines a section of $\mathcal{M}_{c(0)}(\GK(N))$.  If we choose
for $\widehat{\mathcal{T}}_N(-n)$ the representatives of 
\eqref{def:tmhat}, we may conclude by Proposition
\ref{prop:moduli-divisor} and Theorem \ref{borcherdsprod} (iv) that
\begin{align*}
\sum_{n<0} \tilde c(n)\widehat{\mathcal{T}}_N(-n) &=
\sum_{n<0} \tilde c(n) \big( \mathcal{T}_N(-n) , 
\frakg_N(-n)\big) = \big( \dv_N(F) ,\, \frakg_N(F)\big).
\end{align*} 
By \eqref{eq:ca1} the right hand side 
of the latter equality
equals $\cca(\overline{ \mathcal{M}}_{c(0)}(\GK(N)))$.
Using the linearity of the arithmetic Chow groups we obtain \eqref{eq:gs-condition}, hence the 
assertion.
\end{proof}

%If $F$ is a Hilbert modular form of weight $k$ for $\Gamma_K$ with 
%Fourier coefficients in $\ZZ$, we write $F_\Z$ for the corresponding global section of $\mathcal{M}_k(\Gamma_K)$.

\subsection{Arithmetic intersection numbers and Faltings heights}

Recall our convention that for an arithmetic cycle $\alpha\in \cha^3(\widetilde{\calH}(N),\mathcal{D}_\pre)_\Q$ we frequently write $\alpha$ instead of $\dega((h_N)_*\alpha)$.

\begin{lemma}\label{lem:dega-T-raw}  If $p$
  is a prime which is split in $\OK$ or $p=1$, then there is a
  $\gamma_p \in \Q$ such that we have in $\R_N$:
\begin{align*} 
  & \widehat{\mathcal{T}}_N(p) \cdot
    \cca(\overline{\mathcal{M}}_k(\GK(N)))^2 \\& = - k^2 d_N
 \vol(T(p)) \left(2\frac{ \zeta_K'(-1)}{ \zeta_K(-1)} + 2
    \frac{\zeta'(-1)}{\zeta(-1)} +3 + \log(D) \right)
+ \gamma_p \log(p).
%+ \sum_{\substack{\text{$q$ prime}\\ q\mid pN}} \gamma_q \log(q).
\end{align*} 
\end{lemma}
In the proof of the next theorem below 
%which relies on Lemma \ref{lem:dega-T-raw}, 
we will show that in fact $\gamma_p=0$.

\begin{proof} Throughout the proof 
  all equalities which contain the image of the arithmetic degree map $\dega$
  are equalities in $\R_N$.

Without loss of generality we may assume $k$ sufficiently 
  large.  In order to prove the lemma we represent
  $\widehat{\mathcal{T}}_N(p)$ by $(\mathcal{T}_N(p),
  \mathfrak{g}_N(p))$. We also choose integral Borcherds products $F$,
  $G$ of non-zero weight $k$ such that $F(\kappa)=1$ at all cusps
  $\kappa$ and such that all possible intersections on $X(\GK)$ of
  $T(p)$,  $\dv(F)$,  $\dv(G)$ are proper. 
%$$
%T(p)(\CC) \cap \dv (F)(\CC) \cap \dv (G)(\CC) = \emptyset.
%$$
Such Borcherds products exist by Theorem \ref{densitytriple}.  We
take the pairs $(\dv_N(F), \frakg_N(F))$ and $(\dv_N(G),\frakg_N(G))$ as
representatives for $\cca(\overline{\mathcal{M}}_k(\GK(N)))$, where we have put
$\frakg_N(F)= ( 2 \pi i k \omega,
-\log\|F^\sigma(z_1,z_2)\|)_{\sigma \in\Sigma} $ and 
$\frakg_N(G)$ analogously.  
%Let $h_N: \widetilde{\calH}(N) \to S_N$ be the structure morphism. 
With our convention on arithmetic intersection
numbers, we have
\begin{align}\label{eq:star-t(p)} 
&
\widehat{\mathcal{T}}_N(p) \cdot
    \cca(\overline{\mathcal{M}}_k(\GK(N)))^2  \notag\\
%& =
%\dega (h_N)_* \left( \widehat{\mathcal{T}}_N(p) \cdot 
%\cca(\overline{\mathcal{M}}_k(\GK(N)))^2 \right) \notag\\
&= \dega \left(  (h_N)_* \big(\mathcal{T}_N(p) \cdot \dv_N(F) \cdot \dv_N(G)
  \big) \right) + \frac{1}{(2\pi i)^2}
\int\limits_{\widetilde{\calH}(N)_\infty} \mathfrak{g}_N(p) *
\frakg_N(F) * \frakg_N(G).
\end{align}
For the first summand of \eqref{eq:star-t(p)} we obtain by Proposition
\ref{finiteexplicit}:
\begin{align*} 
  \dega &\left( (h_N)_* \big(\mathcal{T}_N(p) \cdot \dv_N(F) \cdot
    \dv_N(G)
    \big) \right)\\
  &= \deg(\pi_N) \cdot \dega\left( g_* \big(\dv(\varphi^*F) \cdot
    \dv(\varphi^*G)\big)\right) + \gamma'_p \log(p)
\end{align*} 
for some $\gamma'_p \in \Q$.

Since $F$ and $G$ are integral Borcherds products, they are invariant
under $\Gal(\C/\Q)$. 
%we have $F^\sigma=F$ and $G^\sigma=G$. 
Therefore, by means of Remark \ref{rem:subgroupstar},
%using the projection formula for
%star products on subgroups of $\GK$ 
we find for the second summand of \eqref{eq:star-t(p)}
\begin{align*} 
\int\limits_{\widetilde{\calH}(N)_\infty} \mathfrak{g}_N(p) *
\frakg_N(F) * \frakg_N(G)
= [\Q(\zeta_N):\Q] [\GK:\Gamma_K(N)] 
   \int\limits_{\widetilde{X}(\Gamma_K)}&
  \mathfrak{g}(p)*\mathfrak{g}(F)*\mathfrak{g}(G).
\end{align*}
The latter integral was calculated in Theorem \ref{thm:hilbertstar}:
\begin{align*} 
  \frac{1}{(2\pi i)^2} \!\!\int\limits_{\widetilde{X}(\Gamma_K)}\!\!
  \mathfrak{g}(p)*\mathfrak{g}(F)*\mathfrak{g}(G) &=-k^2 \vol(T(p))
  \left(2\frac{ \zeta_K'(-1)}{ \zeta_K(-1)}
    + 2 \frac{\zeta'(-1)}{\zeta(-1)} +3 + \log(D) \right)\\
  &\phantom{=}{} - k^2 \vol(T(p)) \frac{p-1}{p+1}\log(p)
  -\big(\div(\varphi^*F), \dv(\varphi^*G) \big)_{
    \mathcal{X}_0(p),\fin}.
\end{align*}
Observe that in this formula the contribution of the finite primes is
calculated with respect to a regular model of $\calX_0(p)$ over $\Spec
\Z$.  Since $ \deg(\pi_N)=[\GK:\Gamma_K(N)]$, and because
\begin{align*}
\dega\left( g_* \big(\dv(\varphi^*F) \cdot 
  \dv(\varphi^*G)\big)\right) 
&=[\Q(\zeta_N):\Q] \big(\div(\varphi^*F), 
\dv(\varphi^*G)
  \big)_{ \mathcal{X}_0(p),\fin},
\end{align*}
 we obtain the assertion.
\end{proof}

\begin{theorem} \label{thm:MAIN}
  Let
  $\overline{\mathcal{M}}_k(\GK(N))$ be the line bundle of modular
  forms of weight $k$ on $\widetilde{\calH}(N)$ equipped with the
  Petersson metric as in Definition \ref{def:petersson}.  Then in
  $\R_N$ we have the following identities of arithmetic intersection
  numbers:
\begin{align*}
&
%\dega \left((h_N)_* \left( \widehat{A}(\tau) \cdot 
%\cca(\overline{\mathcal{M}}_k(\GK(N)))^2 \right)\right)
\widehat{A}_{N}(\tau) \cdot 
\cca(\overline{\mathcal{M}}_k(\GK(N)))^2
= \frac{k^2}{2} d_N  \zeta_K(-1)  \left(
    \frac{\zeta_K'(-1)}{ \zeta_K(-1)} + \frac{\zeta'(-1)}{\zeta(-1)} +
    \frac{3}{2} + \frac{1}{2} \log(D) \right)
\cdot E(\tau).
\end{align*}
Here $\widehat{A}_{N}(\tau)$ is the arithmetic generating series 
\eqref{def:genser} and $E(\tau)$ is the 
holomorphic Eisenstein series of weight $2$ \eqref{eis}.
In particular, for the arithmetic
  self intersection number of
  $\overline{\mathcal{M}}_k(\GK(N))$ we
  have in $\R_N$:
\begin{align}\label{thm:ardeg-Mk} 
  \overline{\mathcal{M}}_k(\GK(N))^3
    & = -k^3 d_N
\zeta_K(-1)  \left(
    \frac{\zeta_K'(-1)}{ \zeta_K(-1)} + \frac{\zeta'(-1)}{\zeta(-1)} +
    \frac{3}{2} + \frac{1}{2} \log(D) \right).
\end{align}
\end{theorem}

\begin{proof}  Throughout the proof 
  all equalities which contain the image of the arithmetic degree map $\dega$
  are equalities in $\R_N$.

Let us start with proving \eqref{thm:ardeg-Mk}. By
  definition we have $$
\overline{\mathcal{M}}_k(\GK(N))^3 =
  \cca(\overline{\mathcal{M}}_k(\GK(N)))^3=
\dega
  \left((h_N)_* \left( \cca(\overline{\mathcal{M}}_k(\GK(N)))^3
    \right)\right).
  $$
  By the multi-linearity of the intersection product we may assume
  $k$ to be sufficiently large. In view of Theorem \ref{densityshim},
  there exists a Borcherds product $F$ of positive weight $k$ with
  divisor $\dv(F)= \sum_p a(p) T(p)$, where $a(p)\neq 0$ only for
  primes $p$ with $\chi_D(p)=1$. 
(Apply Theorem \ref{densityshim} with $N$ equal to any given prime
  $p_0$ with $\chi_D(p_0)=1$.)
We have shown in the proof of
  Theorem \ref{thm:genser} that any such Borcherds product gives rise to a
  decomposition
  $$
  \cca(\overline{\mathcal{M}}_k(\GK(N))) = \sum_p a(p)
  \widehat{\mathcal{T}}_N(p) .$$
  Hence, in view of Lemma 
  \ref{lem:dega-T-raw}, we obtain the equality
\begin{align*} 
&   \cca(\overline{\mathcal{M}}_k(\GK(N)))^3 \\
&= - k^2 d_N 
\sum_p a(p) \vol(T(p)) \left(2\frac{ \zeta_K'(-1)}{ \zeta_K(-1)} + 2
    \frac{\zeta'(-1)}{\zeta(-1)} +3 + \log( D) \right) 
+\sum_{\substack{p \\ a(p)\neq 0}} \gamma_p \log(p)
\end{align*} 
for some $\gamma_p \in \Q$.
There also exists a Borcherds product $G$ 
with divisor  $\dv(G)= \sum_p b(p) T(p)$ where $b(p)\neq 0$ only 
for primes $p$ with $\chi_D(p)=1$ and $a(p)=0$. By means of the 
well-definedness of the arithmetic intersection numbers we find 
$$
  \sum_{\substack{p \,\textrm{prime} \\ a(p)\neq 0}} \gamma_p \log(p)
  = \sum_{\substack{p \,\textrm{prime}\\ b(p)\neq 0}} \gamma_p \log(p)
  \in \R_N ,
$$
which by the unique factorization in $\Z$ in turn implies that all
$\gamma_p$ vanish. 
This shows that the formula of Lemma \ref{lem:dega-T-raw} actually holds with $\gamma_p=0$.
Since $F$ is a Borcherds product of weight $k$ we have
the relation 
$$
 \frac{1}{2} \zeta_K(-1) k= \sum_p a(p) \vol(T(p)).
$$ 
This concludes the proof of    
\eqref{thm:ardeg-Mk}. 

Now let $T(m)$ be any Hirzebruch Zagier divisor.  In view of Theorem
\ref{densityshim} there exists a Borcherds product $F$ whose divisor
is equal to $a(m) T(m) + \sum_p a(p) T(p)$, where $a(m)\neq 0$ and the $a(p)$
are integral coefficients, and $a(p)\neq 0$ only for primes $p$ with
$\chi_D(p)=1$. Therefore
\begin{align*} 
&  \widehat{\mathcal{T}}_N(m) \cdot
    \cca(\overline{\mathcal{M}}_k(\GK(N)))^2 \\
&= \frac{1}{a(m)} \left(  
    \cca(\overline{\mathcal{M}}_k(\GK(N)))^3 
 - \sum_p a(p) 
  \widehat{\mathcal{T}}_N(p) \cdot
    \cca(\overline{\mathcal{M}}_k(\GK(N)))^2\right) \\
&= - k^2 d_N \vol(T(m)) \left(2\frac{ \zeta_K'(-1)}{ \zeta_K(-1)} + 2
    \frac{\zeta'(-1)}{\zeta(-1)} +3 + \log (D) \right).
\end{align*}
In the last equality we have used the relation $a(m)
\vol(T(m))=\frac{\zeta_K(-1)}{2}k- \sum_p a(p) \vol(T(p))$. Hence the
clain follows from \eqref{eq:frankehausmann}.
%
%Finally, since the coefficients $B_D(m)$ of the
% Eisenstein series $E(\tau)$  satisfy 
%$\vol(T(m)) = - \frac{\zeta_K(-1)}{4}  B_D(m)$, 
% we obtain the claim. 
\end{proof}

\begin{theorem}\label{thm:Tmheight} 
If $T(m)$ is an anisotropic Hirzebruch-Zagier divisor, 
then the Faltings height in $\R_N$
of its model  $\mathcal{T}_N(m) \in 
\operatorname{Z}^1_U(\widetilde{\calH}(N))$  
is given by 
\begin{align*}
%\label{eq:fh_tm}
\fh_{\overline{\mathcal{M}}_k(\Gamma_K(N))}
(\mathcal{T}_N(m))&= - (2k)^2 d_N \vol(T(m))
 \left( \frac{\zeta'(-1)}{\zeta(-1)} + \frac{1}{2}
+\frac{1}{2}\frac{\sigma_m'(-1)}{\sigma_m(-1)} 
%\frac{\sum_{d\mid m} d ( \chi_D(d) + \chi_D(m/d)) \log(d)}
%{2 \sum_{d\mid m} d ( \chi_D(d) + \chi_D(m/d))} 
%+\frac{1}{4}\log(m) 
\right).
\end{align*}
\end{theorem}

\begin{proof}  By means of formula \eqref{eq:log-sing-heightpairing} the 
Faltings height equals
\begin{align}\label{eq:fh}
&\fh_{\overline{\mathcal{M}}_k(\GK(N)) }
(\mathcal{T}_N(m)) 
\\
\nonumber
&= \widehat{\mathcal{T}}_N(m)
\cdot \cca( \overline{\mathcal{M}}_k(\GK(N)))^2
- \sum_{\sigma\in \Sigma} \frac{1}{(2 \pi i)^2} 
 \int\limits_{\widetilde{\calH}(N)_\sigma} 
G_m  \land \textrm{c}_1(\overline{\mathcal{M}}_k(\GK(N)))^2 .
\end{align}
The first term is  computed in Theorem \ref{thm:MAIN}. The integral is, by functoriality of the intersection product, equal to
$d_N$ times the 
quantity from \eqref{intphi1} of Corollary \ref{cor:greenint}.
\end{proof}

\begin{remark}
i) Observe that in view of Lemma \ref{lem:intganz} we may calculate
formula \eqref{eq:fh} for all $T(m)$.  We may take this quantity as an
ad hoc definition for the Faltings height of $\mathcal{T}_N(m)$ with
respect to $\overline{\mathcal{M}}_k(\GK(N))$.  For example, if $m$ is
square free and $\chi_D(p)=1$ for all primes $p$ dividing $m$, then
up to boundary components the normalization of $T_N(m)$ is isomorphic to
some modular curve (see Remark \ref{rem:bc}). Although $\mathcal{T}_N(m) \notin 
\operatorname{Z}^1_U(\widetilde{\calH}(N))$, we obtain in $\R_N$:
\begin{align*}
\fh_{\overline{\mathcal{M}}_k(\GK(N))}
({\mathcal{T}_N(m)})
= -(2k)^2 \vol(T_{N}(m)) \left( \frac{\zeta'(-1)}{\zeta(-1)} + \frac{1}{2}
-  \frac{1}{4} 
 \sum_{p|m}  \frac{p-1}{p+ 1} \log (p)   \right).
\end{align*}
Up to the sum over primes dividing the level $m$ this
equals  the
arithmetic self intersection number of the line bundle of modular
forms on that modular curve as computed in \cite{Kue2}.
% (see Corollary 6.2 of \cite{Kue2}).
% The last summand in above formula is easily explained by the change
% of metrics. 
It would be interesting to obtain a geometric interpretation of the
difference using the change of moduli problems over $\Z$ extending
Proposition \ref{prop:tp}.  
At the
primes dividing $m$ the morphism of line bundles \eqref{bundelisom}
is not necessarily an isomorphism anymore.

%
%We expect calculations similar to 
%section 10 of \cite{KRY2}, but in our cases a version of Raynauds theorem
%on the change of Faltings heights of higher dimensional subschemes is needed.

ii) Recall that if $T(m)\subset X(\GK)$ is anisotropic, then its
normalization is isomorphic to some Shimura curve of
discriminant $m$. For example, if $m$ is square free with an even
number of prime factors, and $\chi_D(p)=-1$ for all primes $p$ dividing
$m$, then $T(m)$ is anisotropic and our formula gives in $\R_N$:
\begin{align}\label{eq:kudlaconj}
  \fh_{\overline{\mathcal{M}}_k(\GK(N))}({\mathcal{T}_N(m)})
&= -(2k)^2    \vol(T_N(m)) \left( \frac{\zeta'(-1)}{\zeta(-1)} +
    \frac{1}{2}- \frac{1}{4} \sum_{p|m} \frac{p+1}{p-1} \log (p)
% +    \frac{1}{4}\log n
 \right).
\end{align}
Related formulas have been obtained by Maillot and Roessler
\cite{MaRo} using completely different techniques. 

In recent work, Kudla, Rapoport, and Yang stated a conjecture for the
arithmetic self intersection number of the Hodge bundle on a Shimura
curve (see \cite{KRY2}).  In forthcoming work \cite{KuKue} it will be
shown that \eqref{eq:kudlaconj}, together with versions of Propositions
\ref{prop:pullerback} and \ref{prop:tp} for Shimura curves, can be
used to obtain a proof of their conjecture.

\end{remark}

Up to now we have computed the arithmetic intersection numbers and the
Faltings height in the ring $\R_{N}$ due to the lack of an arithmetic
intersection theory over stacks that prevents us to work directly on
$\widetilde {\calH}(1)$. Nevertheless we can use the whole tower of
schemes $\{\widetilde {\calH}(N)\}_{N\ge 3}$ as a substitute for the
stack $\widetilde {\calH}(1)$. Although it is more canonical to work
with the whole tower, it is enough to pick up two
integers $N,M\ge 3$ with $(N,M)=1$. We first observe that, by the
unique factorization in $\mathbb{Z}$, the sequence
\begin{equation}\label{eq:exactRN}
  0\longrightarrow \R\longrightarrow \R_{N}\oplus
  \R_{M}\longrightarrow  \R_{NM},
\end{equation}
where the first map is the diagonal and the second is the difference,
is exact. For $n|m$ we denote by $\pi _{m,n}:\widetilde
{\calH}(m)\to \widetilde {\calH}(n)$ the natural
mophism. We then define
$\cha^{\ast}(\widetilde{\calH}_{N,M},\mathcal{D}_{\pre})$ as the
kernel of the morphism
\begin{displaymath}
\xymatrix{  
\cha^{\ast}(\widetilde{\calH}(N),\mathcal{D}_{\pre})\oplus
  \cha^{\ast}(\widetilde{\calH}(M),\mathcal{D}_{\pre})
  \ar[rrr]^-{\pi _{NM,N}^*-\pi _{NM,M}^*}&&&
  \cha^{\ast}(\widetilde{\calH}(NM),\mathcal{D}_{\pre})}.
\end{displaymath}
The group
$\cha^{\ast}(\widetilde{\calH}_{N,M},\mathcal{D}_{\pre})_{\Q}$ has a
commutative and associative product defined componentwise. 

Now we define the arithmetic Hirzebruch-Zagier divisors, 
the arithmetic first Chern class of the line bundles of modular
forms and the arithmetic generating series in this group as
\begin{align*}
    \widehat{\mathcal{T}}_{N,M}(m)&=\big( \widehat{\mathcal{T}}_N(m) ,
    \widehat{\mathcal{T}}_M(m) \big), \\
    \cca(\overline{\mathcal{M}}_{k,N,M})&=
    \big(\cca(\overline{\mathcal{M}}_k(\GK(N))),
    \cca(\overline{\mathcal{M}}_k(\GK(M)))\big),\\
    \widehat{A}_{N,M}(\tau)&=\big(\widehat{A}_{N}(\tau),\widehat{A}_{M}(\tau)\big).
\end{align*}

%% Using the
%% projection formula we have, for $x\in
%% \cha^{3}(\widetilde{\calH}(N)\mathcal{D}_{\pre})$, 
%% \begin{displaymath}
%%   \dega (\pi _{NM,N}^{\ast}x)=
%%   \deg (\pi _{NM,N}).\dega(x)=\frac{d_{NM}}{d_{N}}\dega (x). 
%% \end{displaymath}
%% Therefore if $(x,y)\in
%% \cha^{3}(\widetilde{\calH}(N,M),\mathcal{D}_{\pre})$ then 
%% $  \dega(x)/d_{N}=\dega(y)/d_{M} $ in  $\mathbb{R}_{NM}$.
If $(x,y)\in
\cha^{3}(\widetilde{\calH}_{N,M},\mathcal{D}_{\pre})$, 
then 
%in view of the projection formula we have 
by the projection formula  
$  \dega(x)/d_{N}=\dega(y)/d_{M} $ in  $\mathbb{R}_{NM}$.
Thus, by the exactness of the sequence \eqref{eq:exactRN}, the
pair $(\dega(x)/d_{N},\dega(y)/d_{M})$ induces a well defined map
\begin{displaymath}
\dega:\cha^{3}(\widetilde{\calH}_{N,M},\mathcal{D}_{\pre}) 
\longrightarrow \mathbb{R}. 
\end{displaymath}
Then Theorem \ref{thm:MAIN} implies: 

\begin{theorem} \label{thm:BR}
For all coprime integers $N,M\geq 3$ we have the following identity 
of arithmetic intersection numbers in $\R$:
\begin{align*}
&
\widehat{A}_{N,M}(\tau) \cdot 
\cca(\overline{\mathcal{M}}_{k,N,M})^2
= \frac{k^2}{2} \zeta_K(-1)  \left(
    \frac{\zeta_K'(-1)}{ \zeta_K(-1)} + \frac{\zeta'(-1)}{\zeta(-1)} +
    \frac{3}{2} + \frac{1}{2} \log(D) \right)
\cdot E(\tau).
\end{align*}
\end{theorem}

Analogously we can define the group of cycles
$\operatorname{Z}^1_{U}(\widetilde{\calH}_{N,M})$ and a
height pairing 
\begin{displaymath}
(\cdot\mid\cdot):\,\cha^{2}(\widetilde{\calH}_{N,M},\mathcal{D}_{\pre})
\otimes{\rm Z}^{1} _
{U}(\widetilde{\calH}_{N,M})\longrightarrow \mathbb{R}.
\end{displaymath}
The pair $\mathcal{T}_{N,M}(m)= (\mathcal{T}_N(m)
,\mathcal{T}_M(m) )$ belongs to
$\operatorname{Z}^1_{U}(\widetilde{\calH}_{N,M})$ and Theorem
\ref{thm:Tmheight} implies:

\begin{theorem} \label{thm:CR}
If $T(m)$ is an anisotropic Hirzebruch-Zagier divisor, 
then 
%for all coprime integers $N,M\geq 3$ 
the Faltings height
of  $\mathcal{T}_{N,M}(m) \in 
\operatorname{Z}_U^1(\widetilde{\calH}_{N,M})$  
is given by 
\begin{align*}
%\label{eq:fh_tm}
\fh_{\overline{\mathcal{M}}_{k,N,M}}
(\mathcal{T}_{N,M}(m))&= - (2k)^2 \vol(T(m))
 \left( \frac{\zeta'(-1)}{\zeta(-1)} + \frac{1}{2}
+\frac{1}{2}\frac{\sigma_m'(-1)}{\sigma_m(-1)} 
%\frac{\sum_{d\mid m} d ( \chi_D(d) + \chi_D(m/d)) \log(d)}
%{2 \sum_{d\mid m} d ( \chi_D(d) + \chi_D(m/d))} 
%+\frac{1}{4}\log(m) 
\right).
\end{align*}
\end{theorem}

\begin{remark}
Finally, we like to
remark that if there were an arithmetic Chow ring for the stack
$\widetilde\calH(1)$ satisfying the usual functorialities, then it
would inject into
$\cha^{\ast}(\widetilde{\calH}_{N,M},\mathcal{D}_{\pre})$ and the 
arithmetic degrees as well as the values of the height pairing 
would coincide. 
\end{remark}

Mathematisches Institut, Universit\"at zu K\"oln, Weyertal 86-90,
50931 K\"oln, Germany, \texttt{bruinier@math.uni-koeln.de} 

Facultat de Matem\`atiques, Universitat de Barcelona, Gran Via 585,
08007 Barcelona, Spain, \texttt{burgos@mat.ub.es}

Institut f\"ur Mathematik, Humboldt-Universit\"at zu Berlin, Unter den
Linden 6, 10099 Berlin, Germany, \texttt{kuehn@math.hu-berlin.de}

2000 {\em Mathematics subject classification}. 11F41, 14C17, 14G40, 14C20, 11G18

\end{document}